\RequirePackage{fix-cm}

\documentclass[smallextended]{svjour3}       %
\smartqed  
\usepackage{lineno,hyperref}
\usepackage{graphicx}
\usepackage{hyperref}
\usepackage{moreverb}
\usepackage{graphicx,wrapfig}
\usepackage{subfigure}
\usepackage{algorithm}
\usepackage{algpseudocode}
\usepackage{color}
\usepackage{amsmath}
\usepackage{amssymb}
\usepackage{bm}
\usepackage{afterpage}
\usepackage{mathrsfs}
\usepackage{tikz}
\usepackage{pstricks}
\usepackage{diagbox}
\usepackage{arydshln}
\usepackage{multirow}
\usepackage{appendix}
\usepackage{cancel}

\renewcommand{\bm}{\boldsymbol}

\newcommand{\dif}{\mathop{}\!\mathrm{d}}
\newcommand{\btheta}{\bm{\theta}}

\newcommand{\bx}{\bm{x}}
\newcommand{\bk}{\bm{k}}
\newcommand{\bz}{\bm{z}}

\newcommand{\bepsilon}{\bm{\epsilon}}
\newcommand{\bxi}{\bm{\xi}}
\newcommand{\bmu}{\bm{\mu}}

\newcommand{\bsigma}{\bm{\sigma}}

\newcommand{\blambda}{\mathbf{\lambda}}

\newcommand{\bu}{\textbf{u}}
\newcommand{\bd}{\bm{d}_{\text{obs}}}

\newcommand\bphi{\boldsymbol{\phi}}

\newcommand\inth{_{\btheta}}

\let\oldequation\equation
\let\oldendequation\endequation

\renewenvironment{equation}{\linenomathNonumbers\oldequation}{\oldendequation\endlinenomath}

\begin{document}
\title{VI-DGP: A variational inference method with  deep generative prior for solving high-dimensional inverse problems
}


\author{Yingzhi Xia         \and
        Qifeng Liao       \and
        Jinglai Li
}


\institute{Yingzhi Xia \at
              School of Information Science and Technology, ShanghaiTech University, Shanghai 201210, China\\
              Institute of High Performance Computing (IHPC), Agency for Science, Technology and Research (A*STAR), 1 Fusionopolis Way, \#16-16 Connexis, Singapore 138632, Republic of Singapore\\
              \email{Xia\_Yingzhi@ihpc.a-star.edu.sg}           
           \and
           Qifeng Liao \at
           School of Information Science and Technology, ShanghaiTech University, Shanghai 201210, China\\
              \email{liaoqf@shanghaitech.edu.cn}
              \and
           Jinglai Li \at
           School of Mathematics, University of Birmingham, Birmingham B15 2TT, UK\\
              \email{j.li.10@bham.ac.uk}
}

\date{Received: date / Accepted: date}

\maketitle

\begin{abstract}
Solving high-dimensional Bayesian inverse problems (BIPs) with the variational inference (VI) method is promising but still challenging. The main difficulties arise from two aspects. First, VI methods approximate the posterior distribution using a simple and analytic variational distribution, which makes it difficult to estimate complex spatially-varying parameters in practice. Second, VI methods typically rely on gradient-based optimization, which can be computationally expensive or intractable when applied to BIPs involving partial differential equations (PDEs). To address these challenges, we propose a novel approximation method for estimating the high-dimensional posterior distribution. This approach leverages a deep generative model to learn a prior model capable of generating spatially-varying parameters. This enables posterior approximation over the latent variable instead of the complex parameters, thus improving estimation accuracy. Moreover, to accelerate gradient computation, we employ a differentiable physics-constrained surrogate model to replace the adjoint method. The proposed method can be fully implemented in an automatic differentiation manner. Numerical examples demonstrate two types of log-permeability estimation for flow in heterogeneous media. The results show the validity, accuracy, and high efficiency of the proposed method.

\end{abstract}

\section{Introduction}
\label{section_intro}
Inverse problems have extensive applications in science and engineering. Their goal is to determine unknown parameters using indirect and noisy observations. Solving such a challenging problem is a fundamental study in medical imaging, remote sensing, geophysics, and other fields. Identifying parameters from limited observations often involves solving an ill-posed problem that cannot guarantee its stability and uniqueness. To alleviate this problem, many deterministic algorithms solve a penalized least-squares problem with various regularization methods~\cite{engl1996regularization,zhdanov2002geophysical}. Bayesian statistics~\cite{stuart2010inverse,tarantola2005inverse,kaipio2006statistical} provides a framework for inverse problems by treating unknown parameters as random variables and solving them using Bayes' rule. The assigned prior distribution provides a suitable regularization. The estimated posterior distribution determines reasonable solutions along with their uncertainty.

Without a closed-form expression, two types of approximate methods are applied to estimate the posterior distribution in previous studies. The Markov Chain Monte Carlo (MCMC)~\cite{metropolis1953equation,robert1999monte} plays a predominant role, as it is asymptotically exact and easy to implement. The VI method~\cite{blei2017variational,zhang2018advances} estimates the posterior distribution by exploring an optimal approximation within a defined variational distribution family~\cite{barajas2019approximate,povala2022variational,chen2021stein}. However, several common problems still need to be resolved for PDE-constrained BIPs. First, the computationally intensive forward model causes an enormous computational burden. To reduce the computational cost, many surrogates~\cite{mo2019deep} or reduced model-based methods~\cite{chen2021stein,cui2015data} are employed in large-scale problems. Second, a favorable prior distribution should represent all available prior information, but previous methods such as principal component analysis (PCA) and its variants~\cite{liao2019adaptive} have strong assumptions and low accuracy for realistic parameters in prior modeling. Lastly, the curse of dimensionality results in slow convergence and poor approximation, particularly for MCMC, even with the use of advanced methods like sequential Monte Calo~\cite{wan2011bayesian}, Hamiltonian Monte Carlo~\cite{bui2014solving}, and stochastic Newton MCMC~\cite{martin2012stochastic}. 

As an alternative to MCMC, VI methods are widely used in probabilistic machine learning due to their efficiency, flexibility, and scalability, especially in large data scenarios~\cite{blei2017variational}. For high-dimensional BIPs, VI methods can achieve fast convergence and efficient inference by utilizing stochastic gradient-based optimization. Some VI methods, such as mean-field approximation~\cite{jia2021variational,guha2015variational}, can also overcome the curse of dimensionality by assuming independence between different dimensions. However, there are still some bottlenecks when using VI methods to solve PDE-constrained inverse problems with complex parameters. The use of simple and analytical variational distributions, such as multivariate Gaussian~\cite{barajas2019approximate,yang2017bayesian} or Gaussian mixtures~\cite{tsilifis2016computationally}, limits the capability of previous studies to handle complex parameter estimation. Furthermore, the required gradient computation for most VI methods makes them less appealing for solving PDE-constrained inverse problems.  Once these challenges can be addressed or alleviated, the promising VI methods will have broader applications in inverse problems.

In recent years, deep generative models (DGM) have attracted much attention for inverse modeling in various disciplines, such as image processing~\cite{jalal2021robust}, geophysics~\cite{mo2020integration}, compressed sensing~\cite{bora2017compressed}, and material design~\cite{wang2020deep}.  As a data-driven model for prior information representation, DGMs are much more flexible and aim to capture the underlying structure of the given data, enabling the generation of new samples from a learned low-dimensional latent space. The obtained low-dimensional latent variable can serve as the target variable in posterior inference, resulting in dimension reduction. Unlike conventional parameterization methods, non-Gaussian parameters can be well-estimated by various DGMs, such as normalizing flows (NF)~\cite{padmanabha2021solving}, variational autoencoders (VAE)~\cite{xia2022bayesian,laloy2017inversion}, and generative adversarial networks (GANs)~\cite{laloy2018training,patel2022solution}. However, even with advanced strategies, such as domain decomposition~\cite{zhihang2023domain} and multiscale representation and inference~\cite{xia2022bayesian}, sampling methods with DGMs remain computationally intensive. 

The VI approximation is typically solved using gradient-based or Hessian-based optimization~\cite{barajas2019approximate}. In PDE-constrained optimization, the adjoint method is commonly adopted for gradient computation~\cite{warner2015stochastic,wang2018randomized}, but it can be expensive or hard to derive for large and complex physical systems. In contrast, the neural network model is differentiable and can be a potential alternative to the adjoint method~\cite{lye2021iterative,YAN2021114087,wang2021fast}. Unlike conventional surrogates, such as Gaussian process regression~\cite{bilionis2013multi} and polynomial chaos expansion~\cite{xiu2003modeling,marzouk2007stochastic}, the neural network surrogate model can provide a good approximation for high-dimensional parametric PDEs~\cite{zhu2018bayesian,tripathy2018deep}. Recently, physics-informed neural networks (PINN)~\cite{raissi2019physics} have been widely investigated for solving PDEs. The developed physics-constrained surrogates~\cite{zhu2019physics,lu2021learning} can learn the mapping from parameter space to the solution space without using simulation data. Successful applications include Darcy flow~\cite{zhu2019physics}, fluid flows~\cite{sun2020surrogate}, and the Kuramoto-Sivashinsky equation~\cite{geneva2020modeling}. A well-trained neural network surrogate can be used for efficient gradient approximation in PDE-constrained optimization problems, rather than relying solely on forward computation as in other problems.

In this work, we focus on solving high-dimensional inverse problems using VI methods. The main contributions are summarized as follows. First, we propose using deep generative prior (DGP) as the prior model for the VI method, which offers several advantages. As a data-driven model, DGP can encode all prior information from the training data without many assumptions or restrictions. This enables a more informative prior that embodies the underlying complex prior distribution. The VI-DGP method implements posterior estimation for the low-dimensional latent variable, which is more efficient to optimize. By bypassing direct posterior estimation for target parameters, our method leverages the capacity of DGP to mitigate limitations of the variational distribution. Second, we introduce physics-constrained neural networks to address expensive or intractable gradient computation for optimization involving PDEs. We also show that the gradients obtained from neural networks can be effectively applied to stochastic optimization, which can substantially improve the efficiency of the VI-DGP methods. Third, using the asymptotically exact MCMC method as the benchmark, we demonstrate the effectiveness of the proposed method in estimating two types of complex permeability in porous media flow. It should be noted that related complex parameter estimation is very challenging for previous VI methods. With only given three essential components, i.e., prior information (historical data), forward/physical model (PDEs formulation), and noisy observations, the proposed method allows for constructing a complete automatic differentiation workflow that solves high-dimensional BIPs in an efficient (within thousands of iterations) and effective manner. Moreover, the VI-DGP method is easy to implement using existing frameworks like Pytorch and TensorFlow, and can also be applied to other Bayesian inference problems.

The rest of the paper is organized as follows. Section~\ref{sec:Definition} introduces the problem definition of BIPs and explains the difficulties of high-dimensional inverse problems governed by PDEs. Section~\ref{sec:Methodology} gives the proposed  methodology for solving the parameter estimation problems, which includes the DGP for prior information representation in Section~\ref{sec:DGP}, the VI-DGP model for Bayesian inference in Section~\ref{sec:VI}, and the physics-constrained neural networks for the gradient approximation in Section~\ref{sec:PCS}. Section~\ref{sec:Gaussian_Examples} and Section~\ref{sec:Channel_Examples} illustrate two examples of log-permeability estimation in the context of flow in heterogeneous media. Finally, some concluding remarks are provided in Section~\ref{sec:Conclusions}.
\section{Problem setup}
\label{sec:Definition}
Let $\mathcal{D}$ denote a defined spatial domain (in $\mathbb{R}^2$ or $\mathbb{R}^3$), which is bounded,
connected, and with a polygonal boundary $\partial\mathcal{D}$, and $\bx \in \mathcal{D}$ denote a spatial variable. 
In this work, we consider the forward problem governed by physical laws.  Such a physical system can be formulated as PDEs over the spatial domain $\mathcal{D}$ and boundary conditions on the boundary $\partial\mathcal{D}$, e.g., 
\begin{equation}
\label{eq:pde}
\begin{aligned}
\mathcal{N}(\bx, \bk, \bu(\bx,\bk)) &=f(\bx) & & \bx \in \mathcal{D}, \\
\mathfrak{b}(\bx, \bk, \bu(\bx,\bk)) &= g(\bx) & & \bx \in \partial\mathcal{D},
\end{aligned}
\end{equation}
where $\mathcal{N}$ is the partial differential operator and $\mathfrak{b}$ is a boundary operator. $f(\bx)$ denotes the source function, and $g(\bx)$ is the given boundary conditions. Typically, $\bk$ is the spatially-varying parameter (e.g., material property) appearing in the constitutive equations. It can be written as a function $\bk(\bx)$ with respect to the spatial variable $\bx$. $\bu(\bx,\bk)$ is the output or response variable of the physical system.

\par
\subsection{Bayesian inverse problems}
\label{sec:BIPs}
We consider that there is a forward model $\mathcal{F}$ concerning the physical system in Eq.~\eqref{eq:pde}. It maps the unknown parameter $\bk \in \mathbb{R}^{M}$ to the observable output $\bd \in \mathbb{R}^{D}$:
\begin{equation}
\label{eq:forward}
    \bd = \mathcal{F}(\bk) + \bm{\xi},
\end{equation}
where $\bm{\xi} \in \mathbb{R}^{D}$ is the measurement noise. In inverse problems, our interest is to recover the unknown parameter $\bk(\bx)$ from some noisy observations $\bd$. This problem is highly ill-posed, as $D \ll M$,  indicating that the exact parameter is not unique, and its solution is highly sensitive to the measurement noise. To this end, the Bayesian paradigm~\cite{stuart2010inverse} is introduced to highlight the uncertainty of the target parameter $\bk$. One can encode the prior information as a prior distribution $\pi(\bk)$ for the random variable $\bk$. Then, the solution of inverse problems is the posterior distribution with respect to parameter $\bk$ rather than a point estimate. Given the observation data $\bd$, one can calculate the posterior probability $\pi(\bk | \bd)$ via  Bayes' theorem: 
\begin{equation}
\label{eq: Bayes rule}
    \pi(\bk| \bd)=\frac{\pi(\bd | \bk) \pi(\bk)}{\int \pi(\bd | \bk) \pi(\bk) \mathrm{d} \bk},
\end{equation}
where $\pi(\bd | \bk)$ is the likelihood function, which can measure
the discrepancy between the forward predictions and observations, its formulation depends on the type of measurement noise $\bm{\xi}$, i.e., $\pi(\bd | \bk) = \pi_{\bm{\xi}}(\bd-\mathcal{F}(\bk))$. Throughout this work, we assume that $\bm{\xi}$ is a Gaussian distribution with zero mean and diagonal covariance matrix $\hat{\Sigma}$, i.e.,  $\bm{\xi} \sim \mathcal{N}\left(\mathbf{0}, \hat{\Sigma} \right)$. $\hat{\Sigma}$ can define the noise level of each observation. As the parameter $\bk$ is high-dimensional, the denominator in Eq.~\eqref{eq: Bayes rule}, called evidence, involves an intractable high-dimensional integral. Thus, in approximate inference, the evidence in Eq.~\eqref{eq: Bayes rule} mainly serves as a normalization. Then we have
\begin{equation}
    \pi(\bk | \bd) \propto \pi(\bd | \bk) \pi(\bk).
    \label{eq:k_posterior}
\end{equation}

As discussed, two main approaches are applied for posterior approximation. We focus on the VI methods in this work. In practice, the prior information on the spatially-varying parameter $\bk$ is difficult to cast as an analytical distribution. Also, directly introducing an analytical and simple variational distribution to approximate the complex target distribution $\pi(\bk | \bd)$ is not reasonable. In inverse modeling, the parameterization of complex data (e.g., non-Gaussian) is troublesome. However, data collection from historical experiments or prior knowledge is available. To this end,  we consider the data-driven method for modeling the prior information in BIPs. The assumption is that, given training data of $\bk$  before any observation and inference, one can obtain a learned generative model $\bk = \mathcal{G}_{\btheta^{\star}}(\bz)$ and a simple distribution $\pi(\bz)$, where  $\btheta^{\star}$ is parameters of the learned generative model $\mathcal{G}_{\btheta^{\star}}(\cdot)$, and $\bz\in \mathbb{R}^{h}$. Typically, the relatively low-dimensional latent variable $\bz$ can realize dimension reduction since we have  $h\ll M$. One can sample different latent variables from the latent space, and then generate various spatially-varying parameters $\bk$ using the sampled latent variable $\bz$ and pre-trained generative model $\mathcal{G}_{\btheta^{\star}}(\bz)$ correspondingly. With a well-trained generative model, the generated $\bk$ can be considered as samples from the underlying prior distribution $\pi(\bk)$. The inference in Eq.~\eqref{eq:k_posterior} becomes the problem of evaluating the posterior of latent variable $\bz$. Let us write it as
\begin{equation}
\label{eq:PosteriorOfz}
    \pi(\bz | \bd) \propto  \pi(\bd | \bz) \pi(\bz).
\end{equation}
The prior $\pi(\bz)$ is normally a simple distribution (e.g., Gaussian) that we can easily sample and has its closed-form expression. Computing the likelihood function involves $\mathcal{G}_{\btheta^{\star}}(\bz)$ and the forward model $ \mathcal{F}(\bk)$, i.e., $\pi(\bd | \bz)= \pi_{\bm{\xi}}(\bd - \mathcal{F}(\mathcal{G}_{\btheta^{\star}}(\bz)))$, which contains mappings from the latent variable $\bz$ to the spatially-varying parameter $\bk$, and from $\bk$ to observable predictions, respectively. Once $\mathcal{G}_{\btheta^{\star}}(\bz)$ is determined, we are interested in approximating the posterior $\pi(\bz | \bd)$ with the VI method. This method can leverage the capabilities of neural networks for real data generation and gradient computation. Given the analytical prior distribution $\pi(\bz)$, even though we adopt a simple distribution to approximate the posterior distribution $\pi(\bz | \bd)$ in Eq.~\eqref{eq:PosteriorOfz}, the estimated posterior distribution $\pi(\bk | \bd)$ can be very complex due to the representation capacity of the deep generative model. Using the estimated $\pi(\bz | \bd)$ and learned $\mathcal{G}_{\btheta^{\star}}(\bz)$, we can recover samples of the estimated posterior distribution $\pi(\bk | \bd)$. 

We noticed that the VI methods typically define distribution approximation as an optimization problem of statistical distance, e.g., minimizing the Kullback–Leibler (KL) divergence. The defined optimization problems can be solved using standard algorithms like gradient-based methods. However, since the forward model involves PDEs, the VI methods for solving PDE-constrained inverse problems are still restricted by gradient computation. To expedite the VI methods without compromising accuracy, it is worth exploring fast gradient approximation. Importantly, gradient approximation becomes necessary when the gradient is unavailable or computationally expensive for a complex system.

\section{Methodology}
\label{sec:Methodology}

\subsection{Deep generative prior (DGP)}\label{sec:DGP}
In order to bypass direct modeling of the spatially-varying parameter $\bk$, we introduce the DGM for prior modeling. In this sense, the prior information of $\bk$ is cast as a generative model $\mathcal{G}_{\btheta^{\star}}(\bz)$ and a prior distribution $\pi(\bz)$. The DGM aims to learn the underlying distribution from independent and identically distributed (i.i.d) samples. Several popular approaches can accomplish this task, including NF~\cite{rezende2015variational}, VAE~\cite{kingma2013auto}, GAN~\cite{goodfellow2014generative}, etc. NF is constructed using a sequence of invertible transformations. However, its identical-dimensional latent variable cannot favor dimension reduction for the original parameter, leading to complicated inference and expensive computation for BIPs. GAN is notorious for less diversity in generation and unstable training due to its adversarial training nature. In this paper, we adopt VAE to learn a DGP for BIPs because of its desired probabilistic formulation and stable training process.

Given training dataset $ \mathbf{K} = \{ \bk^{(i)}\}_{i=1}^{N}$, where each data point is drawn from the underlying prior distribution $\pi(\bk)$ in Eq.~\eqref{eq: Bayes rule}, $\pi(\bk)$ is the target distribution of the DGP. VAE is a latent variable model that adopts the variational inference method for optimizing model parameters. Introducing the  low-dimensional latent variable $\bz$, the joint distribution $p_{\btheta}(\bk, \bz)$ is factorized as $p_{\btheta}(\bk|\bz) p_{\btheta}(\bz)$, where $p_{\btheta}(\bk|\bz)$ is a probabilistic decoder, $p_{\btheta}(\bz)$ denotes prior distribution of latent variable $\bz$, and $\btheta$ is model parameters. VAE is similar to other likelihood-based models. Its objective is to learn the underlying distribution directly by maximizing the marginal likelihood $p_{\btheta}(\bk)$ of training data.  The direct optimization of $p_{\btheta}(\bk)= \int p_{\btheta}(\bk|\bz) p_{\btheta}(\bz) \mathrm{d} \bz$ involves an intractable integral over $\bz$.  Using Bayes’ rule, the  marginal likelihood can be written as 
\begin{equation}
p_{\btheta}(\bk)=\frac{p_{\btheta}(\bk, \bz)}{p_{\btheta}(\bz|\bk)}=\frac{p_{\btheta}(\bk|\bz) p_{\btheta}(\bz)}{p_{\btheta}(\bz|\bk)},
\label{eqn:vae_likelihood}
\end{equation}
where $p_{\btheta}(\bk|\bz)$ is a probabilistic decoder, and $p_{\btheta}(\bz)$ denotes the prior distribution of the latent variable. $p_{\btheta}(\bz|\bk)$ is the posterior distribution of the latent variable. The computation of $p_{\btheta}(\bz|\bk)$ is also intractable, and thus one can introduce a variational distribution $q_{\bphi}(\bz|\bx)$ to approximate $p_{\btheta}(\bz|\bk)$, where $\bphi$ denotes the encoder model parameters. For any given $q_{\bphi}(\bz|\bk)$, we have
\begin{equation}
\begin{aligned}
\log p_{\btheta}(\bk) &=\mathbb{E}_{q_{\bphi}(\bz|\bk)}\left[\log p_{\btheta}(\bk)\right] \\
&=\mathbb{E}_{q_{\bphi}(\bz|\bk)}\left[\log \left[\frac{p_{\btheta}(\bk, \bz)}{p_{\btheta}(\bz|\bk)}\right]\right] \\
&=\mathbb{E}_{q_{\bphi}(\bz|\bk)}\left[\log \left[\frac{p_{\btheta}(\bk, \bz)}{q_{\bphi}(\bz|\bk)} \frac{q_{\bphi}(\bz|\bk)}{p_{\btheta}(\bz|\bk)}\right]\right] \\
&= \underbrace{\mathbb{E}_{q_{\bphi}(\bz|\bk)}\left[\log \left[\frac{p_{\btheta}(\bk, \bz)}{q_{\bphi}(\bz|\bk)}\right]\right]}_{\mathcal{L}({\btheta,\bphi};\bk)}+ \underbrace{\mathbb{E}_{q_{\bphi}(\bz|\bk)}\left[\log \left[\frac{q_{\bphi}(\bz|\bk)}{p_{\btheta}(\bz|\bk)}\right]\right]}_{D_{KL}\left(q_{\bphi}(\bz|\bk) \| p_{\btheta}(\bz|\bk)\right)}.
\label{eqn:vae_likelihood_elbo}
\end{aligned} 
\end{equation}
Note that the second term above is the KL divergence, which is always non-negative. If and only if $q_{\bphi}(\bz|\bk) = p_{\btheta}(\bz|\bk)$, the KL divergence is equal to zero. Due to the non-negativity of the KL divergence, the first term, called the evidence lower bound (ELBO), provides a lower bound for the marginal log-likelihood. We can rewrite it as 
\begin{equation}
\begin{aligned}
\mathcal{L}({\btheta,\bphi};\bk)=\log p_{\btheta}(\bk)-D_{KL}\left(q_{\bphi}(\bz|\bk) \| p_{\btheta}(\bz|\bk)\right).
\label{eqn:elbo}
\end{aligned} 
\end{equation}
Maximizing $\mathcal{L}({\btheta,\bphi};\bk)$ will maximize the marginal log-likelihood and also make approximation $q_{\bphi}(\bz|\bk)$ close to the true posterior $p_{\btheta}(\bz|\bk)$. So we can maximize $\mathcal{L}({\btheta,\bphi};\bk)$ rather than the marginal log-likelihood for computational convenience~\cite{kingma2013auto}. For the given training dataset $\mathbf{K}$, we can write the ELBO  for any given $\bk^{(i)}$ as 
\begin{equation}
\begin{aligned}
	 \mathcal{L}({\btheta,\bphi};\bk^{(i)}) =& \mathbb E_{q_{\bphi}(\bz|\bk^{(i)})} [\log p\inth(\bk^{(i)},\bz) - \log q_{\bphi}(\bz|\bk^{(i)})] \\
	=& \mathbb E_{q_{\bphi}(\bz|\bk^{(i)})} [\log p\inth(\bk^{(i)}|\bz) ] -  D_{KL}\left(q_{\bphi}(\bz|\bk^{(i)}) || p\inth(\bz)\right).
	\label{eqn:vae_decomposedlowerbound}
\end{aligned}
\end{equation}
Obviously, these two terms play different roles in optimization. The first term is the expected log-likelihood $\log p\inth(\bk|\bz)$, where $\bz$ is sampled from the probabilistic encoder $q_{\bphi}(\bz|\bk)$. Maximizing this term enforces $p\inth(\bk|\bz)$ to assign most of the probability density close to the original $\bk$. The second term aims to minimize the KL divergence between $q_{\bphi}(\bz|\bk)$ and $p\inth(\bz)$, which regularizes the probabilistic encoder $q_{\bphi}(\bz|\bk^{(i)})$ to resemble the prior distribution $p\inth(\bz)$. These two terms are the reconstruction term and the regularization term, respectively.

We still need to specify the distributions for $p\inth(\bk|\bz)$, $q_{\bphi}(\bz|\bk)$, and $p\inth(\bz)$ for computation. Typically, one can assign a simple isotropic Gaussian distribution as the prior distribution, e.g.,
\begin{equation}
    \begin{aligned}
p\inth(\bz)=\mathcal{N}\left(\bz ; \mathbf{0}, \boldsymbol{I}\right).
\label{eqn:prior_z}
\end{aligned}
\end{equation}
Ideally, an appropriate probabilistic encoder $q_{\bphi}(\bz|\bk)$ should be able to approximate the target distribution $p\inth(\bz)$ well. Additionally, a Gaussian distribution with a diagonal covariance can be selected as the variational distribution: 
\begin{equation}
    \begin{aligned}
q_{\bphi}(\bz|\bk) = \mathcal{N}\left(\bz ; \bmu_{\phi}(\bk), \operatorname{diag}(\bsigma_{\bphi}(\bk)^2)\right),  
\label{eqn:encoder_gau}
\end{aligned}
\end{equation}
where $\bmu_{\phi}(\bk)$ and $\bsigma_{\phi}(\bk)$ are computed by the encoder neural networks.  The KL divergence term in Eq.~\eqref{eqn:vae_decomposedlowerbound} has an analytic form~\cite{kingma2013auto} since both $p\inth(\bz)$ and $q_{\bphi}(\bz|\bk)$ are the factorized Gaussian distribution. The distribution $p\inth(\bk|\bz)$ usually depends on the training data. In this paper, we select the Gaussian distribution $\mathcal{N}\left(\bk ; \mathcal{G}_{\btheta}(\bz), \boldsymbol{I}\right)$ for the probabilistic decoder, where $\mathcal{G}_{\btheta}(\bz)$ is the output of the decoder neural networks. The stochastic gradient-based method is applied for large-scale training data to realize the joint optimization for $\{\btheta, \bphi\}$ using the objective function $\mathcal{L}({\btheta,\bphi};\bk)$. The reconstruction term in Eq.~\eqref{eqn:vae_decomposedlowerbound} involves the expectation computation,  which is tackled by Monte Carlo estimation. The gradient $\nabla_{\btheta}\mathbb E_{q_{\bphi}(\bz|\bk)} [\log p\inth(\bk|\bz)]$ can be estimated directly, where the latent variable $\bz$ is randomly sampled from $q_{\bphi}(\bz|\bk)$ for expectation approximation. However, the gradient $\nabla_{\bphi}\mathbb E_{q_{\bphi}(\bz|\bk)} [\log p\inth(\bk|\bz)]$ is difficult to obtain. One cannot swap the gradient and the expectation since the expectation with respect to the distribution $q_{\bphi}(\bz|\bk)$ is a function of $\bphi$. The score function estimator~\cite{blei2017variational,ranganath2014black} can be applied for gradient estimation, but its high variance leads to a slow optimization process. An alternative  differentiable  estimator with low variance is the reparameterization trick~\cite{kingma2013auto,rezende2015variational}, where the latent variable $\bz$ is represented by a deterministic transformation $\bz = \bm{g}_{\bphi} (\bepsilon;\bk)$. 

\begin{figure}[h!]
    \centering
    \includegraphics[width=4in]{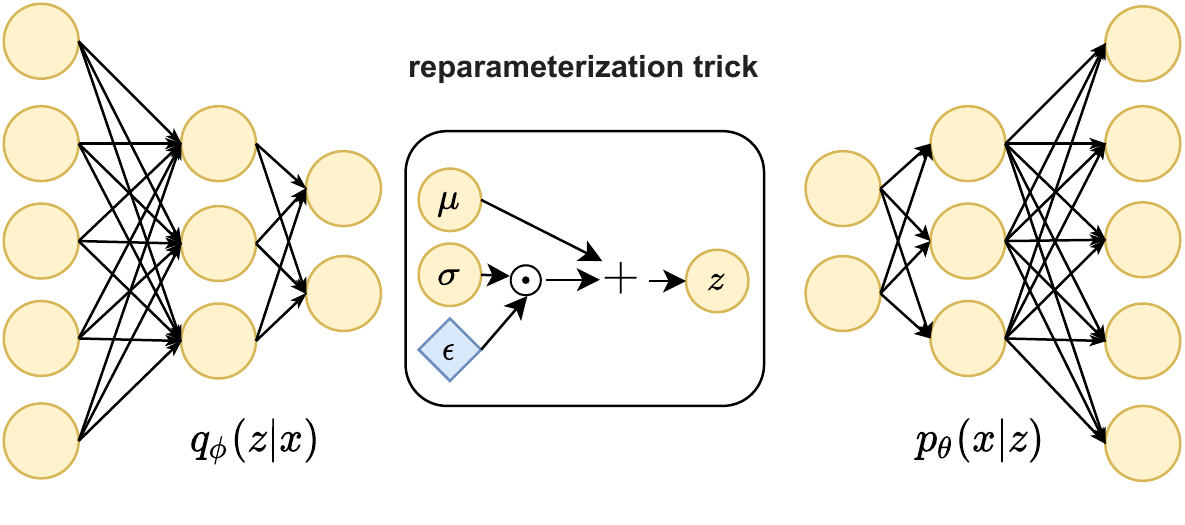}
    \caption{The schematic illustration of the VAE model and the reparameterization trick. The spatially-varying parameter $\bk$ is mapped to a latent variable $\bz$ by the probabilistic encoder $q_{\bphi}(\bz|\bk)$. In turn, the latent variable $\bz$ is mapped to the parameter $\bk$  by the probabilistic decoder $p\inth(\bk|\bz)$}. 
    \label{fig:VAE}
\end{figure}

The differentiable transformation $\bm{g}_{\bphi} (\bepsilon;\bk)$ maps the  auxiliary random noise to the Gaussian distribution in Eq.~\eqref{eqn:encoder_gau} by the following procedure:
\begin{equation}
    \begin{aligned}
    \bz \sim q_{\bphi}(\bz\mid\bk) \quad \Leftrightarrow \quad \bm{g}_{\bphi} (\bepsilon;\bk)=\bmu_{\bphi}(\bk)+\bsigma_{\bphi}(\bk) \odot \bepsilon, \quad \bepsilon \sim \pi(\bepsilon),
\label{eqn:reparameterazation}
\end{aligned}
\end{equation}
where $\odot$  denotes the element-wise product and $\pi(\bepsilon) = \mathcal{N}(\mathbf{0}, \mathbf{I})$. Then the random variable $\bz$ only depends on two deterministic outputs of the encoder neural networks by introducing an auxiliary random variable $\bepsilon$. Since the operators $+$ and $\odot$ are differentiable, the gradient $\nabla_{\bphi}\mathbb E_{q_{\bphi}(\bz|\bk)} [\log p\inth(\bk|\bz)]$ is available. It can be written as
\begin{equation}
    \begin{aligned}
    \nabla_{\bphi}\mathbb E_{q_{\bphi}(\bz|\bk)} [\log p\inth(\bk|\bz)] =&  \mathbb E_{\pi(\bepsilon)}(\nabla_{\bphi}\log p\inth(\bk| \bz)))\\
    =&\mathbb E_{\pi(\bepsilon)}\left[ \frac{\partial \log p\inth(\bk| \bz)}{\partial \bz} \frac{ \partial \bm{g}_{\bphi} (\bepsilon;\bk)}{\partial \bphi}\right]_{\bz =\bm{g}_{\bphi} (\bepsilon;\bk)},
\label{eqn:gradient_theta}
\end{aligned}
\end{equation}
which can be directly estimated by the Monte Carlo method with $L$ samples drawn from $\pi(\bepsilon)$. Then the ELBO in Eq.~\eqref{eqn:vae_decomposedlowerbound} can be rewritten as
\begin{equation}
\begin{aligned}
   \mathcal{L}({\btheta,\bphi};\bk^{(i)}) = \frac{1}{L}\sum_{l=1}^{L} \log p_{\btheta} (\bk^{(i)}|\bz^{(i,l)})  -  D_{KL}\left(q_{\bphi}(\bz|\bk^{(i)}) || p\inth(\bz)\right),
\label{eqn:elbo_loss}
\end{aligned}
\end{equation}
where $\bz^{(i,l)}$ is the $l-$th sample drawn from $q_{\bphi}(\bz|\bk^{(i)})$. To improve computational efficiency, the training of neural networks usually adopts the minibatch stochastic gradient-based method, where the training dataset is divided into many subsets. Each subset contains $n$ data points for each iteration. The optimization objective function in each iteration can be written as
\begin{equation}
\begin{aligned}
   \tilde{\mathcal{L}}({\btheta,\bphi};\bk^n) = \frac{1}{n}\sum_{i=1}^{n} \mathcal{L}({\btheta,\bphi};\bk^{(i)}).
\label{eqn:opt_objective}
\end{aligned}
\end{equation}
One can apply the stochastic gradient-based method, such as Adam~\cite{kingma2014adam}, to optimize the probabilistic encoder $q_{\bphi}(\bz|\bk)$ and the probabilistic decoder $p_{\btheta}(\bk|\bz)$ using the above objective function. Fig.~\ref{fig:VAE} depicts a schematic illustration of the VAE model and the reparameterization trick. The training procedure is outlined in Algorithm~\ref{alg:DGM}. 

For BIPs, we can obtain training dataset $\{ \bk^{(i)}\}_{i=1}^{N}$ based on the history data or the prior knowledge and use Algorithm~\ref{alg:DGM} to learn the DGP that represents the prior information. The underlying prior distribution $\pi(\bk)$ in Eq.~\eqref{eq:k_posterior} can be approximated as $\pi(\bk) \approx \int p_{\btheta}(\bk|\bz) p_{\btheta}(\bz) \dif \bz$, where $p_{\btheta}(\bz) = \mathcal{N}\left(\mathbf{0},\boldsymbol{I}\right)$, and $p_{\btheta}(\bk|\bz)$ is the learned probabilistic decoder. The prior distribution $\pi(\bz)$ in Eq.~\eqref{eq:PosteriorOfz} can be defined as $\pi(\bz) = p_{\btheta}(\bz)$, which can be a simple Gaussian distribution. In BIPs, the  process of generating new prior samples from the underlying distribution $\pi(\bk)$ is as follows:
\begin{equation}
\begin{aligned}
   \bk^{\prime} = \mathcal{G}_{\btheta^{\star}}(\bz^{\prime}), \quad \bz^{\prime} \sim \pi(\bz),
\label{eqn:sample_k}
\end{aligned}
\end{equation}
where $\mathcal{G}_{\btheta^{\star}}$ is the learned decoder neural networks, and $\bk^{\prime} $ can be regarded as sample drawn from $\pi(\bk)$. In this way, the prior information can be cast as the DGP, which includes the simple prior distribution $\pi(\bz)$ and the learned generative model $\mathcal{G}_{\btheta^{\star}}(\bz)$. The method takes advantage of the analytical prior distribution $\pi(\bz)$ and generative model $\mathcal{G}_{\btheta^{\star}}(\bz)$, allowing for sampling from a simple distribution while still being able to generate complex real data by exploiting the high representation capacity of neural networks.

\begin{algorithm}[h!]
	\caption{The training of generative prior}
	\label{alg:DGM}
	\begin{algorithmic}[1]
		\Require Prior data $\{\bk\}_{i=1}^{N}$, training epoch $E$, batch size $n$, learning rate $\eta$, $L=1$.
		\State Initialize $ \bphi,  \btheta \leftarrow \text { Initialize the encoder and decoder parameters }$ 
		\For {$i = 1:E$}
		\For {$j = 1:\frac{N}{n}$}
		\State $\bk^n \leftarrow$ Sample minibatch $n$ data points from $\{\bk\}_{i=1}^{N}$
		\State  $\epsilon^n \leftarrow$ Sample noise from Gaussian distribution $\mathcal{N}(0, I )$
		\State  $\bz^n \leftarrow$ Compute by encoder network with Eq.~\eqref{eqn:reparameterazation}
		\State $\nabla_{ \btheta} \tilde{\mathcal{L}}, \nabla_{\bphi} \tilde{\mathcal{L}} \leftarrow$ \text{Calculate gradients of $\tilde{\mathcal{L}}\left( \btheta, \bphi; \bk^{n} \right)$ w.r.t $\btheta$ and $\bphi$} 
		\State $ \btheta =  \btheta + \eta \nabla_{ \btheta} \tilde{\mathcal{L}}$ 
		\State $\bphi = \bphi + \eta \nabla_{\bphi} \tilde{\mathcal{L}} $
		\EndFor
		\EndFor
		\Ensure probabilistic encoder  $q_{\bphi^{\star}}(\bz|\bk)$, probabilistic decoder  $p_{\btheta^{\star}}(\bk|\bz)$. 
		
	\end{algorithmic}
\end{algorithm}

\subsection{Variational inference with  deep generative prior (VI-DGP)}\label{sec:VI}
Suppose we have learned an appropriate generative model $\mathcal{G}_{\btheta^{\star}}$ via Algorithm~\ref{alg:DGM}. As discussed in Section~\ref{sec:BIPs}, the estimation of the posterior distribution $\pi(\bk | \bd)$ in  Eq.~\eqref{eq:k_posterior} can degenerate into evaluating the posterior of the latent variable, i.e., $\pi(\bz | \bd)$. Without an analytical solution, we adopt the variational inference method for the posterior approximation to emphasize computational efficiency. By introducing a variational distribution $\tilde{q}_{\blambda}(\bz)$ parameterized by $\blambda$, we can determine a good approximation for $\pi(\bz | \bd)$ by minimizing the KL divergence. The KL divergence can be written as 
\begin{equation}
\begin{aligned}
D_{KL}\left(\tilde{q}_{\blambda}(\bz)\| \pi(\bz|\bd)\right) &= \int_{\bz} \tilde{q}_{\blambda}(\bz) \log \frac{\tilde{q}_{\blambda}(\bz)}{\pi(\bz|\bd)}\dif \bz\\  
&= \mathbb{E}_{\tilde{q}_{\blambda}} (\log \tilde{q}_{\blambda}(\bz)) - \mathbb{E}_{\tilde{q}_{\blambda}} \left[\log \pi(\bz|\bd)\right] \\
&=\mathbb{E}_{\tilde{q}_{\blambda}} (\log \tilde{q}_{\blambda}(\bz)) - \mathbb{E}_{\tilde{q}_{\blambda}} \left[\log \frac{\pi(\bz,\bd)}{\pi(\bd)}\right] \\
&=  \mathbb{E}_{\tilde{q}_{\blambda}} \left[\log \tilde{q}_{\blambda}(\bz)\right] - \mathbb{E}_{\tilde{q}_{\blambda}} \left[\log \pi(\bz,\bd)\right] + \log \pi(\bd),
\label{eq:kl_vb1}
\end{aligned}
\end{equation}
where $ \log \pi(\bd)$ is a non-negative constant, so we have 
\begin{equation}
\begin{aligned}
D_{KL}\left(\tilde{q}_{\blambda}(\bz)\| \pi(\bz|\bd)\right) \geq \mathbb{E}_{\tilde{q}_{\blambda}} \left[\log \tilde{q}_{\blambda}(\bz)\right] - \mathbb{E}_{\tilde{q}_{\blambda}} \left[\log \pi(\bz,\bd)\right]. 
\label{eq:kl_vb2}
\end{aligned}
\end{equation}
Minimizing the above KL divergence is equivalent to maximizing the following lower bound:
\begin{equation}
\begin{aligned}
\mathcal{L}_{VI} =  \mathbb{E}_{\tilde{q}_{\blambda}} \left[\log \pi(\bz,\bd)\right] + \mathbb{H}\left[\tilde{q}_{\blambda}(\bz)\right],
\label{eq:kl_vb3}
\end{aligned}
\end{equation}
where $\mathbb{H}\left[\tilde{q}_{\blambda}(\bz)\right] = -\mathbb{E}_{\tilde{q}_{\blambda}} \left[\log \tilde{q}_{\blambda}(\bz)\right]$ is the entropy. One can adopt the Monte Carlo method to approximate expectations. Nevertheless, the gradient $\nabla \mathcal{L}_{VI}$ is also intractable for stochastic gradient-based optimization due to the non-differentiable operator. In this paper, we assume that the variational distribution $\tilde{q}_{\blambda}(\bz)$ is a Gaussian distribution with a diagonal covariance, i.e., $\tilde{q}_{\blambda}(\bz) = \mathcal{N}\left(\bz; \tilde{\bmu}, \operatorname{diag}(\tilde{\bsigma}^2)\right)$, where $\blambda:= \{\tilde{\bmu},\tilde{\bsigma}\}$ is the parameter to be estimated. Such an assumption is reasonable since the latent variable in DGP is the Gaussian distribution with diagonal covariance, while we can still estimate the complex distribution $\pi(\bk)$ with DGP. Similarly, it can be seen that we can approximate the complex posterior $\pi(\bk | \bd)$ using estimated $\tilde{q}_{\blambda}(\bz)$ and the learned generative model. We also need to employ the reparameterization trick to handle the intractable gradient. The VAE proposes the reparameterization trick as an alternative estimator that can resolve the same issue in our posterior estimation problem. This is why we favor VAE for DGP  modeling. Using the transformation $\bm{g}_{\blambda} (\bepsilon) = \tilde{\bmu} + \tilde{\bsigma} \odot \bepsilon, \bepsilon \sim \pi(\bepsilon)$ like Eq.~\eqref{eqn:reparameterazation}, we can make an approximation with the Monte Carlo method, i.e.,
\begin{equation}
\begin{aligned}
\mathcal{L}_{VI} \approx  \frac{1}{M_s}\sum_{i=1}^{M_s}  \left[ \log \pi( \bm{g}_{\blambda} (\bepsilon_i),\bd) - \tilde{q}_{\blambda}(\bm{g}_{\blambda} (\bepsilon_i)) \right],
\label{eq:kl_vb4}
\end{aligned}
\end{equation}
where $\bepsilon_i$ denotes the $i$-th sample drawn from $\pi(\bepsilon) = \mathcal{N}\left(\mathbf{0},\boldsymbol{I}\right)$, and $M_s$ is the number of samples used for approximation. Based on the approximate lower bound, automatic differentiation can be utilized to compute the gradient $\nabla_{\blambda} \mathcal{L}_{VI}$. Then the parameter $\blambda$ can be optimized using the stochastic gradient-based method. We can also write the gradient $\nabla_{\blambda} \mathcal{L}_{VI}$ explicitly as
\begin{equation}
\begin{aligned}
\nabla_{\blambda} \mathcal{L}_{VI} & =  \nabla_{\blambda}\mathbb{E}_{\tilde{q}_{\blambda}} \left[\log \pi(\bz,\bd)\right] + \nabla_{\blambda}\mathbb{H}\left[\tilde{q}_{\blambda}(\bz)\right]\\
& = \nabla_{\blambda}\mathbb{E}_{\pi(\bepsilon)} \left[\log \pi(\bm{g}_{\blambda} (\bepsilon),\bd)\right] + \nabla_{\blambda}\mathbb{H}\left[\tilde{q}_{\blambda}(\bz)\right]\\
&= \mathbb{E}_{\pi(\bepsilon)} \left[ \nabla_{\bm{g}_{\blambda}}\log \pi(\bm{g}_{\blambda} (\bepsilon),\bd) \nabla_{\blambda }\bm{g}_{\blambda} (\bepsilon)\right] + \nabla_{\blambda }\mathbb{H}\left[\tilde{q}_{\blambda}(\bz)\right].
\label{eq:vb_gradient}
\end{aligned}
\end{equation}
Note that the first term is an expectation, which can be approximated with the Monte Carlo method as
\begin{equation}
\begin{aligned}
\mathbb{E}_{\pi(\bepsilon)} \left[ \nabla_{\bm{g}_{\blambda}}\log \pi(\bm{g}_{\blambda} (\bepsilon),\bd) \nabla_{\blambda }\bm{g}_{\blambda} (\bepsilon)\right] \approx \frac{1}{M_s}\sum_{i=1}^{M_s}  \left[ \nabla_{\bm{g}_{\blambda}}\log \pi(\bm{g}_{\blambda} (\bepsilon_i),\bd) \nabla_{\blambda }\bm{g}_{\blambda} (\bepsilon_i)\right].
\label{eq:vb_gradient_term1}
\end{aligned}
\end{equation}
 With the differentiable transformation $\bm{g}_{\blambda} (\bepsilon)$, the second term in Eq.~\eqref{eq:vb_gradient} can be written as  
\begin{equation}
\begin{aligned}
\nabla_{\blambda} \mathbb{H}\left[\tilde{q}_{\blambda}\right] &=-\nabla_{\blambda} \mathbb{E}_{\pi(\bepsilon)}\left[\log \tilde{q}_{\blambda}\left(\bm{g}_{\blambda}(\bepsilon)\right)\right] \\
&=-\mathbb{E}_{\pi(\bepsilon)}\left[\nabla_{\blambda} \log \tilde{q}_{\blambda}\left(\bm{g}_{\blambda}(\bepsilon)\right)\right] \\
&=-\mathbb{E}_{\pi(\bepsilon)}\left[\nabla_{\bm{g}_{\blambda}} \log \tilde{q}_{\blambda}\left(\bm{g}_{\blambda}(\bepsilon)\right) \nabla_{\blambda} \bm{g}_{\blambda}(\bepsilon)\right].
\label{eq:entropy_gradient}
\end{aligned}
\end{equation}
Since we assume that the variational distribution $\tilde{q}_{\blambda}(\bz)$ is a Gaussian distribution, its normalization constant also depends on the variational parameters $\blambda$. The third line above involves the expectation of the score function, i.e., $\mathbb{E}_{\tilde{q}_{\blambda}} [\nabla_{\blambda} \log \tilde{q}_{\blambda}(\bz)]$. However, the expectation $\mathbb{E}_{\tilde{q}_{\blambda}} [\nabla_{\blambda} \log \tilde{q}_{\blambda}(\bz)]$ is always zero~\cite{ranganath2014black}. Therefore, we can obtain the expectation in the third line directly.
The expectation in Eq.~\eqref{eq:entropy_gradient} can be approximated with the Monte Carlo method as
\begin{equation}
\begin{aligned}
\nabla_{\blambda} \mathbb{H} \left[ \tilde{q}_{\blambda}\right] \approx -\frac{1}{M_s}\sum_{i=1}^{M_s}  \left[ \nabla_{\bm{g}_{\blambda}}\log \tilde{q}_{\blambda}\left(\bm{g}_{\blambda}(\bepsilon_i)\right) \nabla_{\blambda} \bm{g}_{\blambda}(\bepsilon_i) \right],
\label{eq:entropy_term2}
\end{aligned}
\end{equation}
where the sampled noise $\bepsilon_i$ is the same as Eq.~\eqref{eq:vb_gradient_term1} in each optimization iteration. By using $M_s$ random samples, we can write the estimated gradient of the lower bound $\nabla_{\blambda} \mathcal{L}_{VI}$ for stochastic optimization as 
\begin{equation}
\begin{aligned}
\nabla_{\blambda} \mathcal{L}_{VI} & = \frac{1}{M_s}\sum_{i=1}^{M_s} \left[ \nabla_{\bm{g}_{\blambda}}\log \pi(\bm{g}_{\blambda} (\bepsilon_i),\bd) \nabla_{\blambda }\bm{g}_{\blambda} (\bepsilon_i)-\nabla_{\bm{g}_{\blambda}}\log \tilde{q}_{\blambda}\left(\bm{g}_{\blambda}(\bepsilon_i)\right) \nabla_{\blambda}\bm{g}_{\blambda}(\bepsilon_i) \right]\\
&=\frac{1}{M_s}\sum_{i=1}^{M_s} \left[ (  \nabla_{\bm{g}_{\blambda}}\log \pi(\bm{g}_{\blambda} (\bepsilon_i),\bd)  -\nabla_{\bm{g}_{\blambda}}\log \tilde{q}_{\blambda}(\bm{g}_{\blambda}(\bepsilon_i))\nabla_{\blambda}\bm{g}_{\blambda}(\bepsilon_i)\right].
\label{eq:entropy_gradient_final}
\end{aligned}
\end{equation}

With the above gradient, maximizing the lower bound $\mathcal{L}_{VI}$ with stochastic gradient ascent will obtain an appropriate approximation for $\pi(\bz|\bd)$. Note that the first term in Eq.~\eqref{eq:entropy_gradient_final} involves the gradient $\nabla_{\bm{g}_{\blambda}} \log \pi(\bz,\bd)|_{\bz = \bm{g}_{\blambda}(\bepsilon_i)}$, where $\pi(\bz,\bd) =  \pi(\bd | \bz) \pi(\bz)$ and $\pi(\bd | \bz)= \pi_{\bm{\xi}}(\bd - \mathcal{F}(\mathcal{G}_{\btheta^{\star}}(\bz)))$. It is easy to see that the optimization requires two necessary gradient computations, i.e., the gradient $\frac{\partial \mathcal{F}(\bk)}{\partial \bk}$ and the gradient  $\frac{\partial \mathcal{G}_{\btheta^{\star}}(\bz)}{\partial \bz}$. $\frac{\partial \mathcal{G}_{\btheta^{\star}}(\bz)}{\partial \bz}$ is easy to compute by adopting automatic differentiation since $\mathcal{G}_{\btheta^{\star}}$ is constructed by neural networks. Unfortunately, $\frac{\partial \mathcal{F}(\bk)}{\partial \bk}$ is often not available. Note that the forward model in most applications is the physical model involving numerical PDEs. Then the potential difficulties are two-fold: first, the gradient computation associated with the complex physics model is very challenging to obtain; second, even though the adjoint method for some models is available, the computation cost is not affordable if the gradient estimation in Eq.~\eqref{eq:entropy_gradient_final} requires a large $M_s$ to ensure stable optimization. These issues dramatically decrease the advantages and popularity of solving the BIPs with the VI methods.  

In the next section, we will introduce the neural network surrogate for gradient approximation. The neural networks can act as an alternative to the adjoint method. One can construct a complete neural network model and directly ask the automatic differentiation to tackle stochastic gradient-based optimization, which is easy to implement and highly efficient. We also show in the numerical experiments that a small $M_s$ (even $M_s = 1$) can bring a stable optimization process under the reparameterization trick. In this way, the VI-DGP method can guarantee efficiency and solve complex parameter estimation problems by exploiting the representation capability of the deep generative model.
\begin{figure}[h!]
    \centering
    \includegraphics[scale=1]{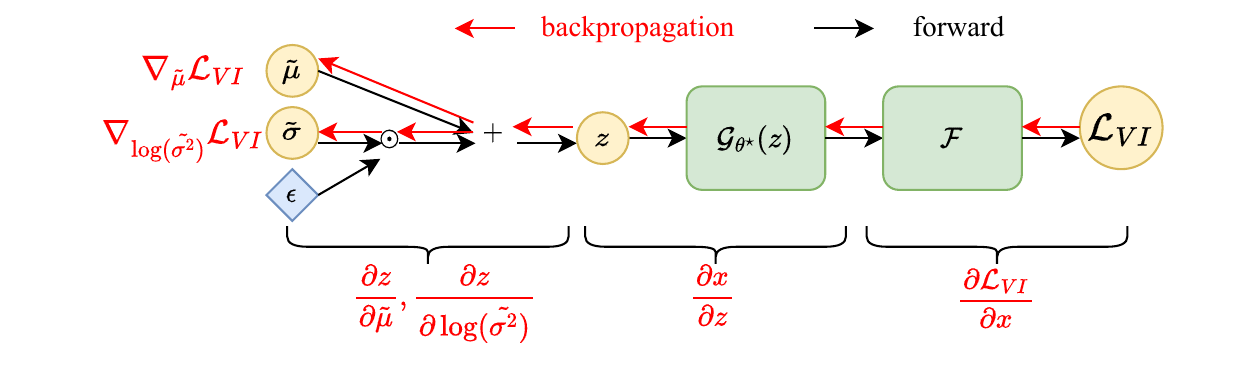}
    \caption{The complete workflow of the VI-DGP for BIPs. The black arrows illustrate the forward computation for the lower bound $\mathcal{L}_{VI}$ in Eq.~\eqref{eq:kl_vb3}. The red arrows indicate the gradient computation with respect to $\tilde{\bmu}$ and $\log(\tilde{\bsigma}^2)$. The differentiable operator for $\bz$ is constructed under the reparameterization trick.} 
    \label{fig:AD_ELBO}
\end{figure}

The demonstration of the forward and backward computation in optimization is given in Fig.~\ref{fig:AD_ELBO}. The forward model $\mathcal{F}$ can be the finite element method solver or neural network surrogate model. The detailed procedure is shown in Algorithm~\ref{alg:VIM}.

\begin{remark}
Directly optimizing $\bsigma$ will cause unstable convergence even if gradient clipping is introduced. To ensure the stable optimization for parameter $\blambda$,  we have to optimize log variance $\log \bsigma^2$ instead of variance or standard deviation $\bsigma$ in implementation. This is also a training trick used in VAE. $\bsigma$ is typically positive and close to $0$. However, poor floating-point arithmetic and unstable gradient computation around $0$ lead to numerical instability. The $\log \bsigma^2$ can transform the narrow feasible domain into a broader space, making stochastic gradient-based optimization more stable and easier to converge.
\end{remark}
\begin{remark}
Directly approximating the entropy term using the Monte Carlo method in Eq.~\eqref{eq:kl_vb4} and applying automatic differentiation techniques to optimize the variational lower bound $\mathcal{L}_{VI}$ may lead to high variance. \cite{roeder2017sticking} illustrates that one can remove the gradient with respect to the variational parameters that correspond to the score function, resulting in an unbiased gradient estimator. The introduced implementation tricks are also applied in our experiment, which makes a stable convergence.
\end{remark}

\begin{algorithm}[h!]
	\caption{The VI-DGP mothod}
	\label{alg:VIM}
	\begin{algorithmic}[1]
	\Require generative model $\mathcal{G}_{\btheta^{\star}}$, forward model $\mathcal{F}$, optimization iteration $N_{opt}$, sampling number $M_{s}$, learning rate $\eta_{\tilde{\bmu}}, \eta_{\tilde{\bsigma}}$, number of posterior samples $N_{s}$.
	
	\State Initialize $ \tilde{\bmu} , \log(\tilde{\bsigma}^2) $
	\While {$i < N_{opt}$}
	\State sample $M_s$ latent variables $\bz^{M_s}$ using the reparameterization trick
	\State   compute $M_s$ spatially-varying parameters $\bx^{M_s}$: $\bx^{M_s} = \mathcal{G}_{\btheta^{\star}}(\bz^{M_s})$ 
	\State compute $M_s$ predictions using the forward model $\mathcal{F}(\bx^{M_s})$    
	\State compute the gradient w.r.t. $\tilde{\bmu}, \log(\tilde{\bsigma}^2)$: $ \nabla_{\tilde{\bmu}}\mathcal{L}_{VI}, \nabla_{\log(\tilde{\bsigma}^2)} \mathcal{L}_{VI}$ 
	\State $ \tilde{\bmu} =  \tilde{\bmu}+ \eta_{\tilde{\bmu}} \nabla_{\tilde{\bmu}}\mathcal{L}_{VI}$ 
	\State $\log(\tilde{\bsigma}^2) = \log(\tilde{\bsigma}^2) + \eta_{\tilde{\bsigma}} \nabla_{\log(\tilde{\bsigma}^2)} \mathcal{L}_{VI}$
	\EndWhile
	
	\State Let $q(\bz;\blambda^{\star}) =  \mathcal{N}(\tilde{\bmu}^{\star},  \operatorname{diag}(\tilde{\bsigma}^{\star 2}))$, $\tilde{\bmu}^{\star}, \tilde{\bsigma}^{\star}$ are obtained parameters.
	\State Sample $N_s$ posterior samples by 
	$$\bx^{(i)} =  \mathcal{G}_{\btheta^{\star}}(\bz^{(i)}), \quad \bz^{(i)} \sim \mathcal{N}(\tilde{\bmu}^{\star}, \operatorname{diag}(\tilde{\bsigma}^{\star 2}))$$
	\Ensure posterior samples $\{\bx^{(i)}\}_{i = 1}^{N_s}$ of BIPs. 

	\end{algorithmic}
\end{algorithm}

\subsection{Gradient approximation with neural networks}
\label{sec:PCS}
As a model for universal function approximation, deep neural networks dominate various high-dimensional tasks~\cite{lu2021learning,zhu2018bayesian,khoo2019switchnet,fan2019solving}. The study of solving PDEs with deep neural networks is promising in science and engineering~\cite{raissi2019physics,li2023deep}. We are interested in training a surrogate model with deep neural networks due to its inherent automatic differentiation~\cite{wang2021fast}. A well-trained neural network surrogate can provide gradient approximations for stochastic gradient descent/ascent in VI methods.

Data-driven and model-driven are two primary methods for training the neural networks for physical models~\cite{zhu2018bayesian,raissi2019physics}. The model-driven method follows the physics model and can incorporate physical constraints into the loss function to learn the surrogate model without the need for simulation data. Its loss function typically includes residual loss regarding PDEs and boundary conditions. If we need to construct a surrogate for the parametric PDEs given in Eq.~\eqref{eq:pde}, we can use $\bu(\bx,\bk,\Theta)$ as the neural networks with parameters $\Theta$ and write the loss function as follows:
\begin{equation}
\begin{aligned}
J(\bu(\bx,\bk ; \Theta) = J_{\text{pde}}(\bu(\bx,\bk ; \Theta))+\gamma J_{\text{b}}(\bu(\bx, \bk ; \Theta)),
\label{eq:pcs_loss}
\end{aligned}
\end{equation}
where $\gamma$ is the hyperparameter in training,  $J_{\text{pde}}(\cdot)$ and $J_{\text{b}}(\cdot)$ denote the residual loss for PDEs and boundary conditions, respectively. Although $\bu(\bx,\bk,\Theta)$ can be the mesh-free model based on random samples $\bx^{(i)}$ in the defined domain~\cite{lu2021learning}, to take advantage of the computational efficiency and fast convergence of convolutional neural networks (CNNs)~\cite{zhu2019physics}, one can choose uniformly distributed collocation points $ \{\bx_{\mathcal{D}}^{(i)}\}_{i=1}^{n_p}$ and $ \{\bx_{\partial\mathcal{D}}^{(i)}\}_{i=1}^{n_b}$ for PDEs loss and boundary loss, respectively. In this paper, since the spatially-varying parameter has been discretized, we can adopt the uniformly distributed points for $\bx$, similar to the finite element or finite difference method. For each iteration, given the training data $\{\bk^{(j)}(\bx)\}_{j=1}^{N_k}$, we can rewrite the two terms in Eq.~\eqref{eq:pcs_loss} as
\begin{equation}
\begin{aligned}
J_{\text{pde}}(\bu(\bx, \bk ; \Theta)) = \frac{1}{n_s n_p}\sum_{j=1}^{n_s}\sum_{i=1}^{n_p}   
\|\mathcal{N}(\bx_{\mathcal{D}}^{(i)},\bk^{(j)}(\bx_{\mathcal{D}}^{(i)}),\bu(\bx_{\mathcal{D}}^{(i)},\bk^{(j)}(\bx_{\mathcal{D}}^{(i)})))-f(\bx_{\mathcal{D}}^{(i)})\|^2,  \\
J_{\text{b}}(\bu(\bx, \bk ; \Theta)) = \frac{1}{n_s n_b}\sum_{j=1}^{n_s}\sum_{i=1}^{n_b} \|\mathfrak{b}(\bx_{\partial\mathcal{D}}^{(i)},\bk^{(j)}(\bx_{\partial\mathcal{D}}^{(i)}),\bu(\bx_{\partial\mathcal{D}}^{(i)},\bk^{(j)}(\bx_{\partial\mathcal{D}}^{(i)})))- g(\bx_{\partial\mathcal{D}}^{(i)})\|^2,
\label{eq:empirical_loss}
\end{aligned}
\end{equation}
where $n_s$ is the batch size of the training data in the training procedure. With the training data $\{\bk^{(j)}(\bx)\}_{j=1}^{N_k}$ and the above discretization form, one can obtain a good approximation $\Theta^{\star}$ by minimizing the loss function in Eq.~\eqref{eq:pcs_loss}, i.e.,
\begin{equation}
\begin{aligned}
\Theta^{\star} = \underset{\Theta}{\arg \min } J(\bu(\bx,\bk ; \Theta)).
\label{eq:argmin loss}
\end{aligned}
\end{equation}

The training procedure for the physics-constrained surrogate model is summarized in Algorithm~\ref{alg:PCS}. Suppose we have obtained a good approximation for the forward model, the $\frac{\partial \mathcal{L}_{VI}}{\partial \bk}$ in Fig.~\ref{fig:AD_ELBO} can be computed by the neural networks with automatic differentiation. It can bypass the expensive computation of the adjoint method.

\begin{algorithm}[h!]
	\caption{The training of physics-constrained surrogate model}
	\label{alg:PCS}
	\begin{algorithmic}[1]
	\Require Dataset $\{\bk^{(j)}(\bx)\}_{j=1}^{N_k}$, neural networks $u(\bx,\bk ; \Theta)$, hyperparameter $\gamma$,  training epoch $E_s$, batch size $n_s$,  learning rate $\eta_{\Theta}$
	
	\State Initialize $ \Theta $ 
	\For {$i = 1:E_s$}
		\For {$j = 1:\frac{N_k}{n_s}$}
	\State $\bk^{n_s} \leftarrow$ Sample minibatch $n_s$ data points from $\{\bk^{(j)}(\bx)\}_{j=1}^{N_k}$
	\State $\nabla_{ \Theta} J(u(\bx,\bk^{n_s} ; \Theta) \leftarrow$ \text{compute gradients of $J(u(\bx,\bk^{n_s}; \Theta))$ w.r.t $\Theta$} 
	\State $ \Theta =  \Theta - \eta_{\Theta} \nabla_{ \Theta} J(u(\bx,\bk^{n_s} ; \Theta)$ 
    \EndFor
    \EndFor
    \Ensure surrogate model $u(\bx,\bk ; \Theta^{\star})$
	\end{algorithmic}
\end{algorithm}

\section{Numerical study}
In this section, we consider the problem of estimating the log-permeability field in a single-phase, steady-state Darcy flow. Given a log-permeability field $\bk$, the pressure field $\bm{p}$ and velocity field $\bm{v}$ are governed by the equations:  
\begin{equation}
\begin{aligned}
\bm{v}(\bx) &=-\exp(\bk(\bx)) \nabla \bm{p}(\bx), \quad \bm{\bx} \in \mathcal{D}, \\
\nabla \cdot \bm{v}(\bx) &= f(\bx), \quad \bx \in \mathcal{D},  
\label{eq:darcy}
\end{aligned}
\end{equation}
\noindent  
with boundary conditions
\begin{equation}
\begin{aligned}
\bm{v}(\bx) \cdot \hat{\bm{n}}&=0, \quad \bx \in \Gamma_{N}, \\
  p(\bx) &= 1, \quad \bx \in \Gamma_{D_{l}},\\
  p(\bx) &= 0, \quad \bx \in \Gamma_{D_{r}},
  \label{eq:darcy_boundary}
\end{aligned}
\end{equation}
where $\mathcal{D}$ denotes a 2D unit square domain $\mathcal{D}=[0,1]^2$, and $\hat{\bm{n}}$ is the unit normal vector to the Neumann boundary $\Gamma_{N}$. The Neumann boundary $\Gamma_{N}$ consists of the top boundary $\Gamma_{D_{t}}$ and bottom boundary $\Gamma_{D_{b}}$, and the Dirichlet boundary consists of the left boundary $\Gamma_{D_{l}}$ and right boundary $\Gamma_{D_{r}}$. We set the source term $f(\bx) = 3$. The spatial domain is discretized into uniform $64 \times 64$ grids. In BIPs, we need to estimate the unknown log-permeability field based on collected noisy observations from the pressure field. Two types of log-permeability field estimation are used to demonstrate the performance of the proposed method. The Gaussian random field (GRF) is a typical example in many previous works~\cite{liao2019adaptive,mo2019deep}. The assumed GRF with fixed mean, covariance, and correlation length can be parameterized by the truncated Karhunen-Lo\`{e}ve expansion (KLE). However, this assumption is unrealistic as spatially-varying parameters typically involve a nontrivial correlation structure. We study the GRF with uncertain correlation length in the numerical example to validate the advantage of DGP in parameter representation. The other example is the complex channelized random field~\cite{wan2011bayesian,laloy2017inversion}, a common geological media in the groundwater flow. Note that the BIPs regarding the non-Gaussian random field still have difficulties in parameterization and inference. Furthermore, the performance of the neural network surrogate for gradient approximation will greatly affect the estimation results. By using the proposed VI-DGP method with a surrogate model to solve the discontinuous random field estimation problem, we can demonstrate the feasibility and robustness of gradient approximation. To test the gradient approximation and estimation performance, we use a binary channelized field as an example due to its sharp permeability discontinuity on the channel edge.

\subsection{ GRF with uncertain correlation lengths}\label{sec:Gaussian_Examples}
\begin{figure}[h!]
    \centering
    \includegraphics[width=4.5in]{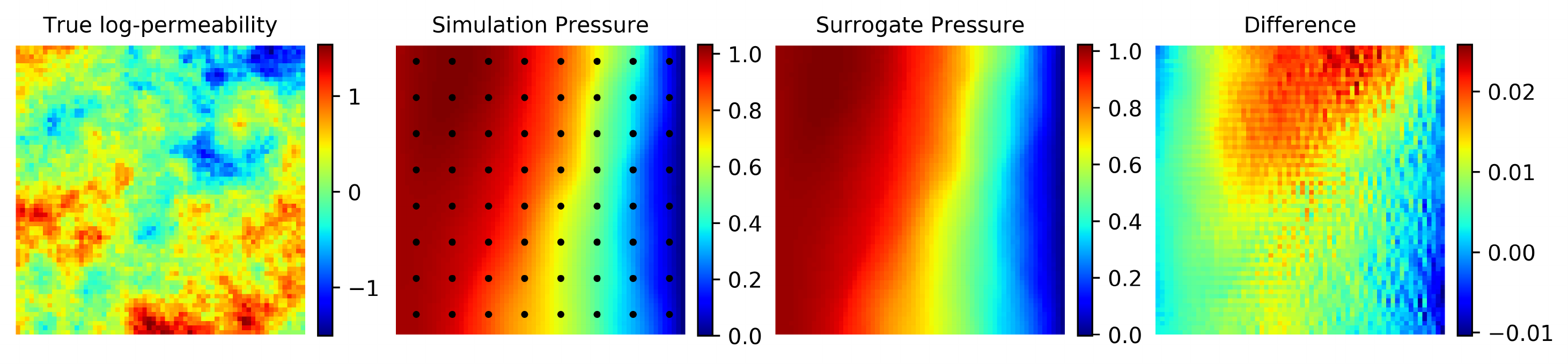}
    \caption{Illustration of the test example for the GRF. The four figures from left to right are the true log-permeability to be estimated, the corresponding pressure computed by the simulator, the corresponding pressure computed by the physics-constrained surrogate model using 4096 training data, and the difference between the two pressure results, respectively. The black dots in the second figure represent the observation locations used in BIPs.}
    \label{fig:Gau_truth}
\end{figure}

In this example, we assume that the log-permeability field is a GRF with the $L_2$ norm exponential covariance function, i.e.,
$\bk(\bx) \sim \mathcal{GP} \left( m(\bx), \operatorname{Cov}\left(\bx, \bx^{\prime}\right)\right)$, where $ m(\bx)$ and $\operatorname{Cov}\left(\bx, \bx^{\prime}\right)$ denote the mean and covariance function, respectively. $\bx = (x_1, x_2)$ and $\bx^{\prime} = (x_1^{\prime}, x_2^{\prime})$ are two arbitrary spatial locations. The $L_2$ norm exponential covariance function is
\begin{equation}
    \operatorname{Cov}\left(\bx, \bx^{\prime}\right)= \sigma_{k}^2 \exp \left(-\sqrt{\left(\frac{x_1-x_1^{\prime}}{l_{1}}\right)^{2}+\left(\frac{x_2-x_2^{\prime}}{l_{2}}\right)^{2}}\right),
    \label{eq_cov}
\end{equation}
where $ \sigma_{k}^2$ is the variance, $l_{1}$ and $l_{2}$ are the correlation lengths along the horizontal and vertical directions, respectively. We set $ m(\bx) = 0$ and $ \sigma_{k}^{2} = 0.5$. Since the KLE method cannot handle varying correlation lengths, we consider uncertain correlation lengths sampled from the uniform distribution $\mathcal{U}[0.1,0.4]$ to highlight the advantage of DGP in the prior information representation. For $10$ sampled correlation lengths, we generate $1000$ GRF samples for each correlation length. The training dataset $\{\bk^{(i)}\}_{i=1}^{N}$ for DGP naturally embodies all assumptions or prior information, where $N=10000$. The test example for the GRF case is given in Fig.~\ref{fig:Gau_truth}. The first image is the true log-permeability field to be estimated, and it is not in the training dataset of the DGP model and the surrogate model. The black dots on the  second image illustrate the collected 64 observations that are uniformly located in the pressure field. These observation locations can be denoted by the tensor product $\{x_1^i\} \otimes \{x_2^j\}$ of the one-dimensional grids: $x_1^i= 0.0625 + 0.125 i, i=0,1, \ldots, 7, x_2^j= 0.0625 + 0.125 j, j=0,1, \ldots, 7$. Our goal is to estimate the log-permeability given the noisy observations and prior information.
\subsubsection{DGP results}
Given the training dataset $\{\bk^{(i)}\}_{i=1}^{N}$, where log-permeability $\bk \in \mathbb{R}^{64\times 64}$, one can train the DGP with Algorithm~\ref{alg:DGM}. We set the latent variable $\bz \in \mathbb{R}^{256}$ to be a $256-$dimensional vector. The network architectures applied for DGP are given in Appendix~\ref{app:nn_VAE}. In this paper, all the training of neural networks and inference using surrogate models are implemented on a GPU. The GPU card used for training and inference is a single NVIDIA GeForce GTX $1080$ Ti GPU card. For the training hyperparameters in Algorithm~\ref{alg:DGM}, we set the batch size in the loss function to $n = 64$. In the optimization, the Adam optimizer~\cite{kingma2014adam} is employed with a learning rate $\eta = 0.0001$. The neural networks are trained with $300$ epochs. The training procedure takes about 15 minutes. Once the DGP is obtained, one can first sample a latent variable $\bz^{\prime}$ from Gaussian distribution $\mathcal{N}\left( \mathbf{0}, \boldsymbol{I}\right)$, and then generate the corresponding log-permeability random field $\bk^{\prime}$ by the learned decoder model $\mathcal{G}_{\btheta^{\star}}(\cdot)$, i.e., $\bk^{\prime} = \mathcal{G}_{\btheta^{\star}}(\bz^{\prime})$. With a well-trained DGP model, we can assume that generated $\bk^{\prime}$ is sampled from the underlying prior distribution $\pi(\bk)$. 

\begin{figure}[h!]
    \centering
    \includegraphics[width=4in]{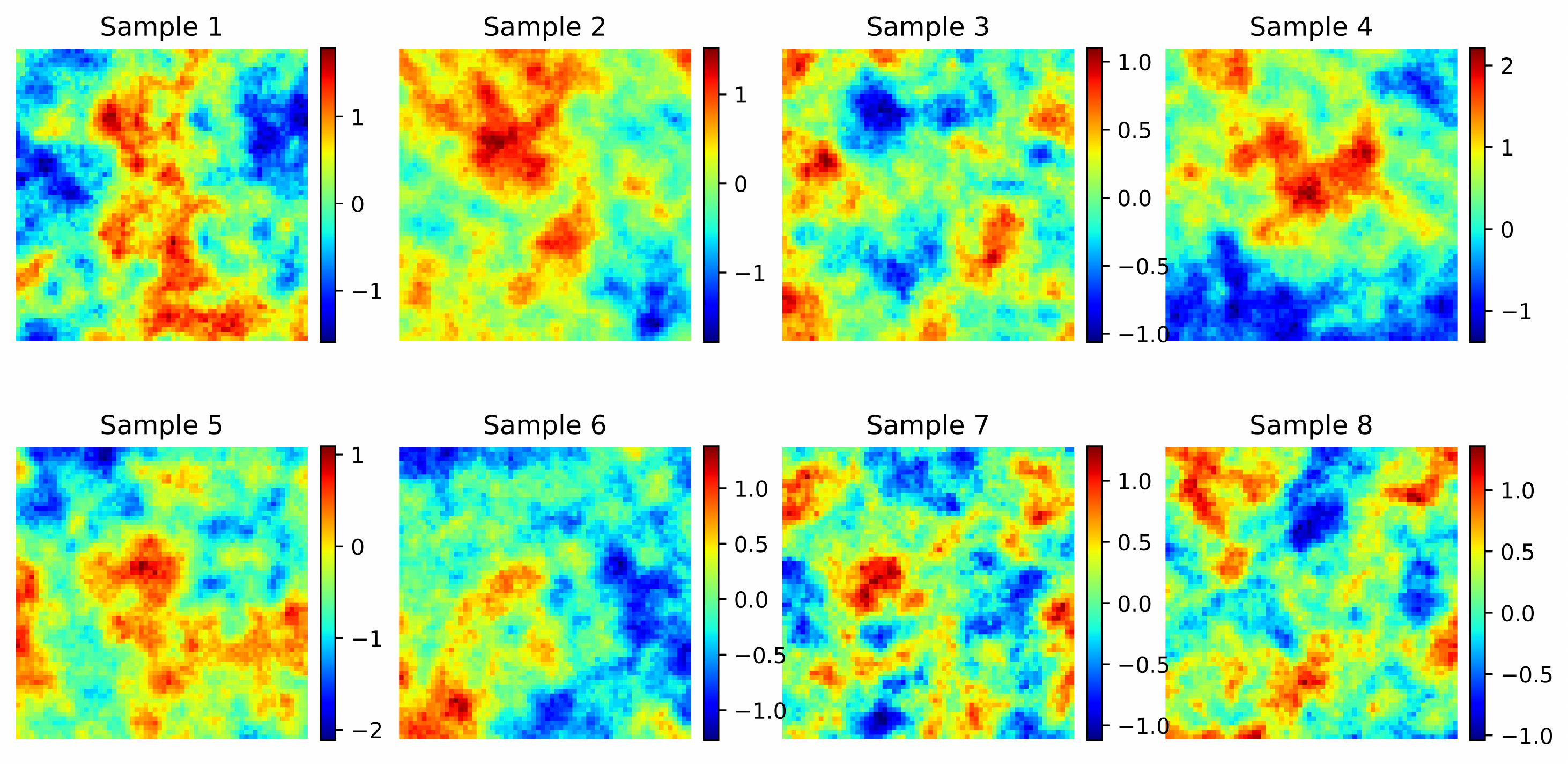}
    \caption{The prior samples generated by the learned DGP for the GRF. All samples are generated by $\bk = \mathcal{G}_{\btheta}^{\star}(\bz)$, where $\bk  \in \mathbb{R}^{64 \times 64}$, and latent variables $\bz\in \mathbb{R}^{256}$ are sampled from $\mathcal{N}(\bm{0},\bm{I})$.}
    \label{fig:Gau_vae_samples}
\end{figure}

The prior samples generated by the learned DGP are shown in Fig.~\ref{fig:Gau_vae_samples}. It is easy to find that the DGP has successfully captured the prior information of the log-permeability $\bk$ based on two facts. One is that it generates various GRF realizations that resemble those given in the training dataset. Moreover, for the given 8 samples in Fig.~\ref{fig:Gau_vae_samples}, it is obvious that sample 2 and sample 4 have long correlation lengths, while sample 7 and sample 8 show short correlation lengths. Their diverse correlation lengths are consistent with our setup that uncertain correlation lengths are sampled from the uniform distribution $\mathcal{U}[0.1,0.4]$. So the well-trained DGP can learn the features of the varying correlation lengths. Using this well-trained DGP, one can estimate the posterior distribution of the latent variable in Eq.~\eqref{eq:PosteriorOfz}, then generate posterior samples of the log-permeability $\bk$ with  the posterior distribution $\pi(\bz | \bd)$ and generative model $\mathcal{G}_{\btheta^{\star}}(\cdot)$.

\subsubsection{Gradient approximation results}
The most computational cost in the VI-DGP method involves the forward computation and its corresponding gradient computation. To accelerate the inference, we propose a gradient approximation method using the neural network surrogate to replace the adjoint method in the VI-DGP. The neural network surrogate has automatic differentiation and is extremely fast with deep learning frameworks like Pytorch. For the test problem, given $\{\bx_{\mathcal{D}}^{(i)}\}_{i=1}^{n_p}$, $ \{\bx_{\partial\mathcal{D}}^{(i)}\}_{i=1}^{n_b}$, and training data $\{\bk^{(j)}(\bx)\}_{j=1}^{N_k}$, we can rewrite the loss function $J_{\text{pde}}(u(\bx, \bk ; \Theta))$ and $J_{\text{b}}(u(\bx, \bk ; \Theta))$ into a discretized form to learn the surrogate model. The detailed discretized loss functions for Darcy flow and the network architectures applied in this paper are given in Appendix~\ref{app:nn_pcs}. 

When the 2D unit square domain $\mathcal{D}=[0,1]^2$ in Eq.~\eqref{eq:darcy} is discretized into uniform $64 \times 64$ grids in advance, $\{\bx_{\mathcal{D}}^{(i)}\}_{i=1}^{n_p}$ and $ \{\bx_{\partial\mathcal{D}}^{(i)}\}_{i=1}^{n_b}$ are naturally defined. Using Algorithm~\ref{alg:PCS}, we can train the surrogate model using only  the input data, i.e., the log-permeability dataset $\{\bk^{(j)}(\bx)\}_{j=1}^{N_k}$. The penalty  parameter $\gamma$ in Eq.~\eqref{eq:pcs_loss} is 10, which requires predictions to satisfy the boundary conditions. The batch size $n_s$ is 32. The networks are trained for 300 epochs using the Adam optimizer paired with one cycle policy (learning rate scheduler), where the maximum learning rate is 0.001.  We train the surrogate with different numbers of training data to test its effect on the gradient approximation. When the number of training data $N_k$ is 1024, 2048, and 4096, their corresponding training time is about $8.3$, $16.2$, and $31.9$ minutes, respectively. Unlike the previous surrogate model, whose evaluation emphasizes the error or relative error between the surrogate predictions and simulation outputs, we focus on gradient approximation using neural networks. It relates to whether the surrogate model can replace the adjoint method in optimization. We adopt the stochastic gradient descent/ascent in the VI-DGP method to reduce the computational burden.  Since stochastic gradient descent/ascent only requires an appropriate descent/ascent direction rather than an exact gradient, it relaxes strict constraints on gradient accuracy in the VI-DGP method. As long as the approximate gradient can provide an appropriate direction for the lower bound optimization in each iteration, it will converge and obtain a good approximation.

\begin{figure}[h!]
    \centering
    \includegraphics[width=2.5in]{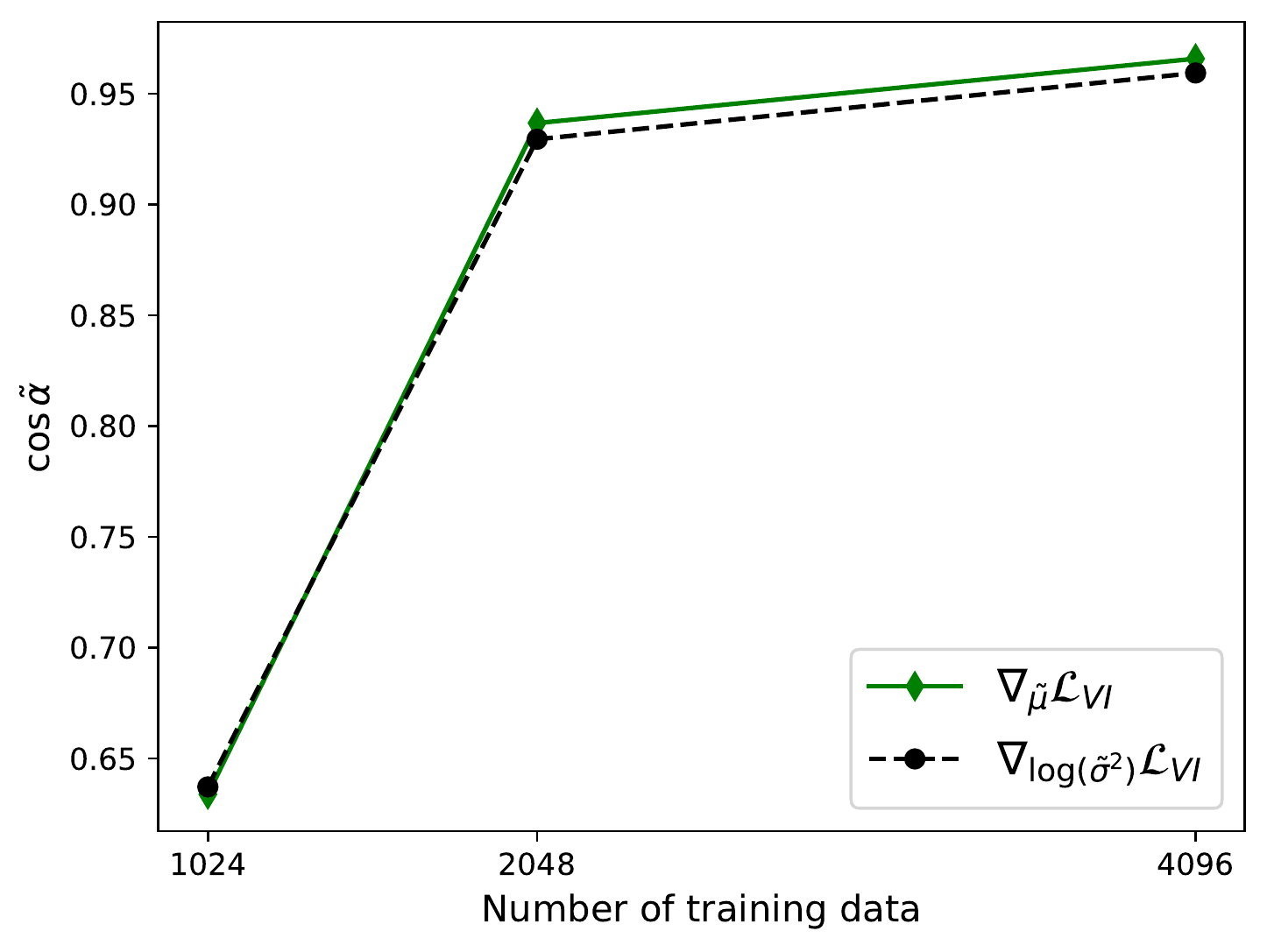}
    \caption{The computed $\cos\tilde{\alpha}$ with different surrogates in the GRF case. The surrogates are trained using $1024$, $2048$, and $4096$ training data, respectively. The green solid line shows the results of $\nabla_{\tilde{\bmu}}\mathcal{L}_{VI}$. The black dashed line shows the results of $\nabla_{\log(\tilde{\bsigma}^2)} \mathcal{L}_{VI}$. }
    \label{fig:Gau_cos}
\end{figure}

\begin{figure}[h!]
    \centerline{
        \begin{tabular}{cc}
    \includegraphics[scale=0.35]{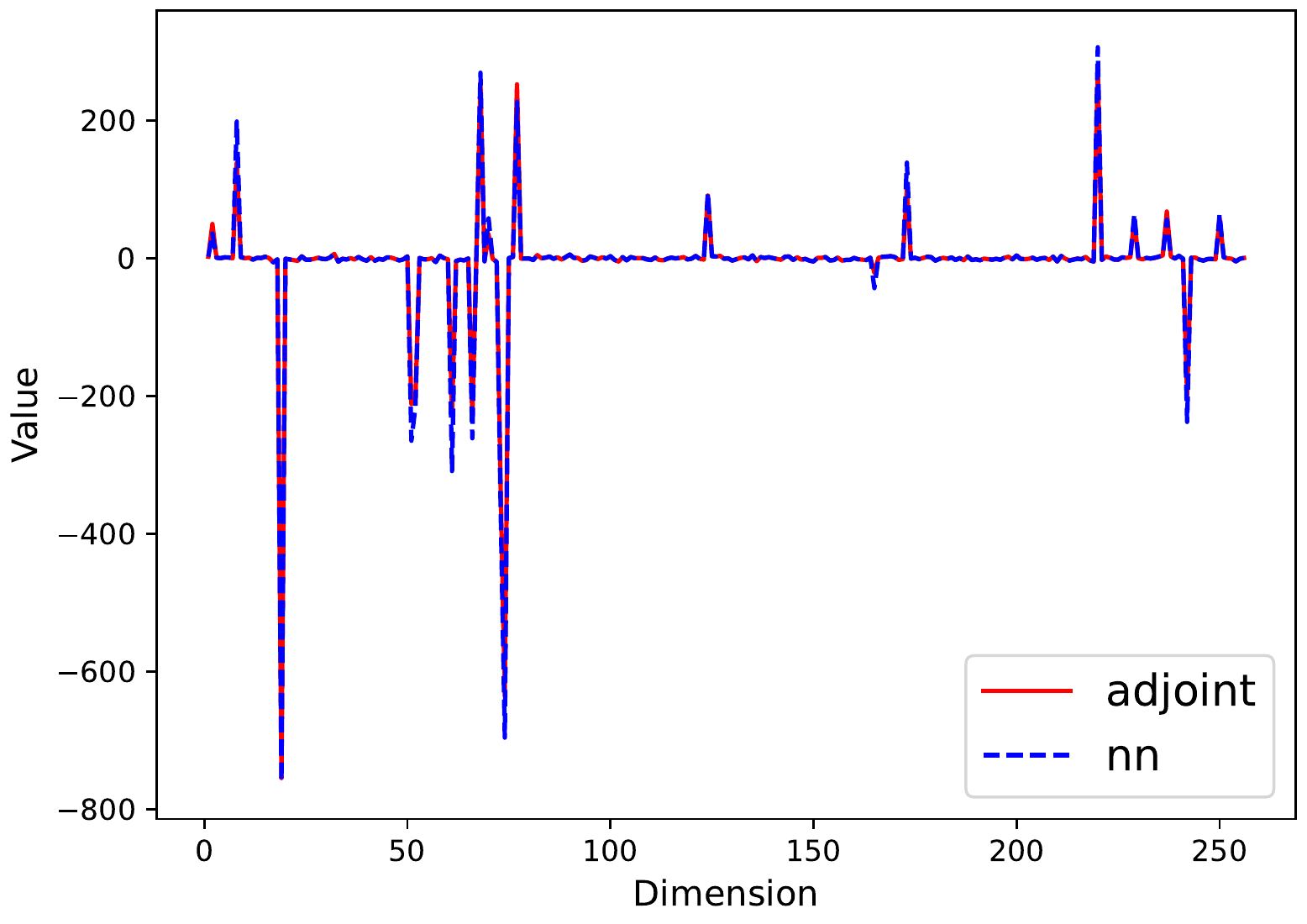}
    & 
    \includegraphics[scale=0.35]{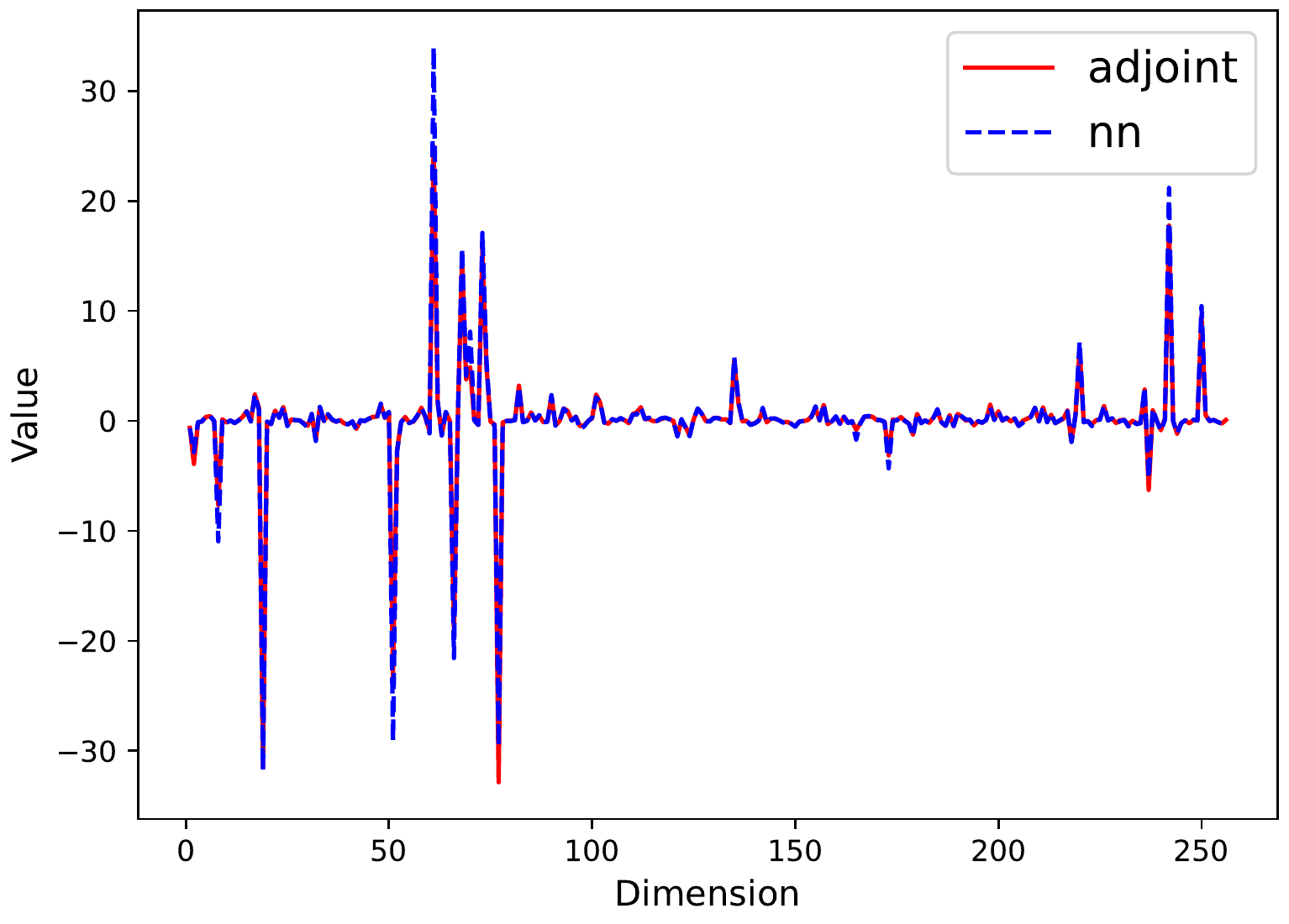}
    \\
    (a) $\nabla_{\tilde{\bmu}}\mathcal{L}_{VI}$& 
    (b) $\nabla_{\log(\tilde{\bsigma}^2)} \mathcal{L}_{VI}$
\end{tabular}}
\caption{An example of the gradient computation by the adjoint method and neural network in the GRF case. Given the $\tilde{\bmu}$ and $\log(\tilde{\bsigma}^2)$, the vector (a)$\nabla_{\tilde{\bmu}}\mathcal{L}_{VI}$ (b)$\nabla_{\log(\tilde{\bsigma}^2)} \mathcal{L}_{VI}$ are computed by the workflow in Fig.~\ref{fig:AD_ELBO}, where the $\frac{\partial \mathcal{L}_{VI}}{\partial \bk}$ in the workflow are computed by the neural networks (blue dashed line) and the adjoint method (red solid line), respectively.}
\label{fig:Gau_cos_example}
\end{figure}

Although the noisy gradient leads to the optimization not being the steepest descent/ascent, it is a trade-off between the convergence rate and computational cost. To assess the feasibility of using gradient approximation from the learned neural networks to accelerate the computation of $\frac{\partial \mathcal{L}_{VI}}{\partial \bk}$ in Fig.~\ref{fig:AD_ELBO}, for any given parameter $\tilde{\bmu}$ and $\log(\tilde{\bsigma}^2)$, we compute the $\nabla_{\tilde{\bmu}}\mathcal{L}_{VI}$ and $\nabla_{\log(\tilde{\bsigma}^2)} \mathcal{L}_{VI}$ in Algorithm~\ref{alg:VIM} using the adjoint method and the learned neural networks simultaneously. Inspired by the computation of the angle between two vectors with respect to the Euclidean norm, we define the evaluation metric for gradient approximation as follows:
\begin{equation}
    \cos\tilde{\alpha}=\frac{1}{N_g}\sum_{i=1}^{N_g}\frac{\bm{g}_{nn}^{(i)}\cdot \bm{g}_{a}^{(i)}}{\|\bm{g}_{nn}^{(i)}\|_{2}\|\bm{g}_{a}^{(i)}\|_{2}},
    \label{eq:cos}
\end{equation}
where $N_g$ denotes the number of samples for evaluation, given the $i-$th parameter sample, $\bm{g}_{nn}^{(i)}$ and $\bm{g}_{a}^{(i)}$ are the gradient computed by the neural network surrogate and the adjoint method, respectively. If $N_g\to \infty$ and $\cos\tilde{\alpha} = 1$, almost all of the gradients computed by the neural network surrogate will keep the same direction as those computed by the adjoint method. To test the gradient approximation, we consider the test example in Fig.~\ref{fig:Gau_truth} where observations are corrupted with $5\%$ independent additive Gaussian random noise. Using the pre-trained DGP, we sample 1000 pairs of $\tilde{\bmu}$ and $\log(\tilde{\bsigma}^2)$ from Gaussian distribution $\mathcal{N}(\bm{0},\bm{I})$, then compute their corresponding stochastic gradient $\nabla_{\tilde{\bmu}}\mathcal{L}_{VI}$ and $\nabla_{\log(\tilde{\bsigma}^2)} \mathcal{L}_{VI}$ with sampling number $M_s=1$ using  the neural network surrogate and the adjoint method, respectively. The computed $\cos\tilde{\alpha}$ is shown in Fig.~\ref{fig:Gau_cos}. The results reflect that choosing an appropriate number of training data for the surrogate model is essential for gradient approximation in the VI-DGP method.  When $N_k =4096$, the $\cos\tilde{\alpha}$ is around $0.95$, which means most of the gradients computed by the surrogate model keep a relatively consistent direction with the gradients computed by the adjoint method. The third image in Fig.~\ref{fig:Gau_truth} shows the pressure prediction using the surrogate model trained with $4096$ training data, and the fourth image suggests its good performance. Given a certain pair of pairs $\tilde{\bmu}$ and $\log(\tilde{\bsigma}^2)$, Fig.~\ref{fig:Gau_cos_example} depicts their corresponding gradient, where blue dashed line and red solid line are computed by the surrogate model trained with $4096$ training data and the adjoint method, respectively. It shows that two vectors keep coincident in most dimensions for both $\nabla_{\tilde{\bmu}}\mathcal{L}_{VI}$ and $\nabla_{\log(\tilde{\bsigma}^2)} \mathcal{L}_{VI}$. Based on above results, it is reasonable to employ the surrogate model trained with $4096$ training data for the following Bayesian inversion task.

\subsubsection{Bayesian inversion results}
In this section, we will discuss the performance of the proposed VI-DGP method for solving BIPs. We will present the results in three aspects. First, we will compare the estimated results obtained using various methods, including the VI-DGP method with the neural network surrogate (VI-NN), the VI-DGP method with the adjoint method (VI-adjoint), the MCMC method with the neural network surrogate (MCMC-NN), and the MCMC method with the finite element method (MCMC-FEM). We will present and analyze their corresponding results in terms of accuracy and efficiency. Second, as discussed in Section~\ref{sec:VI}, a good sampling number $M_s$ requires the trade-off between convergence rate and computational cost. We will show the convergence and estimated results under different $M_s$. Lastly, we will investigate the robustness of the VI-DGP method under different noise levels.

\emph{Comparisons.} Four methods are applied for the test problem given in  Fig.~\ref{fig:Gau_truth}, where the observations are added with $5\%$ independent Gaussian random noise. We shall see the performance of the VI-DGP method and the impact of the trained surrogate model. Using the pre-trained generative model $\bk = \mathcal{G}_{\btheta}(\bz)$, where $\bk \in \mathbb{R}^{64 \times 64}$, and $\bz \in \mathbb{R}^{256}$, we implement the VI-DGP method with Algorithm~\ref{alg:VIM}. For the GRF case, we set the optimization iteration $N_{opt}$ to $5000$, the sampling  number $M_s$ to $1$, and the number of posterior samples $N_s$ to $10000$. We adopt the SGD optimizer in the Pytorch library with the learning rate $\eta_{\tilde{\bmu}}=\eta_{\tilde{\bsigma}} = 0.0008$. The initial values for $\tilde{\bmu}$ and $\log(\tilde{\bsigma}^2)$ are both zero vectors. For the MCMC, we use preconditioned Crank–Nicolson (pCN)  algorithm~\cite{cotter2013mcmc,hairer2014spectral} for the posterior approximation. The specific details of the algorithm can be found in Appendix~\ref{app:pcn}. We run a Markov chain for $50000$ steps and use the last $10000$ steps as the posterior samples. Table~\ref{table:Gaussian cost} presents the computational cost for the posterior approximation with the four implemented methods. Even though we use the first-order element for fast simulation in this experiment, the proposed VI-DGP method is still faster than other methods, both with and without the surrogate model. With the GPU acceleration, it only takes $94$ seconds to run $5000$ iterations using the pre-trained surrogate model. If the simulation involves a complex physics system or a large-scale problem with a high-order element, the computational cost will be unaffordable for the MCMC-FEM method and the VI-adjoint method. For inference efficiency, the VI-DGP method using gradient approximation has significant advantages.
\begin{table}[h!]
	\caption{Computational cost of estimation with different methods in GRF case.} 
	\centering	
	\begin{tabular}{ccccc}  
		\hline
	Methods	&VI-NN & VI-adjoint   & MCMC-NN &  MCMC-FEM\\\hline
	 	Iterations & $5000$ & $5000$ & $50000$ & $50000$  \\ 
	 	Inference time (s) & $\bm{94}$ & $2093$  & $353$ & $12894$  \\
		\hline
	\end{tabular}
	\label{table:Gaussian cost}
\end{table}

\begin{figure}[h!]
    \centering
    \includegraphics[width=4.5in]{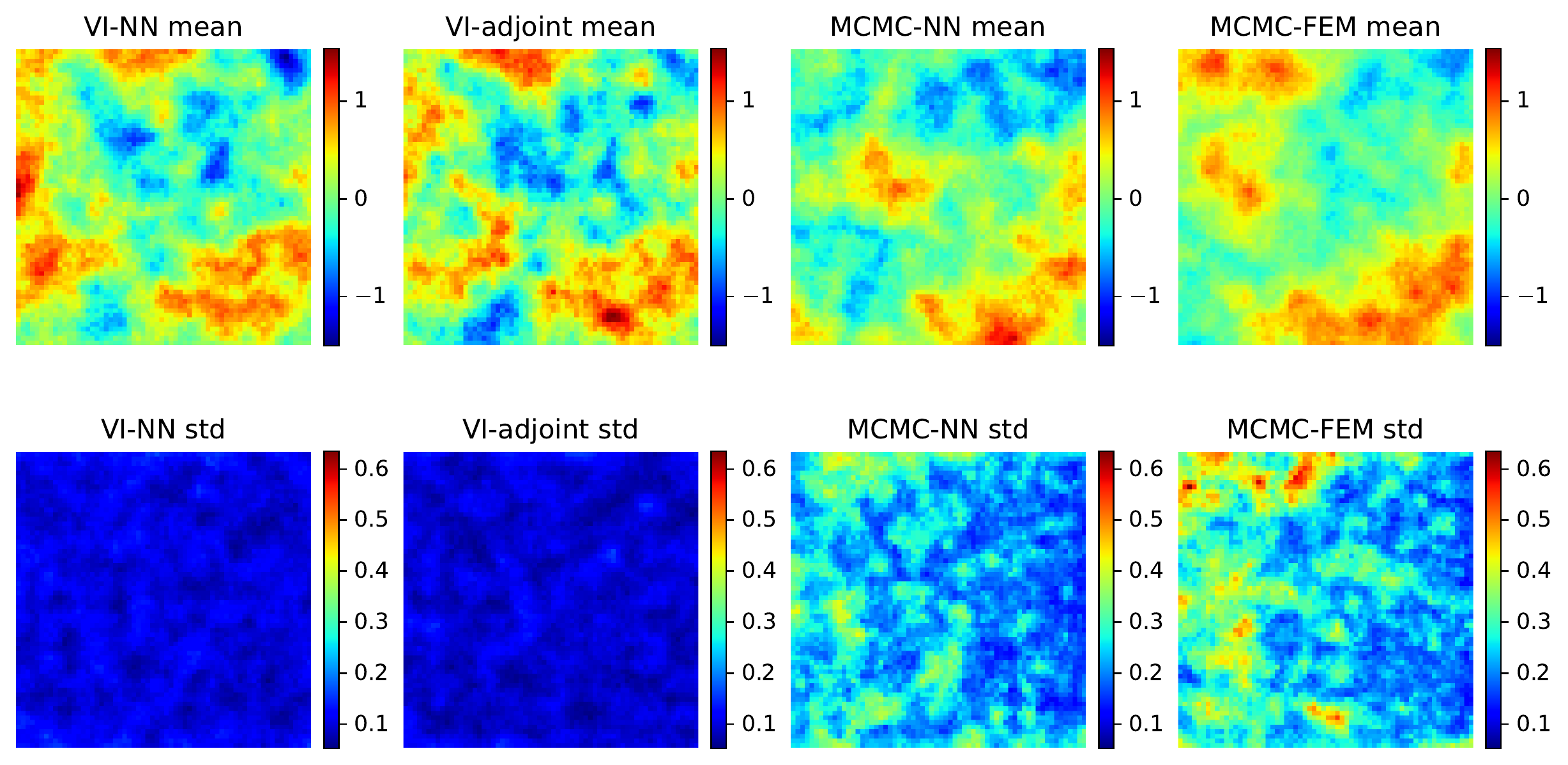}
    \caption{The posterior estimation results for the GRF with different methods. The first row shows the estimated mean of the log-permeability field, and the second row gives the corresponding standard deviation (std).}
    \label{fig:Gau_posterior}
\end{figure}

Fig.~\ref{fig:Gau_posterior} provides the estimated results obtained by the above four methods. The computed mean and the standard deviation using the posterior samples are given in the first and the second row, respectively. Four methods produce comparably good mean results on the right region in comparison to the true log-permeability, as their main features are captured by them. However, the VI-DGP method can achieve a better mean result on the left region. The standard deviation results generated by the MCMC method are significantly higher than those computed by the VI-DGP method, indicating high uncertainty in the posterior estimation when using the MCMC method. Moreover, we know that variational inference tends to underestimate the uncertainty of the posterior distribution, this is a result of its objective function~\cite{blei2017variational}. With a good mean result, underestimating the variance may also be acceptable. When the strategy for optimization stability discussed in Section~\ref{sec:VI} is applied, the results of the VI-NN and the VI-adjoint show that a small sampling number (even $M_s=1$) can still realize a good estimation. For the smooth GRF, the similar results of the VI-NN and the VI-adjoint indicate that a well-trained surrogate model can replace the adjoint method on the gradient approximation even using a small $M_s$. 

\begin{figure}[h!]
    \centerline{
        \begin{tabular}{cc}
    \includegraphics[width=2.2in]{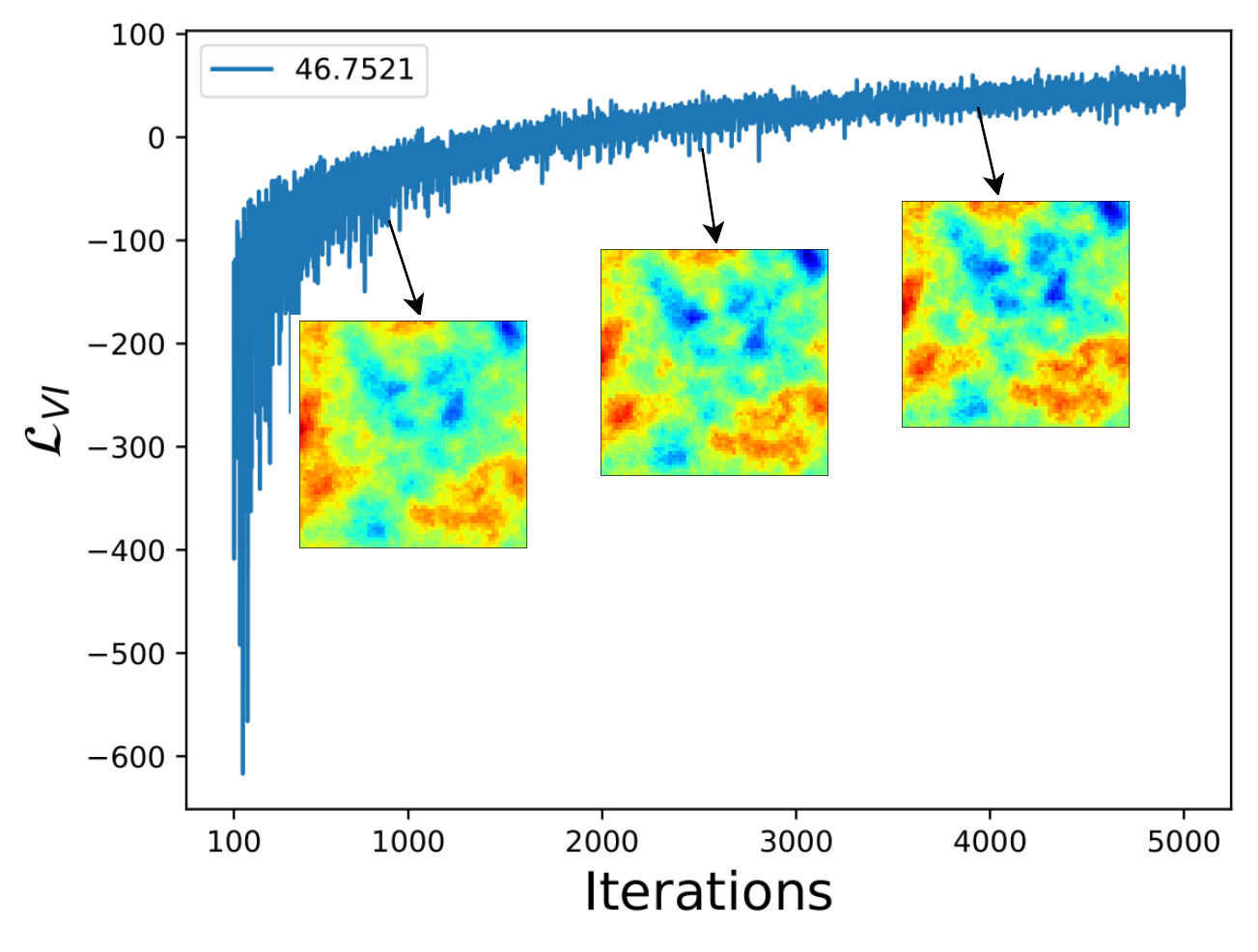}
    & 
    \includegraphics[width=2.2in]{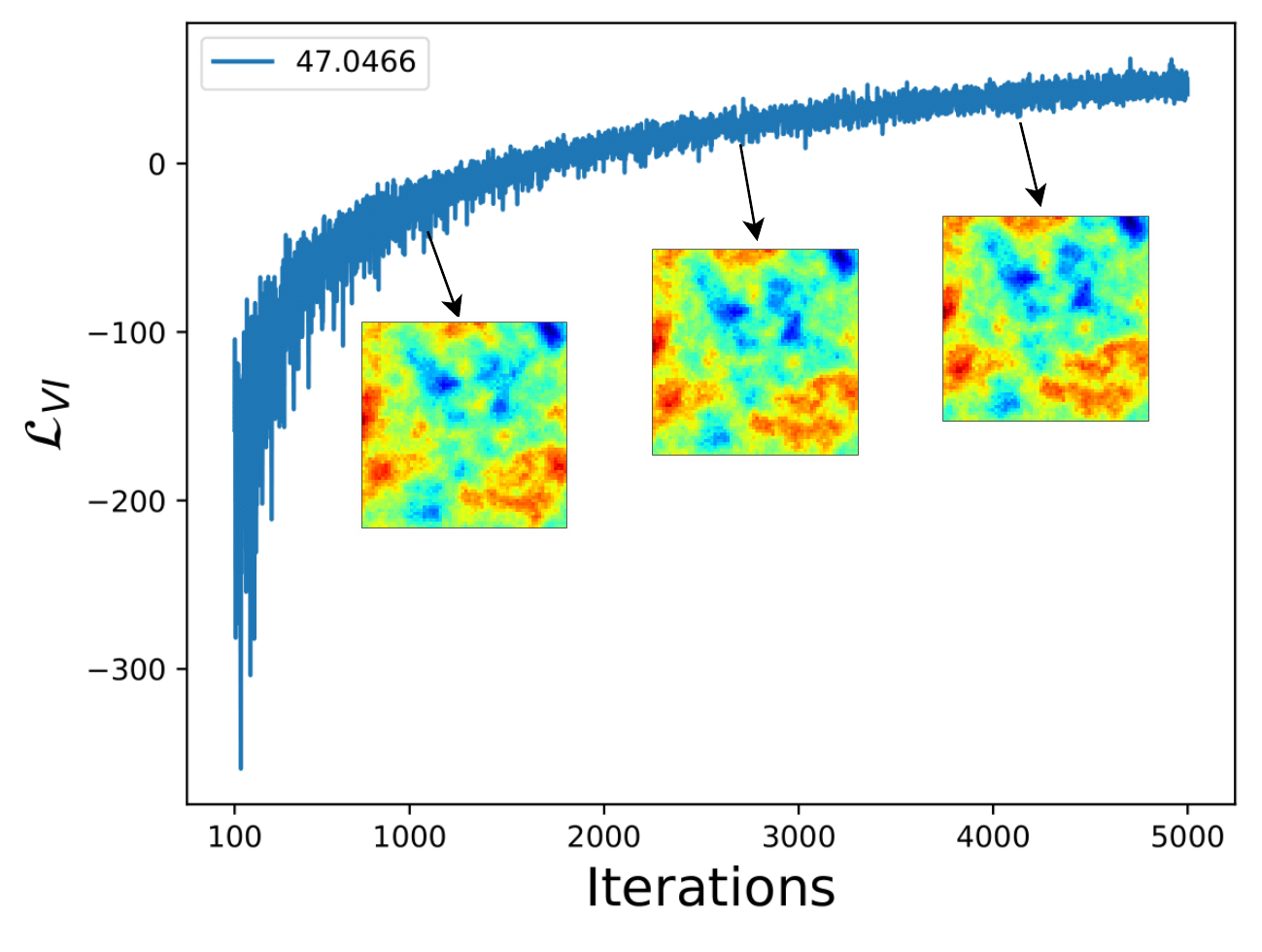}
    \\
    (a) $M_s=1$& 
    (b) $M_s=3$\\   
     \includegraphics[width=2.2in]{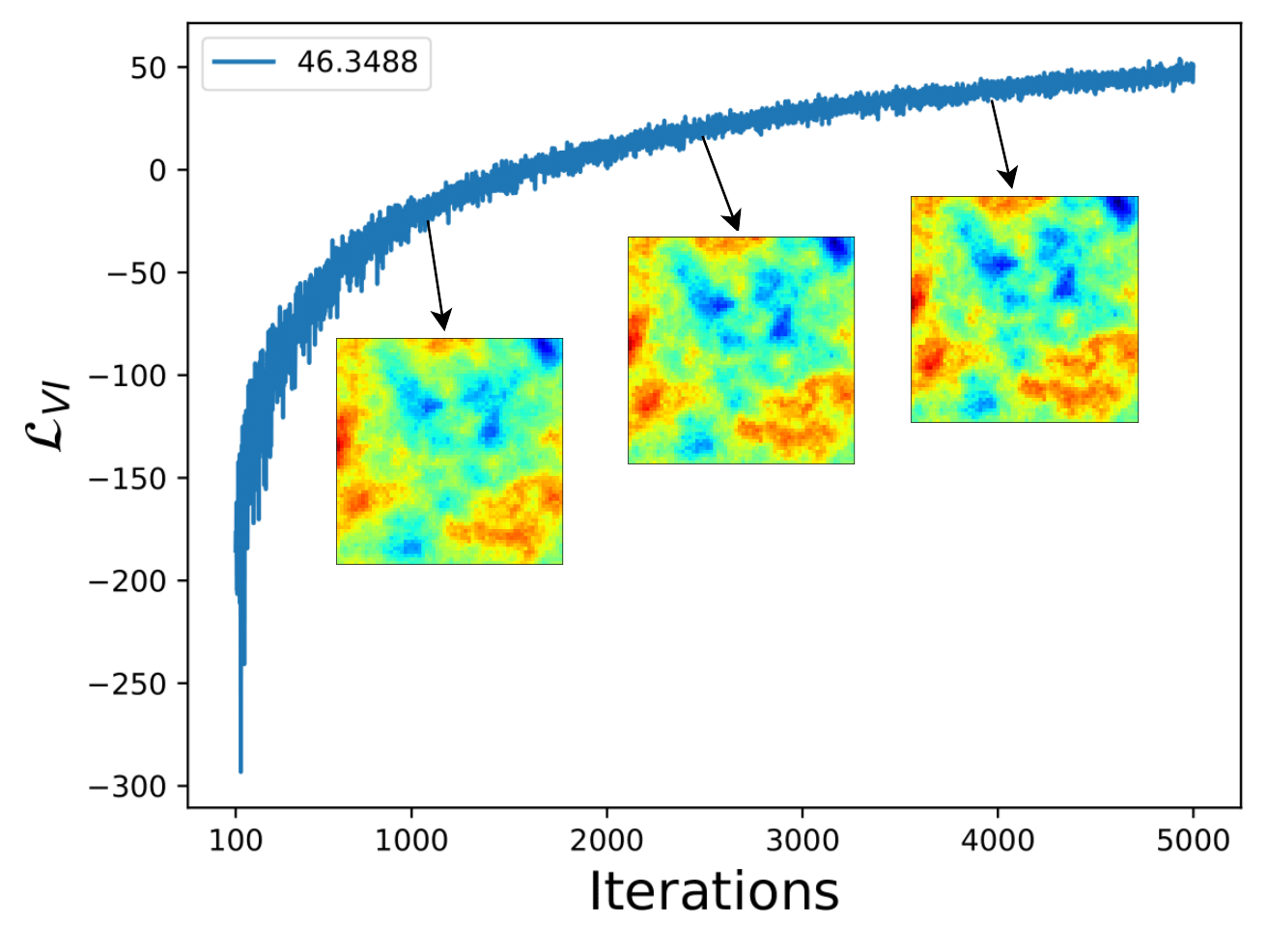} &
    \includegraphics[width=2.2in]{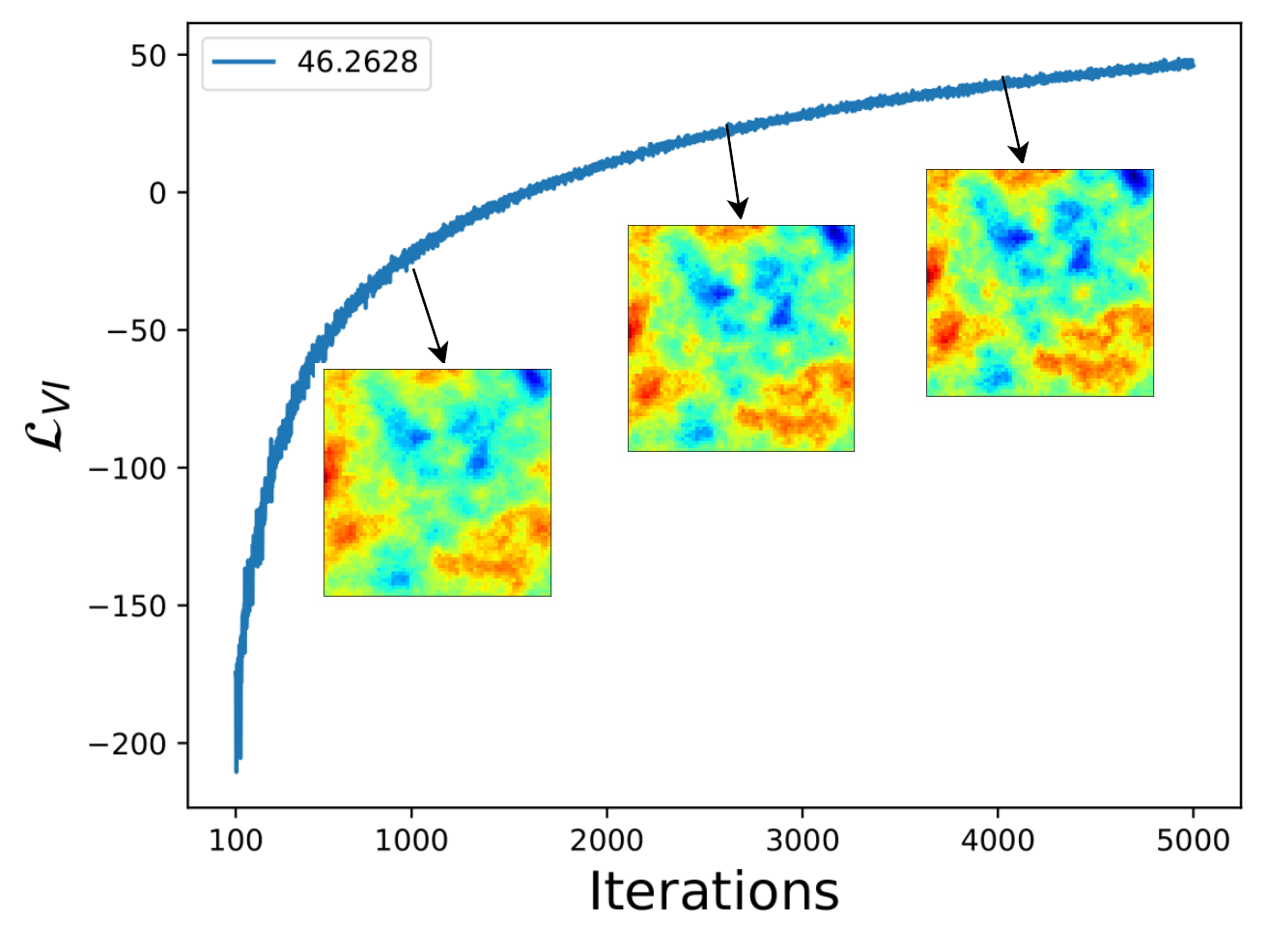}\\
    (c) $M_s=10$& 
    (d) $M_s=100$\\  
\end{tabular}}
\caption{The convergence of the variational lower bound $\mathcal{L}_{VI}$ with varying sample numbers $M_s$ in stochastic optimization. The three log-permeability fields below the black arrows are the estimated mean at the $1000$-th, $2500$-th, and $4000$-th iteration, respectively.}
	\label{fig:Gau_elbo}
\end{figure}

\emph{Effect of the sampling number $M_s$.} 
In order to investigate the impact of the sampling number $M_s$ on the optimization convergence of the VI-DGP method, we implement the above VI-NN experiment for the given test problem. All of the configurations are the same, except for the sampling number $M_s$. Four different sampling numbers are considered in our experiments. The results of the variational lower bound $\mathcal{L}_{VI}$ and the estimated mean at some specific iterations are shown in Fig.~\ref{fig:Gau_elbo}. It is clear that a larger sampling number leads to a more stable convergence of the variational lower bound $\mathcal{L}_{VI}$. However, after $5000$ iterations, the values of $\mathcal{L}_{VI}$ of four experiments are similar, around $46$. This suggests that the choice of a small sampling number only affects the convergence process. Once the optimization has converged, the obtained estimation results are similar. Note that using a sampling number of $M_s=100$ for the VI-DGP method will result in $100$ times the computational cost compared to when $M_s=1$ is used. Based on the estimation results and computational cost, a small sampling number is a feasible and better choice for implementation.

\begin{figure}[h!]
    \centerline{
        \begin{tabular}{cc}
    \includegraphics[width=2.2in]{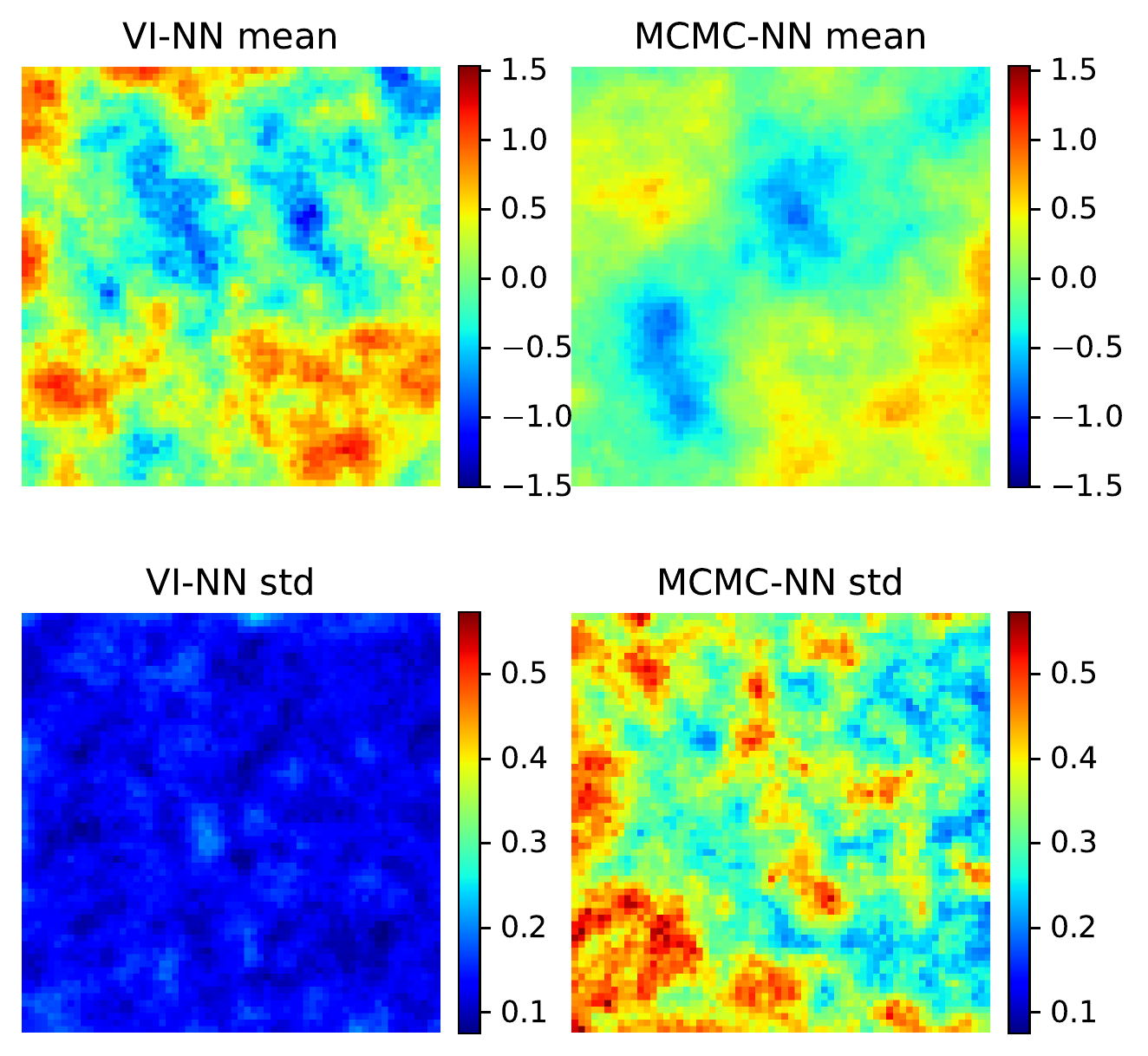}
    & 
    \includegraphics[width=2.2in]{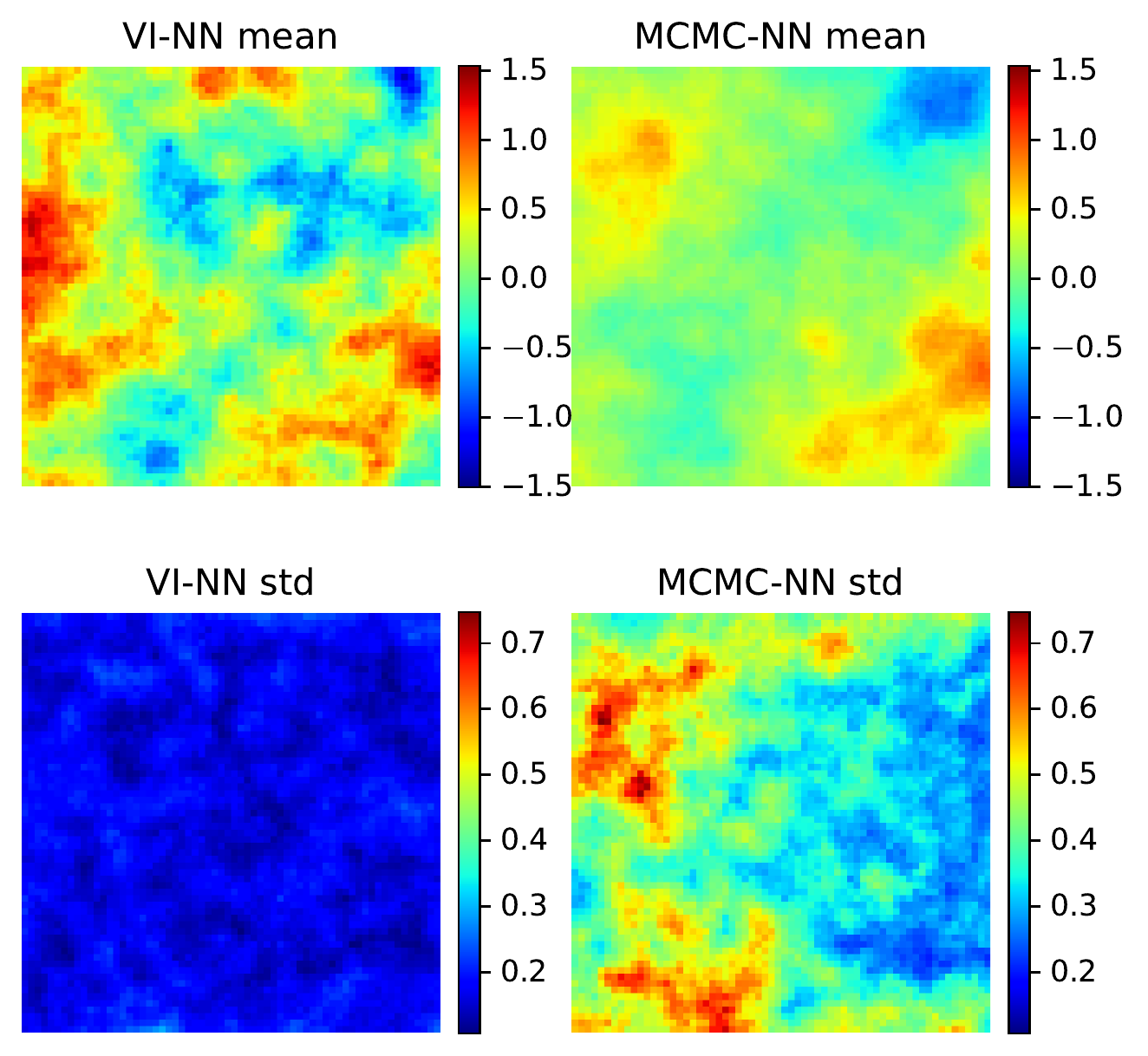}
    \\
    (a) $7\%$ noise& 
    (b) $10\%$ noise\\   
       
\end{tabular}}
\caption{The posterior estimation results in GRF case using VI-NN and MCMC-NN method under (a) $7\%$ noise and (b) $10\%$ noise observations. The first row shows the estimated mean of the log-permeability field, and the second row gives the corresponding standard deviation (std).}
	\label{fig:Gau_noise}
\end{figure}

\emph{Effect of the observation noise.}
Keeping the same configurations as the experiments in Fig.~\ref{fig:Gau_posterior}, we evaluate the robustness of the proposed method using two additional observation setups with higher levels of noise. $7\%$ and $10\%$ independent Gaussian random noise are imposed on the 64 pressure observations. Using these observations, we infer the log-permeability field using the VI-NN and MCMC-NN methods. The estimated results are shown in Fig.~\ref{fig:Gau_noise}. It is clear that the VI-DGP method still achieves a good estimation even though high noise is provided. In contrast, the MCMC method almost failed on such a difficult task. Although their estimated mean can still capture the feature on the right region with much lower values, the results are still substantially different from the true log-permeability field. The estimated standard deviation with MCMC is very high, while the VI-DGP method results present much lower uncertainty.

\subsection{Binary channelized field}\label{sec:Channel_Examples}

In this test example, we are focused on the estimation of the non-Gaussian log-permeability parameter. The challenges are two-fold: first, the parameterization for the complex non-Gaussian parameters is still challenging and requires further development. Second, inferring these non-Gaussian parameters is challenging due to their spatially correlated properties, even when using methods with high computational costs.  We use the binary channelized field to demonstrate the capabilities of the DGP representation and evaluate the efficiency and accuracy of the VI-DGP method for non-Gaussian parameter estimation. Additionally, using a neural network surrogate gradient approximation for complex and discontinuous field estimation may cause additional issues. We can examine its performance in terms of gradient computation and estimation.
 
\begin{figure}[h!]
    \centering
    \includegraphics[width=4.5in]{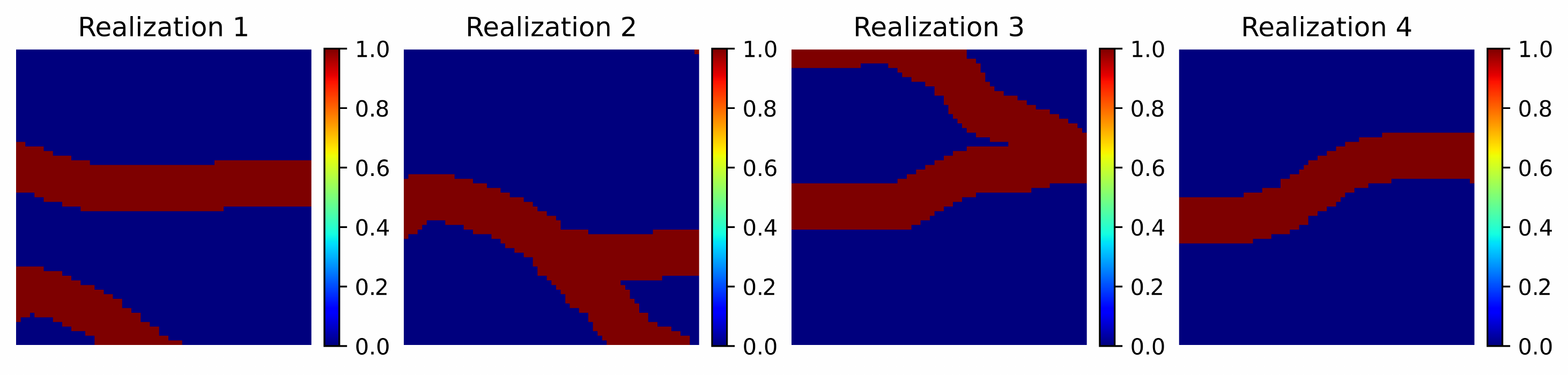}
    \caption{The randomly sampled realizations in the dataset for the DGP training.}
    \label{fig:Chan_training_samples}
\end{figure}

Suppose that the prior information of the binary channelized field is based on the historical data, which is a large image~\cite{laloy2018training}  of  size $2500 \times 2500$. One can crop small images, the size of $64 \times 64$, from this large image using a fixed $16-$pixel stride in both the horizontal and vertical directions. To obtain sufficient training data, we flip the entries in each row of the image in the left/right direction using the \verb|fliplr| operation~\footnote{\href{https://numpy.org/doc/1.18/reference/generated/numpy.fliplr.html}{https://numpy.org/doc/1.18/reference/generated/numpy.fliplr.html}} in the Numpy package to obtain a new image, and then we crop this image in the same way.  We use $40000$ images out of $46208$ cropped samples as the training dataset for the DGP model. Fig.~\ref{fig:Chan_training_samples} depicts four examples from the training dataset $\{\bk\}_{i=1}^{N}$. 

\begin{figure}[h!]
    \centering
    \includegraphics[width=4.5in]{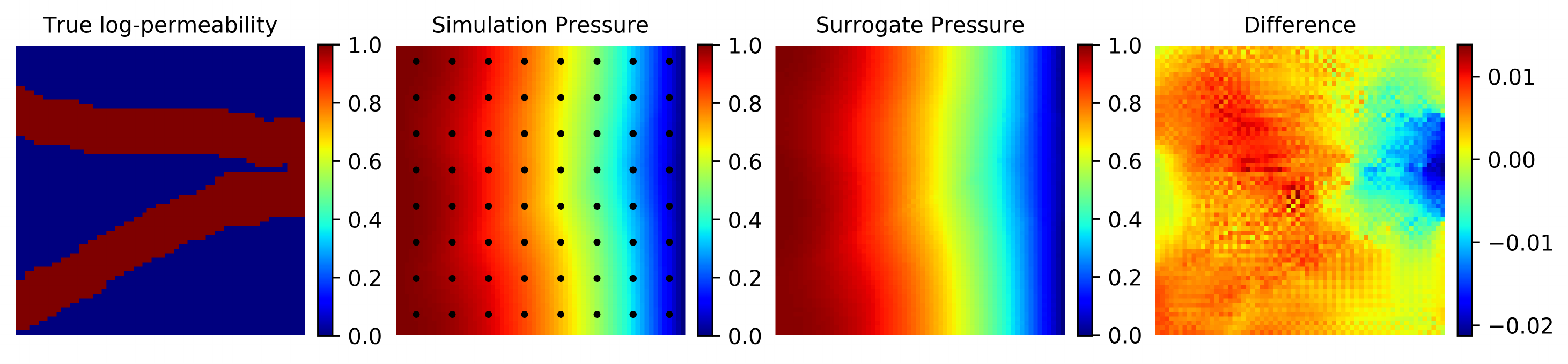}
    \caption{Illustration of the test example for the binary channelized field. The four figures from left to right are the true log-permeability to be estimated, the corresponding pressure computed by the simulator, the corresponding pressure computed by the physics-constrained surrogate model using 4096 training data, and the difference between the two pressure results, respectively. The black dots in the second figure represent the observation locations used in BIPs.}
    \label{fig:chan_truth}
\end{figure}

For the Bayesian inversion task, the unknown true log-permeability is not included in the training dataset. Fig.~\ref{fig:chan_truth} presents the test example of the binary channelized field. The red and blue regions in the first image represent the high- and low-permeability values, respectively. The observations located on the pressure field are computed by the simulator, which is shown in the second image. The inversion task is to estimate the true channels based on these noisy observations.
\subsubsection{DGP results}
\begin{figure}[h!]
    \centering
    \includegraphics[width=4in]{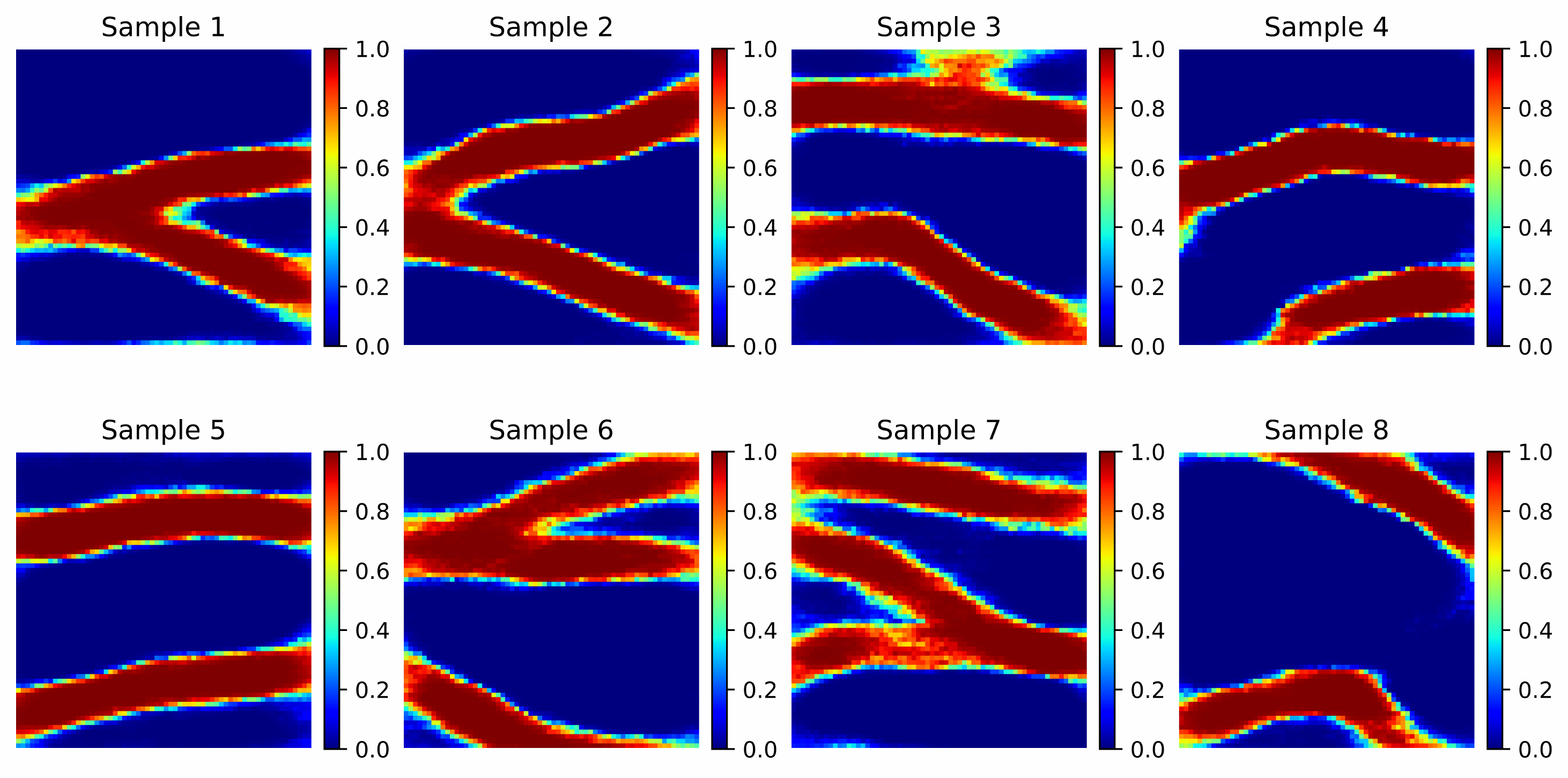}
    \caption{The prior samples generated by the learned DGP for the binary channelized field. All samples are generated by $\bk = \mathcal{G}_{\btheta}^{\star}(\bz)$, where $\bk  \in \mathbb{R}^{64 \times 64}$, and latent variables $\bz\in \mathbb{R}^{512}$ are sampled from  $\mathcal{N}(\bm{0},\bm{I})$.}
    \label{fig:Chan_vae_samples}
\end{figure}

Using the cropped $40000$ images as the prior information, we train the DGP with Algorithm~\ref{alg:DGM} and the network architectures described in Appendix~\ref{app:nn_VAE}. Here, the hyperparameters are the same as in the GRF case. The only difference is the latent variable $\bz$, where $\bz \in \mathbb{R}^{512}$.  To keep continuous channels and capture the diversity, we choose a higher dimension to relieve information compression. The DGP training for binary channelized fields takes approximately $62$ minutes. Fig.~\ref{fig:Chan_vae_samples} shows $8$ random samples generated by the learned DGP model. These prior samples keep continuous channels and resemble the training dataset realizations as shown in Fig.~\ref{fig:Chan_training_samples}, even though the values on the field are not binary, especially on the channel edge. The learned DGP provides enough prior information and can be applied in the posterior estimation. 

\subsubsection{Gradient approximation results}
In the binary channelized case, the training dataset $\{\bk^{(j)}(\bx)\}_{j=1}^{N_k}$ for surrogate training is a subset of the training dataset for DGP training. We also choose $1024$, $2048$, and $4096$ for $N_k$ to test the relationship between training data and gradient approximation. The setups and hyperparameters are the same as in the GRF case, except for the learning rate. Based on the discontinuous features of the log-permeability field, we adopt a small learning rate for the Adam optimizer, where the maximum learning rate is $0.0001$. The training time is about $8.4$, $16.3$, and $32.2$ minutes for three training data setups, respectively. Using the learned surrogate mode with $4096$ training data, we predict the pressure field of the given true log-permeability, as shown in the third image in Fig.~\ref{fig:chan_truth}. The maximum absolute error between the simulation output and surrogate prediction is only about $0.02$. It indicates that the learned surrogate model can make a good prediction for forward computation.

\begin{figure}[h!]
    \centering 
    \includegraphics[width=2.5in]{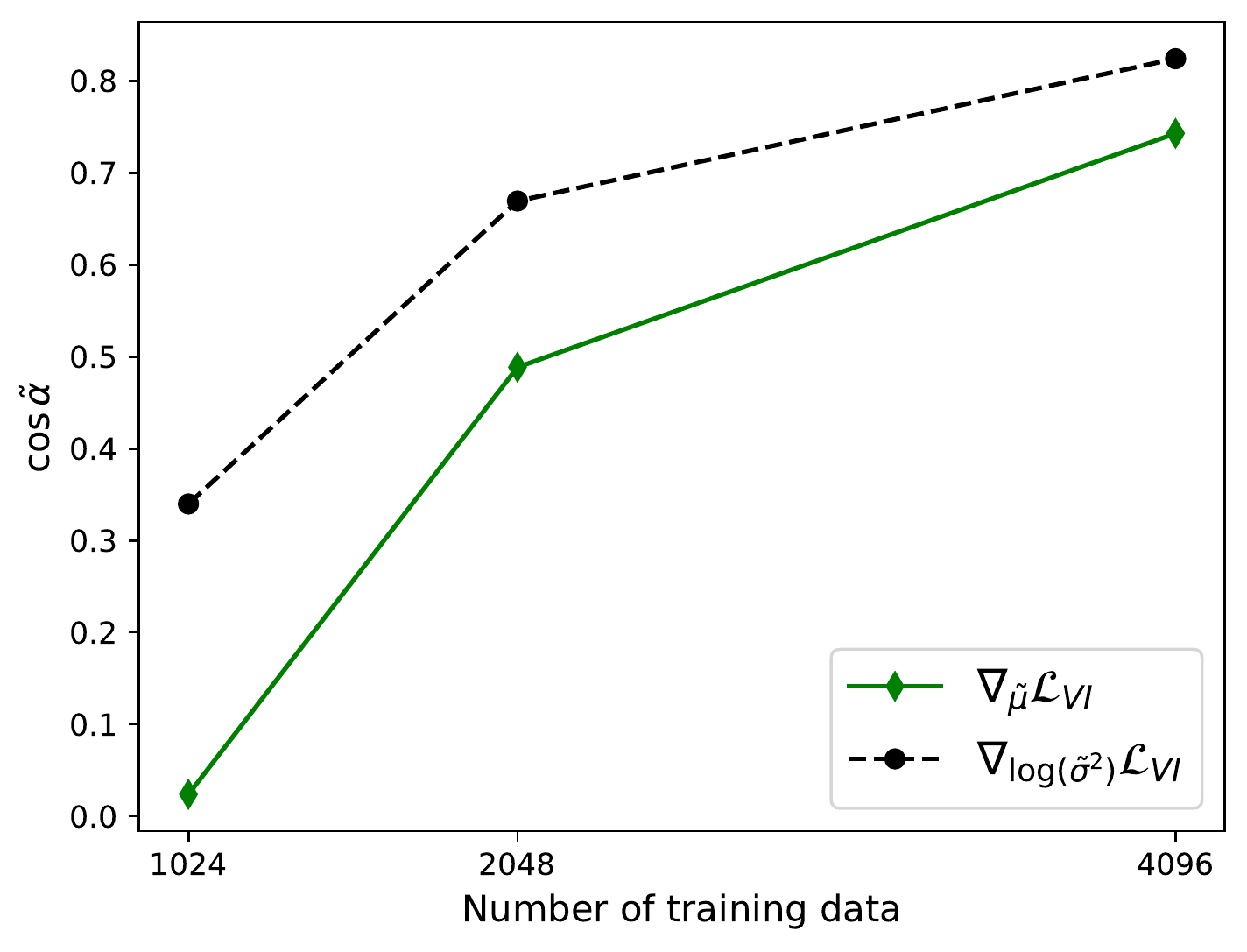}
    \caption{The computed $\cos\tilde{\alpha}$ with different surrogates in the binary channelized field  case. The surrogates are trained using $1024$, $2048$, and $4096$ training data, respectively. The green solid line shows the results of $\nabla_{\tilde{\bmu}}\mathcal{L}_{VI}$. The black dashed line shows the results of $\nabla_{\log(\tilde{\bsigma}^2)} \mathcal{L}_{VI}$.}
    \label{fig:chan_cos}
\end{figure}

\begin{figure}[h!]
    \centerline{
        \begin{tabular}{cc}
    \includegraphics[scale=0.35]{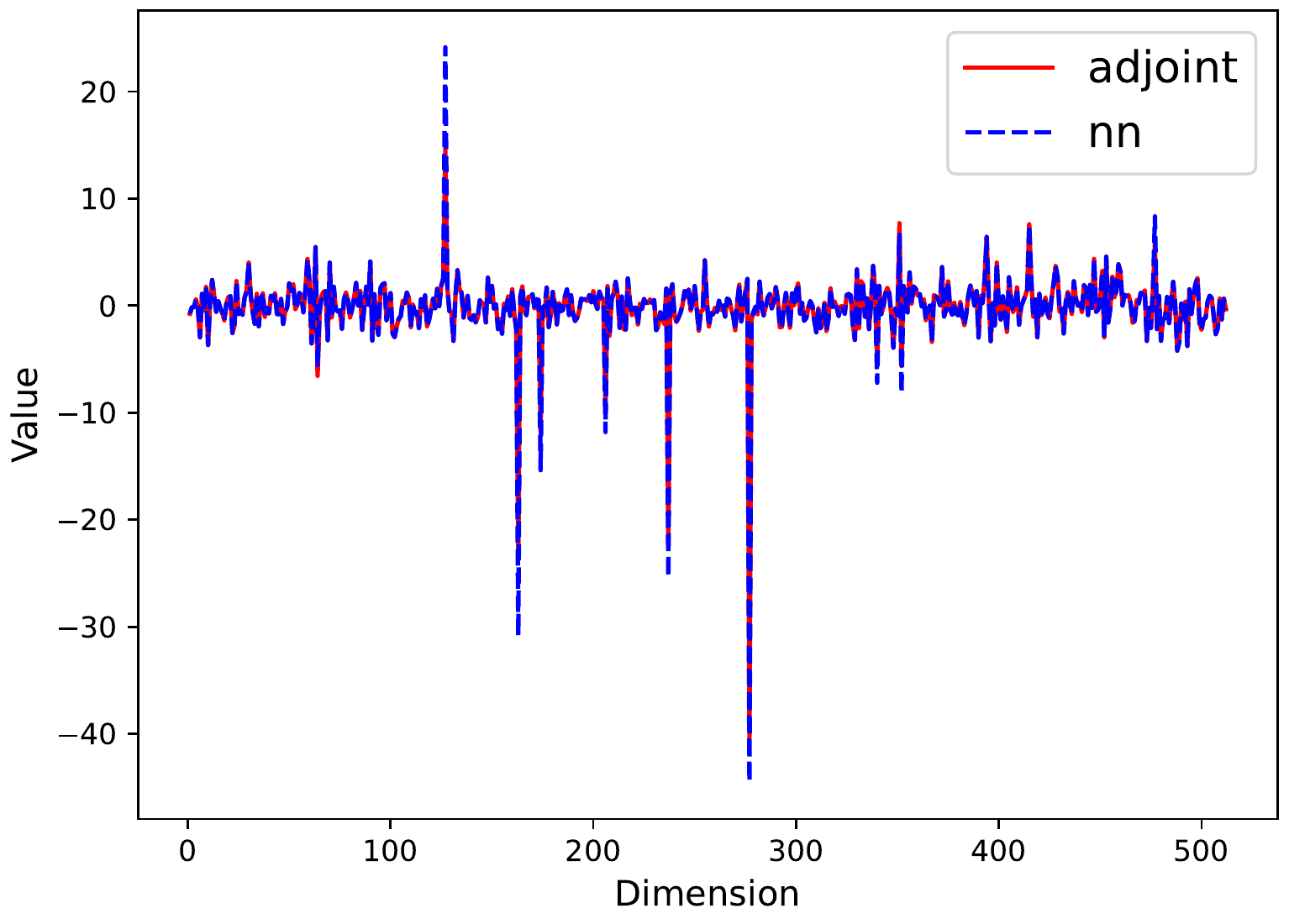}
    & 
    \includegraphics[scale=0.35]{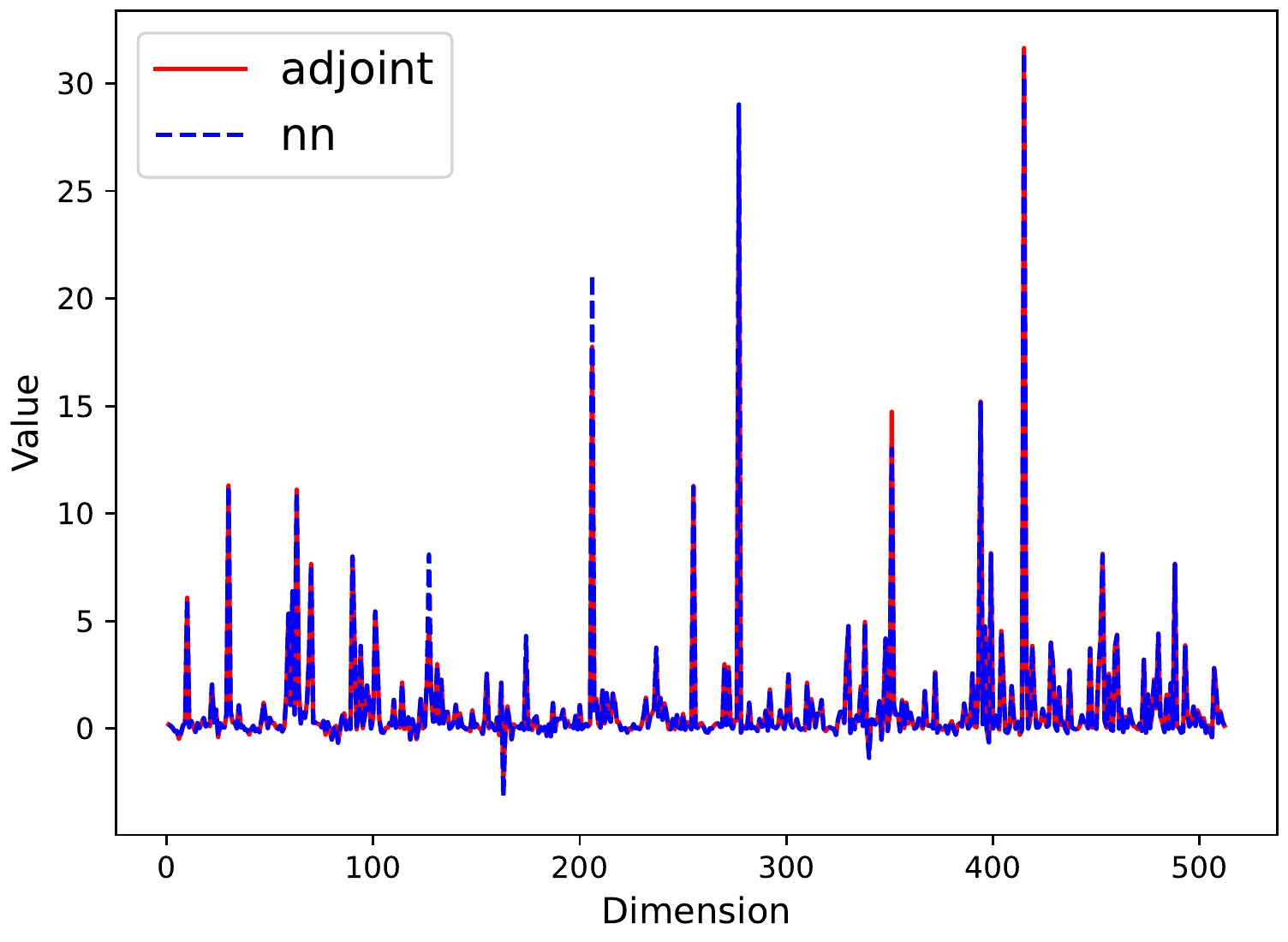}
    \\
    (a) $\nabla_{\tilde{\bmu}}\mathcal{L}_{VI}$& 
    (b) $\nabla_{\log(\tilde{\bsigma}^2)} \mathcal{L}_{VI}$\\   
     
\end{tabular}}
\caption{An example of the gradient computation by the adjoint method and neural network in the binary channelized field case. Given the $\tilde{\bmu}$ and $\log(\tilde{\bsigma}^2)$, the vector (a)$\nabla_{\tilde{\bmu}}\mathcal{L}_{VI}$ (b)$\nabla_{\log(\tilde{\bsigma}^2)} \mathcal{L}_{VI}$ are computed by the workflow in Fig.~\ref{fig:AD_ELBO}, where the $\frac{\partial \mathcal{L}_{VI}}{\partial \bk}$ in the workflow are computed by the neural networks (blue dashed line) and the adjoint method (red solid line), respectively.}
\label{fig:chan_cos_example}
\end{figure}

We also use Eq.~\eqref{eq:cos} to evaluate the gradient approximation. For the test example in Fig.~\ref{fig:chan_truth} with $5\%$ independent Gaussian random noise on the $64$ observations, we can compute the corresponding gradient for any given parameters.  With $N_g=1000$ pairs of $\tilde{\bmu}$ and $\log(\tilde{\bsigma}^2)$ sampled from Gaussian distribution $\mathcal{N}(\bm{0},\bm{I})$, the computed $\cos\tilde{\alpha}$ is given in Fig.~\ref{fig:chan_cos} for three training data scenarios. The gradient approximation is worse compared to the GRF case. This is mainly because the gradient approximation is much more sensitive to the discontinuous log-permeability field, while the GRF is much smoother. Note that two computed $\cos\tilde{\alpha}$ with respect to $\nabla_{\tilde{\bmu}}\mathcal{L}_{VI}$ and $\nabla_{\log(\tilde{\bsigma}^2)} \mathcal{L}_{VI}$ using $4096$ training data are close to $0.8$, we can adopt this learned surrogate model for the VI-DGP method to replace the adjoint method.  Fig.~\ref{fig:chan_cos_example} gives an example of gradient computed by the surrogate model trained with $4096$ training data and the adjoint method, which are similar to those obtained in the GRF case even though the dimension is $512$.

\subsubsection{Bayesian inversion results}
For the non-Gaussian parameter estimation, previous sampling methods have employed advanced strategies such as multiple chains~\cite{laloy2017inversion}, multiscale representation~\cite{xia2022bayesian}, and ensemble-based data assimilation methods~\cite{mo2020integration} to obtain proper results and reduce computational cost. Using the VI methods for non-Gaussian parameter estimation is still a problem to be explored. Typically, the VI methods are restricted by the analytical variational distribution, which leads to large approximation errors for complex non-Gaussian parameter estimation. In the VI-DGP method, we only need to estimate the posterior distribution of the latent variable. Based on this example, we can verify whether the VI-DGP method can still recover the non-Gaussian parameter with uncertainty and get rid of the curse of dimensionality, although the latent variable is high-dimensional. As discussed in the GRF case, we test the accuracy and efficiency of the VI-NN and VI-adjoint methods compared to the referenced MCMC-NN and MCMC-FEM methods. Also, we test the influence of the sampling number $M_s$ and the noise level for the estimation.

\emph{Comparisons.}
The four methods are applied to the test example in Fig.~\ref{fig:chan_truth} with $5\%$ independent Gaussian random noise. We choose an optimization iteration $N_{opt} = 8000$ for the VI method in Algorithm~\ref{alg:VIM}. The other inputs are the same as the GRF case, such as the posterior samples $N_s=10000$, sampling number $M_s =1$, the SGD optimizer with learning rate $\eta_{\tilde{\bmu}}=\eta_{\tilde{\bsigma}} = 0.0008$, and zero initial states. For MCMC, we use a long Markov chain with a length of $300000$ to guarantee convergence. The last $10000$ states are used as posterior samples. The inference time of the four methods is given in Table~\ref{Table:Channel cost}. From this example, one can find significant differences in computational time between the VI method and MCMC.

\begin{table}[h!]
	\caption{Computational cost of estimation with different methods in the binary channelized case.} 
	\centering	
	\begin{tabular}{ccccc}  
		\hline
	Methods	&VI-NN & VI-adjoint   & MCMC-NN &  MCMC-FEM\\\hline
	 	Iterations & $8000$ & $8000$ & $300000$ & $300000$  \\ 
	 	Inference time (s) & $\bm{147}$ & $3410$ & $2094$ &$72078$  \\
		\hline
	\end{tabular}
	\label{Table:Channel cost}
\end{table}

\begin{figure}[h!]
    \centering
    \includegraphics[width=4.5in]{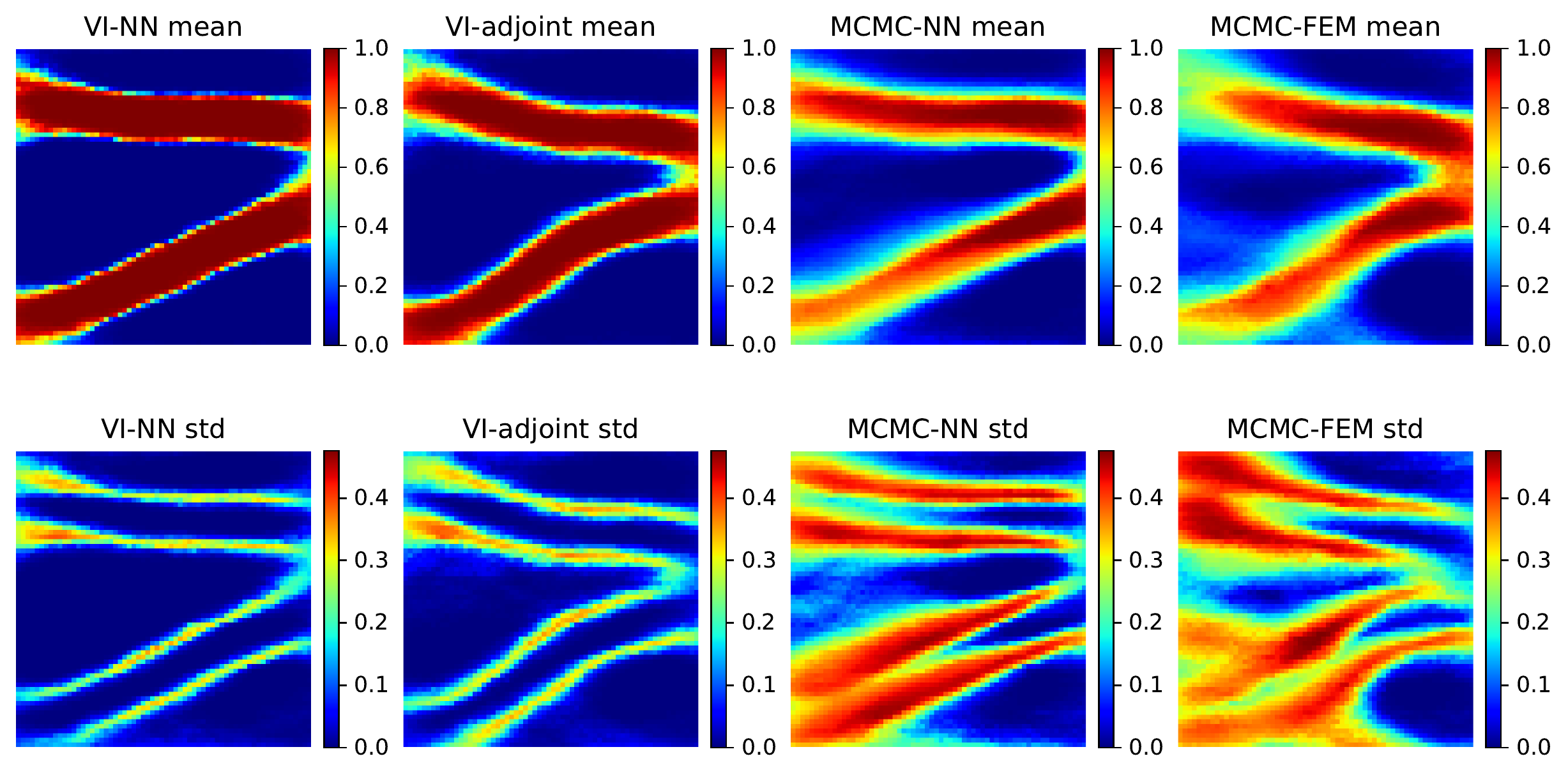}
    \caption{The posterior estimation results for the binary channelized field with different methods. The first row shows the estimated mean of the log-permeability field, and the second row gives the corresponding standard deviation (std).}
    \label{fig:chan_posterior}
\end{figure}

Fig.~\ref{fig:chan_posterior} shows the estimated results by the four methods. It is obvious that the uncertainty of MCMC results is higher than those of the VI methods, and their posterior samples are more diverse. The estimated mean results using the VI method are much better than those estimated by the MCMC method, although the MCMC method can still achieve relatively valid estimations in such a high-dimensional problem. The estimated results also illustrate that the DGP can capture channelized features and generate similar realizations for non-Gaussian parameters, which helps inference acquire appropriate results for both the VI and MCMC methods. Note that even though the sampling number $M_s$ is $1$, the desired accuracy and efficiency can still be realized in the non-Gaussian case.

\begin{figure}[h!]
    \centerline{
        \begin{tabular}{cc}
    \includegraphics[width=2.2in]{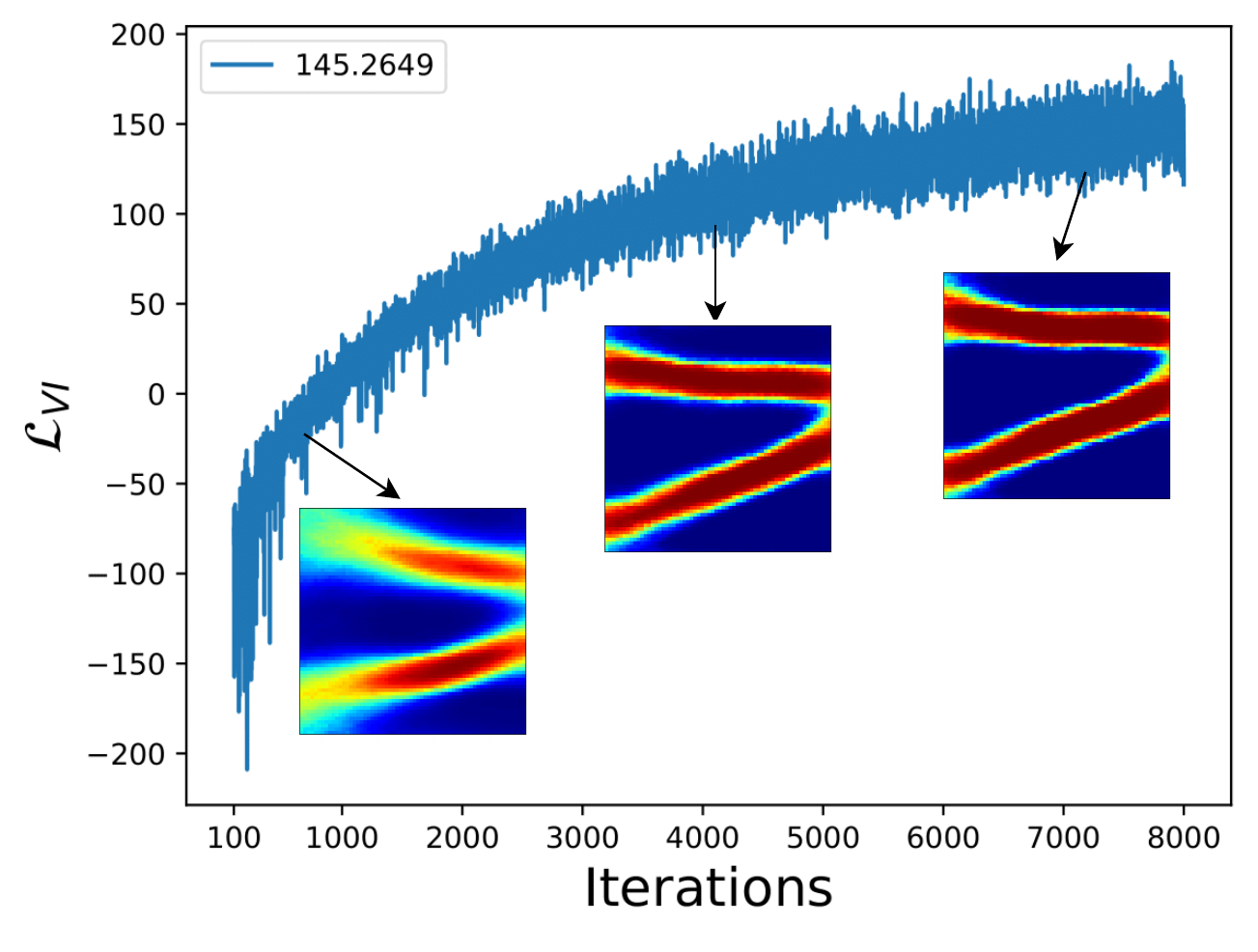}
    & 
   \includegraphics[width=2.2in]{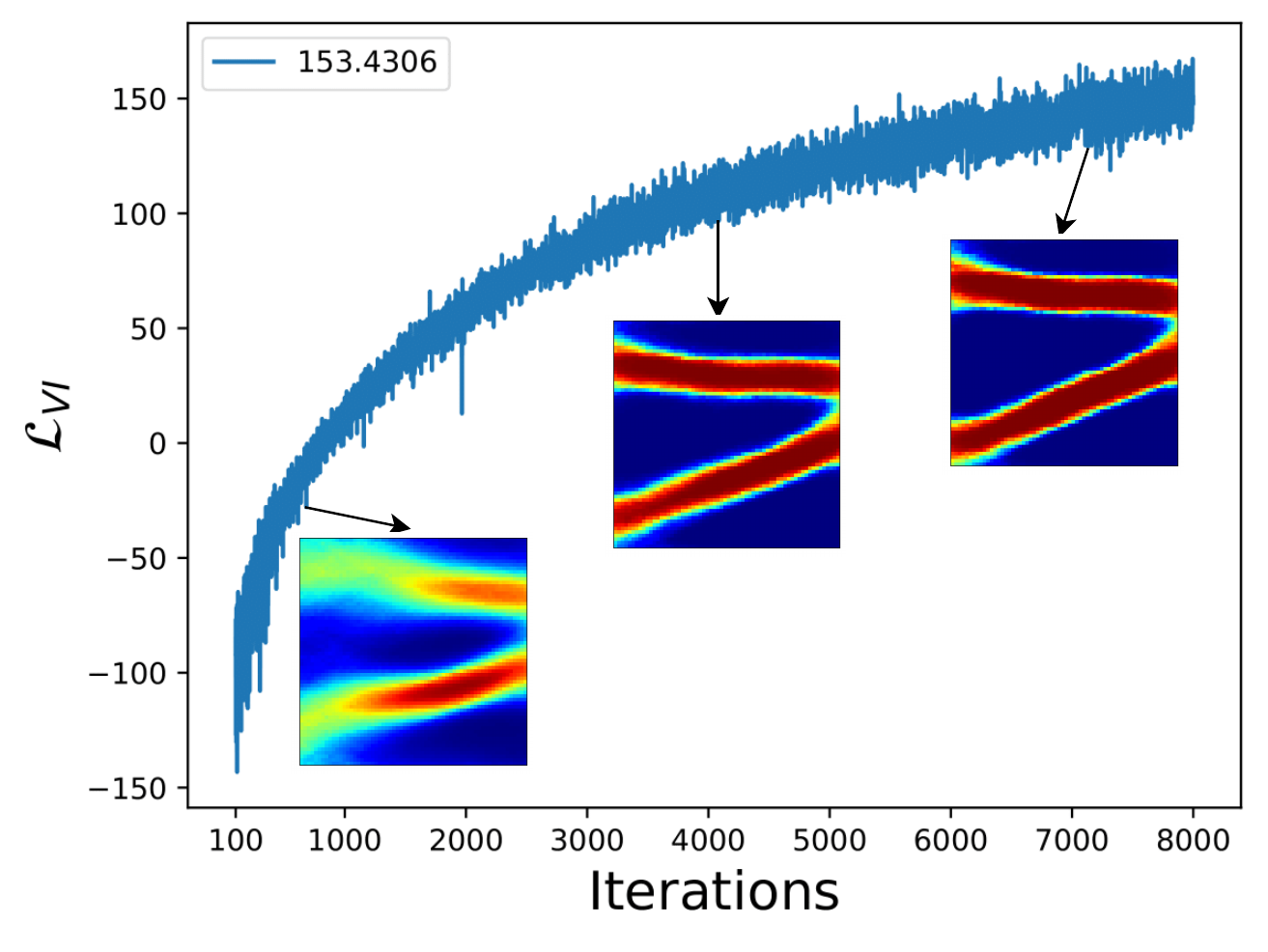}
    \\
    (a) $M_s=1$& 
    (b) $M_s=3$\\   
     \includegraphics[width=2.2in]{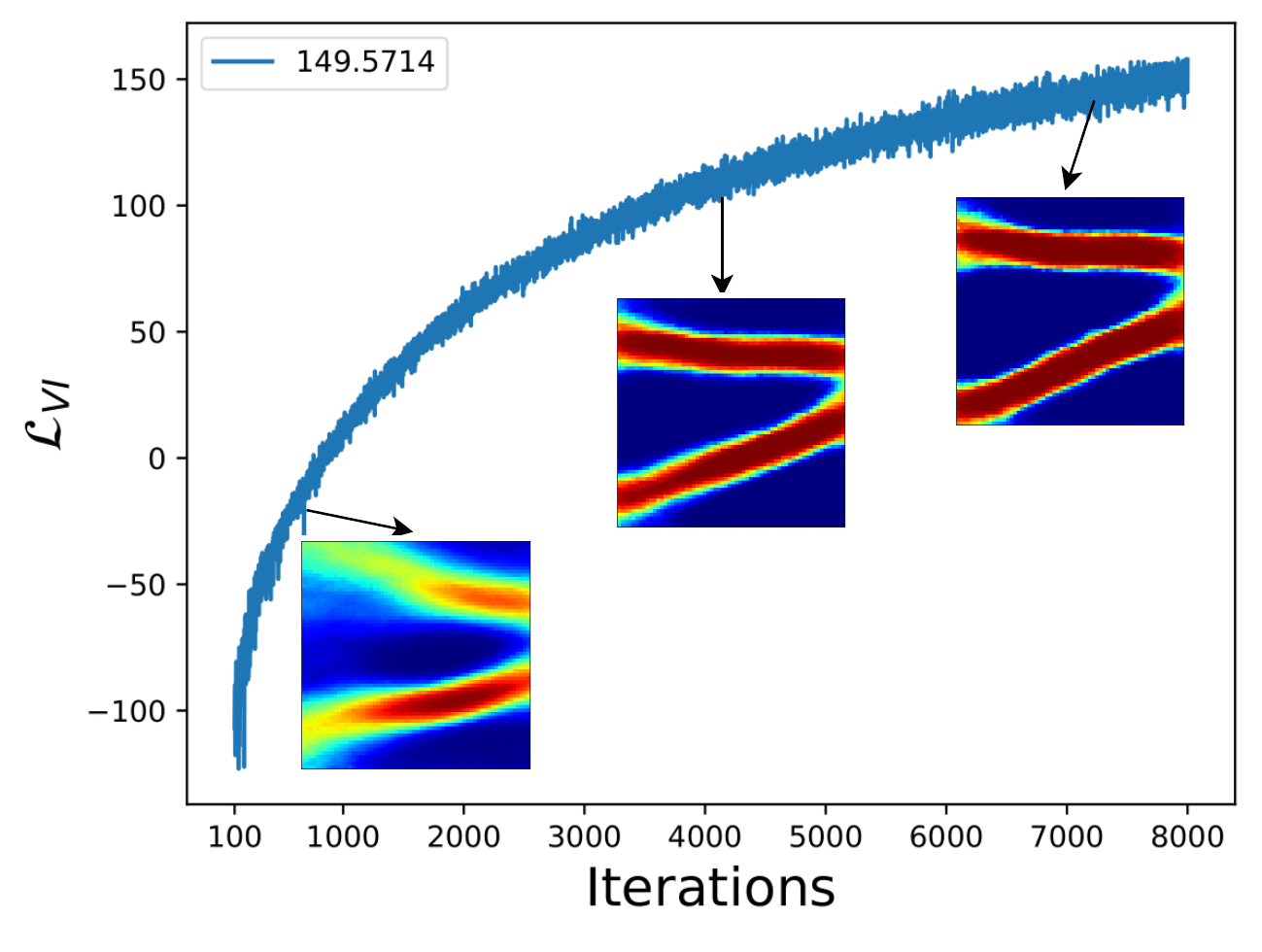} &
    \includegraphics[width=2.2in]{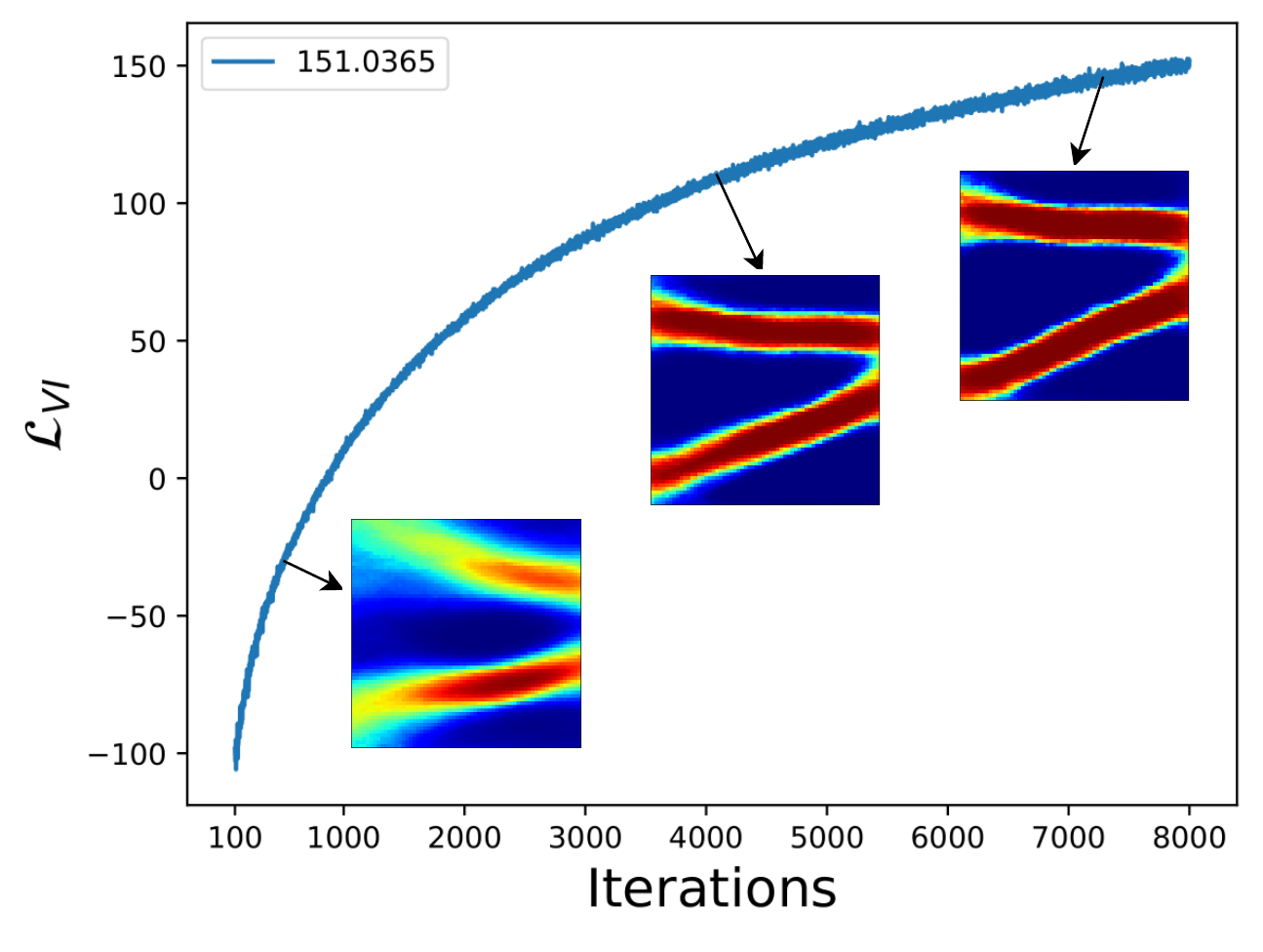}\\
    (c) $M_s=10$& 
    (d) $M_s=100$\\  
\end{tabular}}
\caption{The convergence of the variational lower bound $\mathcal{L}_{VI}$ with varying sample numbers $M_s$ in stochastic optimization. The three log-permeability fields below the black arrows are the estimated mean at the $500$-th, $4000$-th, and $7000$-the iteration, respectively. }
\label{fig:chan_elbo}
\end{figure}

\emph{Effect of the sampling number $M_s$.}
Here, we also test the influence of the sampling number for convergence in the non-Gaussian case. Fig.~\ref{fig:chan_elbo} shows the convergence of the variational lower bound $\mathcal{L}_{VI}$ and the estimated mean results at certain iterations. The four results provide similar convergence trends and estimated mean at those iterations. Even at the $500$-th iteration, the optimization algorithm can capture the important features (channel locations) of the underlying true log-permeability. The main difference between the four experiments is the stability of convergence, where a larger sampling number can give more stable convergence (like $M_s=100$). Correspondingly, the increased computational cost is a significant burden for applications.

\begin{figure}[h!]
    \centerline{
        \begin{tabular}{cc}
    \includegraphics[width=2.2in]{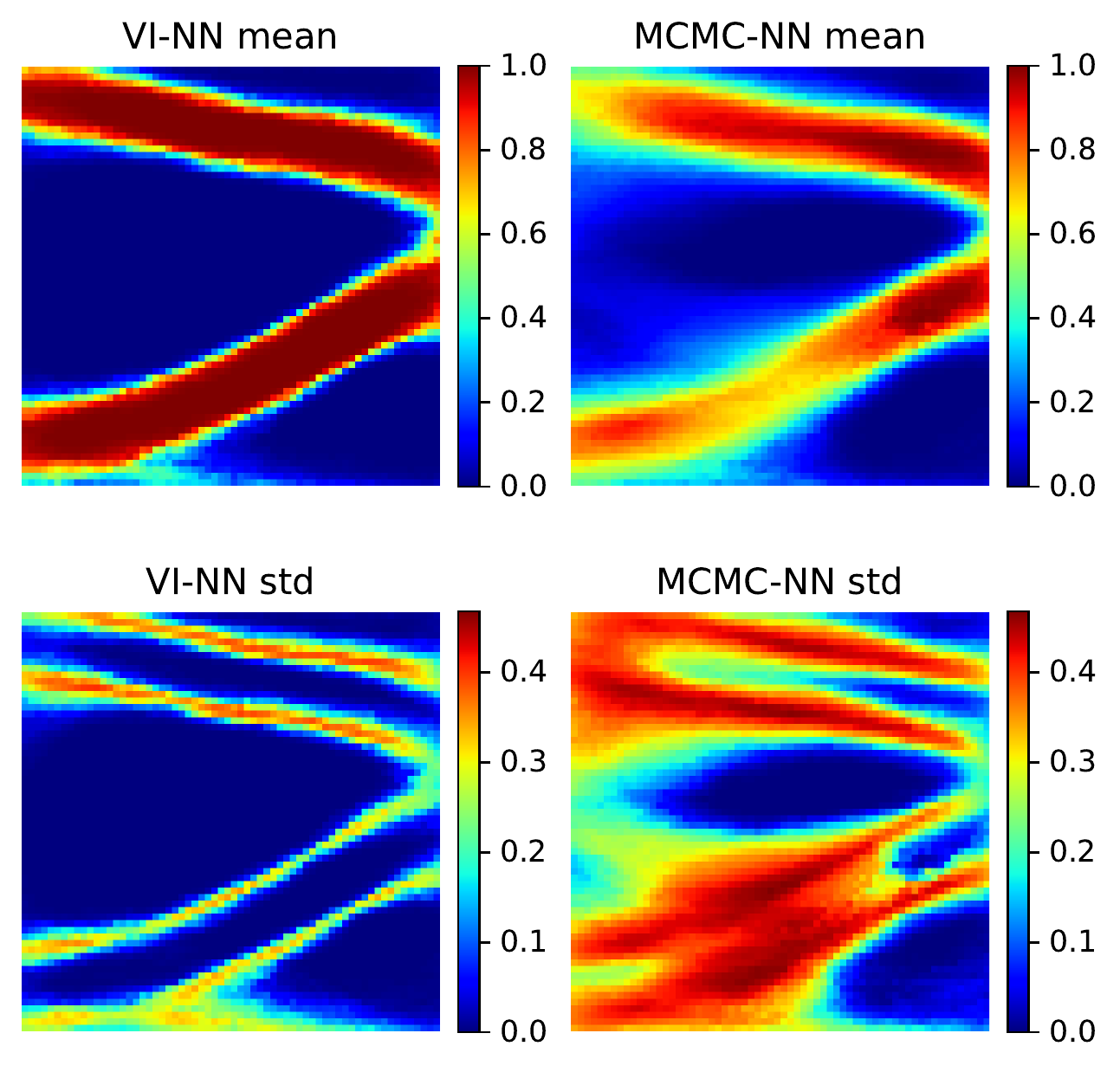}
    & 
    \includegraphics[width=2.2in]{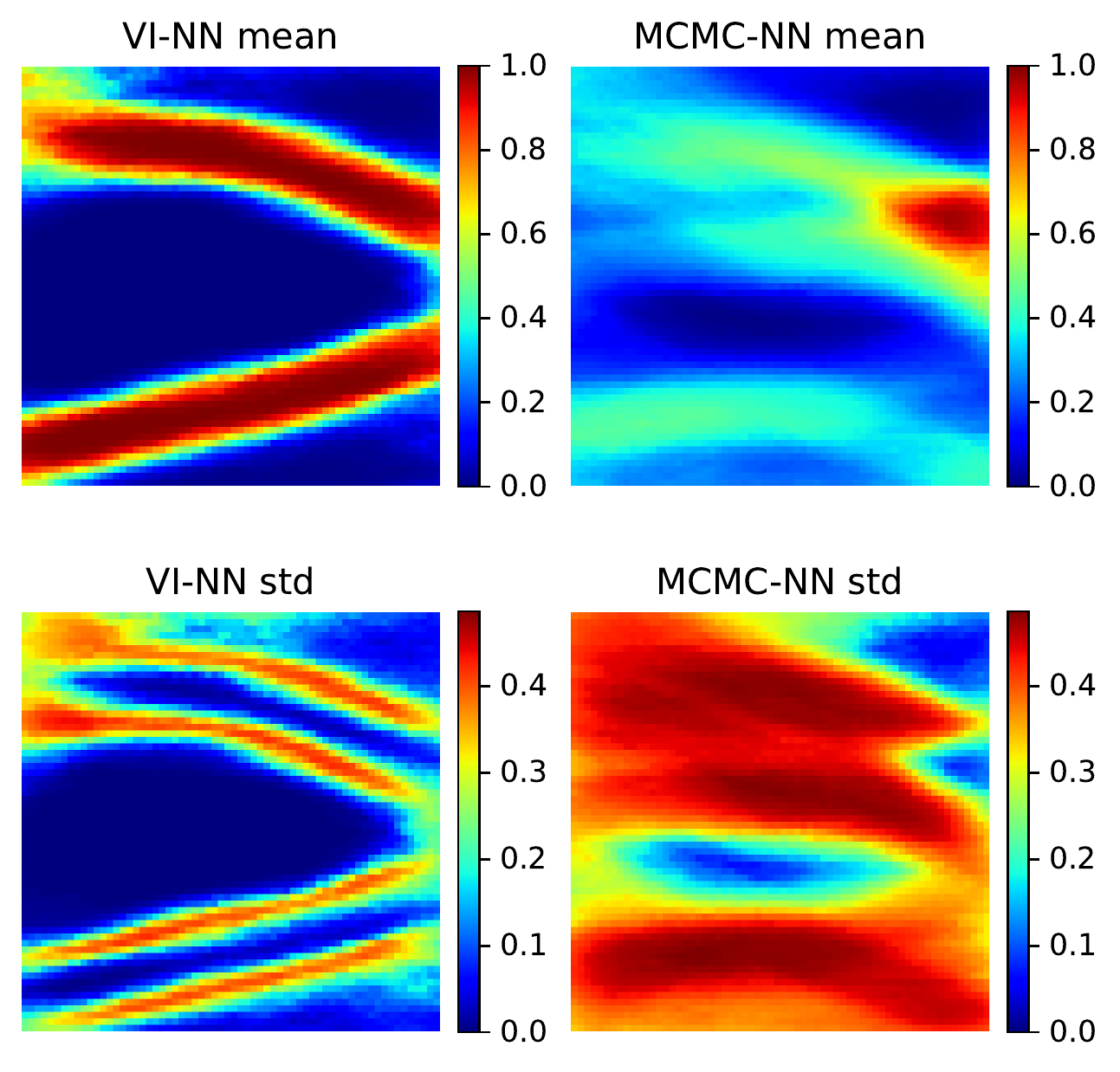}
    \\
    (a) $7\%$ noise& 
    (b) $10\%$ noise\\   
      
\end{tabular}}
\caption{The posterior estimation results in binary channelized field case using VI-NN and MCMC-NN method under (a) $7\%$ noise and (b) $10\%$ noise observations. The first row shows the estimated mean of the log-permeability field, and the second row gives the corresponding standard deviation (std).}
\label{fig:chan_noise}
\end{figure}

\emph{Effect of the observation noise.}
The high noise level, together with the discontinuous parameter, may pose challenges for the estimation with the VI-DGP method. We test two additional examples with $7\%$ and $10\%$ independent Gaussian random noise. The estimated results using VI-NN and MCMC-NN are shown in Fig.~\ref{fig:chan_noise}. The estimated results using the VI-DGP method can obtain reasonable mean results with low uncertainty, although they are worse than the results under $5\%$ noise in Fig.~\ref{fig:chan_posterior}. In contrast, the estimated results using the MCMC method are much worse, especially in the $10\%$ noise case. These results demonstrate the good performance of the VI-DGP method for non-Gaussian parameters.

\section{Conclusions}\label{sec:Conclusions}
Performing efficient inference for probabilistic models is a fundamental problem in machine learning and Bayesian statistics. For BIPs, efficiency and accuracy are the primary influences of their popularity in science and engineering. In this work, we propose a novel method for solving high-dimensional inverse problems applied in spatially-varying parameter estimation. Unlike sampling methods, VI methods typically approximate the posterior distribution through optimization, which favors scalability and acceleration using GPUs. However, their limited choice of variational distribution can restrict the capacity to approximate complex distributions. To overcome this limitation, we propose the VI-DGP method, which exploits the generation ability of the DGM in prior modeling and posterior approximation. Our data-driven prior model can incorporate various prior information, and the obtained latent variable can be leveraged for dimension reduction and posterior approximation. Additionally, we use physics-constrained neural networks and their inherent automatic differentiation to avoid the need for the adjoint method and make our method easy to implement and transfer to various problems. Our numerical experiments show that the proposed VI-DGP method outperforms the referenced method in terms of both efficiency and accuracy.

Although the proposed VI-DGP method provides a very general and flexible framework for BIPs, there are still many issues that need to be investigated and discussed. The use of neural network surrogates is promising in PDE-constrained optimization problems. However, further theoretical analysis and comparison with the adjoint method are needed for scientific computing tasks. Furthermore, advanced VI methods and auxiliary latent variables can also be employed to improve the flexibility of approximations and inference capacity.


\begin{appendices}

\section{The network architectures for the encoder and decoder in VAE}\label{app:nn_VAE}
In this work, we use fully-connected neural networks as the encoder and decoder for both Gaussian and channel cases. Table~\ref{VAE_nn} illustrates the implemented neural networks for the encoder and decoder.
For the decoder, we use \verb|ReLU| and \verb|Sigmoid| as the activation function for the Gaussian and channel cases, respectively. Additionally, for the channel case, we apply an extra \verb|Sigmoid| activation function for the last layer of the decoder model, which ensures that the output values are within the interval $[0,1]$. $h$ denotes the number of neurons in the encoder's hidden layer, which will also define the dimensionality of the latent variable $\bz$. We set $h$ to $256$ and $512$ for the Gaussian and channel cases, respectively.
\begin{table}[h!]
\caption{The employed network architectures of the VAE model. $\text{Linear}  (H_{in}, H_{out})$ denotes the linear operator, where $H_{in}$ and $H_{out}$ are the parameter size of the input and output, respectively. } 
\centering
\begin{tabular}{|c|c|c|c|}
\hline  \multicolumn{2}{|c|}{\text { Encoder }} & \multicolumn{2}{c|}{\text { Decoder }}  \\
\hline  \multicolumn{2}{|c|}{\text { Input: $\bk$ }}  &  \multicolumn{2}{c|}{\text { Input: $\bz$ }} \\
\hline  \multicolumn{2}{|c|}{\text { Linear $(4096, h)$ }}  &  \multicolumn{2}{c|}{\text { Linear $(h, 4096)$ }} \\
\hline  \multicolumn{2}{|c|}{\text { ReLU  } }  &  \multicolumn{2}{c|}{\text { ReLU/Sigmoid }} \\
\hline  \multicolumn{2}{|c|}{\text { Linear $(h, h)$ }}  &  \multicolumn{2}{c|}{\text { Linear (4096, 4096) }} \\
\hline  \multicolumn{2}{|c|}{\text { ReLU } }  &  \multicolumn{2}{c|}{\text { ReLU/Sigmoid }} \\
\hline \text { Linear $$(h, h)$$ }  & \text { Linear $(h, h)$} & \multicolumn{2}{c|}{\text { Linear (4096, 4096) }}\\
\hline \text { ReLU }   &  \text { ReLU }   &  \multicolumn{2}{c|}{\text { ReLU/Sigmoid }} \\
\hline \text { Linear $(h, h)$ }  & \text { Linear $(h, h)$ } &\multicolumn{2}{c|}{\text { Linear (4096, 4096) }}\\
\hline  \text { output: $\bmu$ } & \text { output: $\log(\bsigma)$ }  &  \multicolumn{2}{c|}{\text { output: $\bk$ }} \\
\hline

\end{tabular}
\label{VAE_nn}
\end{table}

\section{The network architectures for the physics-constrained surrogate model}\label{app:nn_pcs}
We can rewrite the loss function in discretization form for the given PDEs in Eq.~\eqref{eq:darcy} and Eq.~\eqref{eq:darcy_boundary}. The PDEs loss and boundary loss in Eq.~\eqref{eq:empirical_loss} can be written as
\begin{equation}
\begin{aligned}
J_{\text{pde}}(u(\bx, \bk ; \Theta)) =& \frac{1}{n_s n_p}\sum_{j=1}^{n_s}\sum_{i=1}^{n_p} ( 
  \|  \nabla \cdot \bm{v}(\bx_{\mathcal{D}}^{(i)}) - f(\bx_{\mathcal{D}}^{(i)}) \|^2\\
 +&\| \bm{v}(\bx_{\mathcal{D}}^{(i)}) + \exp(\bk^{(j)}(\bx_{\mathcal{D}}^{(i)})) \odot \nabla p(\bx_{\mathcal{D}}^{(i)})\|^2 ),
\\
J_{\text{b}}(u(\bx, \bk ; \Theta)) =& \frac{1}{n_{bl}}\sum_{i=1}^{n_{bl}}\| p(\bx_{\mathcal{D}_l}^{(i)}) -1\|^2
+\frac{1}{n_{br}}\sum_{i=1}^{n_{br}}\| p(\bx_{\mathcal{D}_r}^{(i)}) \|^2\\
+&\frac{1}{n_{bt}}\sum_{i=1}^{n_{bt}}\| \bm{v}(\bx_{\mathcal{D}_t}^{(i)}) \|^2
+\frac{1}{n_{bb}}\sum_{i=1}^{n_{bb}}\| \bm{v}(\bx_{\mathcal{D}_b}^{(i)}) \|^2,\\
\label{eq:darcy_loss}
\end{aligned}
\end{equation}
respectively, where $n_b$ boundary samples include $n_{bl}$ samples of left boundary $\mathcal{D}_l$, $n_{br}$ samples of right boundary $\mathcal{D}_r$, $n_{bt}$ samples of top boundary $\mathcal{D}_t$, and $n_{bb}$ samples of bottom boundary  $\mathcal{D}_b$.

The network architectures applied in this paper are based on previous works~\cite{zhu2018bayesian,zhu2019physics}. These works perform greatly in uncertainty quantification tasks for the flow in heterogeneous media. The main architectures are shown in Table~\ref{pcn_nn}. The number of dense layers in the three dense blocks is $6, 8, 6$, with a growth rate of 16. Each dense layer contains a Conv block (\verb|Batch-ReLU-Conv|). \text {Encoding 1}, \text{Decoding 1}, and \text{Decoding 2} have $2, 2, 3$  Conv blocks, respectively. The \verb|nearest| mode is used for the \verb|upsampling| operator in the decoding layers.
  
\begin{table}[h!]
\caption{The network architectures for the physics-constrained surrogate in this paper. } 
\centering
\begin{tabular}{|c|c|}
\hline \text {Networks} &  \text {Feature maps} \\
\hline \text {Input} &  \text {$1\times64\times64$}  \\
\hline \text {Conv layer} &  \text {$48\times32\times32$}  \\
\hline \text {Dense Block} &  \text {$144\times32\times32$}  \\
\hline \text {Encoding 1} &  \text {$76\times16\times16$}  \\
\hline \text {Dense Block} &  \text {$200\times16\times16$}  \\ 
\hline \text {Decoding 1} &  \text {$100\times32\times32$}  \\ 
\hline \text {Dense Block} &  \text {$196\times32\times32$}  \\
\hline \text {Decoding 2} &  \text {$3\times64\times64$}  \\
\hline \text {Output} &  \text {$3\times64\times64$}  \\ 
\hline
\end{tabular}
\label{pcn_nn}
\end{table}

\section{The pCN algorithm for MCMC simulation}\label{app:pcn}
We employ the pCN algorithm to explore the posterior distribution, which is the reference method for the proposed approach. The details are shown in the Algorithm~\ref{agl:pcn}, where the forward model $\mathcal{F}(\cdot)$ can be either the learned neural network surrogate or the finite element method. These correspond to MCMC-NN and MCMC-FEM in the experiments, respectively.
\begin{algorithm}[h!]
	\caption{pCN algorithm with the DGP}
	\label{agl:pcn}
	\begin{algorithmic}[1]
    	\Require  the likelihood $\pi(\bd | \bz)$, the forward model $\mathcal{F}(\cdot)$, chain length $N_{ite}$, generative model $\mathcal{G}_{\btheta^{\star}}( \bz)$, burn-in length $N_b$, $\beta = 0.15$
		\State Initialize $\bz^{(1)},  \bz^{(1)} \sim \mathcal{N}(\bm{0},\bm{I})$
		\For {$j = 1:N_{ite}$}
		\State Draw $\bz^{\prime}$ via
		$$\bz^{\prime} = \sqrt{1-\beta^{2}} \bz^{(i)}+\beta \hat{\bxi}, \quad \text { where } \hat{\bxi} \sim \mathcal{N}\left(0, \boldsymbol{I} \right)$$
		\State Compute the likelihood function by solving the forward model $\mathcal{F}(\bk^{\prime})$, where $\bk^{\prime} = \mathcal{G}_{\btheta^{\star}}(\bz^{\prime})$
		\State Compute the acceptance ratio
		$$\alpha = min\left(1,\frac{\pi(\bd | \bz^{\prime})}{\pi(\bd | \bz^{(j)})}\right)$$
		\State Draw $\rho$ from the uniform distribution $\mathcal{U}[0,1]$
 		\If {$\rho < \alpha$}
		\State $\text { Let } \bz^{(j+1)}=\bz^{\prime}, \bk^{(j+1)}=\bk^{\prime}$
		\Else
		\State {$\text { Let }\bz^{(j+1)}=\bz^{j}, \bk^{(j+1)}=\bk^{(j)}$}
		\EndIf
		\EndFor
		\Ensure posterior samples $\{\bk^{(i)}\}_{i=N_b}^{N_{ite}}$
	\end{algorithmic}
	
\end{algorithm}

\end{appendices}


%
%

\bibliographystyle{spmpsci}      

\bibliography{references}

\begin{thebibliography}{10}
\providecommand{\url}[1]{{#1}}
\providecommand{\urlprefix}{URL }
\expandafter\ifx\csname urlstyle\endcsname\relax
  \providecommand{\doi}[1]{DOI~\discretionary{}{}{}#1}\else
  \providecommand{\doi}{DOI~\discretionary{}{}{}\begingroup
  \urlstyle{rm}\Url}\fi

\bibitem{barajas2019approximate}
Barajas-Solano, D.A., Tartakovsky, A.M.: Approximate bayesian model inversion
  for pdes with heterogeneous and state-dependent coefficients.
\newblock Journal of Computational Physics \textbf{395}, 247--262 (2019)

\bibitem{bilionis2013multi}
Bilionis, I., Zabaras, N., Konomi, B.A., Lin, G.: Multi-output separable
  gaussian process: Towards an efficient, fully bayesian paradigm for
  uncertainty quantification.
\newblock Journal of Computational Physics \textbf{241}, 212--239 (2013)

\bibitem{blei2017variational}
Blei, D.M., Kucukelbir, A., McAuliffe, J.D.: Variational inference: A review
  for statisticians.
\newblock Journal of the American Statistical Association \textbf{112}(518),
  859--877 (2017)

\bibitem{bora2017compressed}
Bora, A., Jalal, A., Price, E., Dimakis, A.G.: Compressed sensing using
  generative models.
\newblock In: International Conference on Machine Learning, pp. 537--546. PMLR
  (2017)

\bibitem{bui2014solving}
Bui-Thanh, T., Girolami, M.: Solving large-scale pde-constrained bayesian
  inverse problems with riemann manifold hamiltonian monte carlo.
\newblock Inverse Problems \textbf{30}(11), 114014 (2014)

\bibitem{chen2021stein}
Chen, P., Ghattas, O.: Stein variational reduced basis bayesian inversion.
\newblock SIAM Journal on Scientific Computing \textbf{43}(2), A1163--A1193
  (2021)

\bibitem{cotter2013mcmc}
Cotter, S.L., Roberts, G.O., Stuart, A.M., White, D.: Mcmc methods for
  functions: modifying old algorithms to make them faster.
\newblock Statistical Science \textbf{28}(3), 424--446 (2013)

\bibitem{cui2015data}
Cui, T., Marzouk, Y.M., Willcox, K.E.: Data-driven model reduction for the
  bayesian solution of inverse problems.
\newblock International Journal for Numerical Methods in Engineering
  \textbf{102}(5), 966--990 (2015)

\bibitem{engl1996regularization}
Engl, H.W., Hanke, M., Neubauer, A.: Regularization of inverse problems, vol.
  375.
\newblock Springer Science \& Business Media (1996)

\bibitem{fan2019solving}
Fan, Y., Ying, L.: Solving inverse wave scattering with deep learning.
\newblock arXiv preprint arXiv:1911.13202  (2019)

\bibitem{geneva2020modeling}
Geneva, N., Zabaras, N.: Modeling the dynamics of pde systems with
  physics-constrained deep auto-regressive networks.
\newblock Journal of Computational Physics \textbf{403}, 109056 (2020)

\bibitem{goodfellow2014generative}
Goodfellow, I., Pouget-Abadie, J., Mirza, M., Xu, B., Warde-Farley, D., Ozair,
  S., Courville, A., Bengio, Y.: Generative adversarial nets.
\newblock Advances in Neural Information Processing Systems \textbf{27} (2014)

\bibitem{guha2015variational}
Guha, N., Wu, X., Efendiev, Y., Jin, B., Mallick, B.K.: A variational bayesian
  approach for inverse problems with skew-t error distributions.
\newblock Journal of Computational Physics \textbf{301}, 377--393 (2015)

\bibitem{hairer2014spectral}
Hairer, M., Stuart, A.M., Vollmer, S.J.: Spectral gaps for a
  metropolis-hastings algorithm in infinite dimensions.
\newblock The Annals of Applied Probability \textbf{24}(6), 2455--2490 (2014)

\bibitem{jalal2021robust}
Jalal, A., Arvinte, M., Daras, G., Price, E., Dimakis, A.G., Tamir, J.: Robust
  compressed sensing mri with deep generative priors.
\newblock Advances in Neural Information Processing Systems \textbf{34},
  14938--14954 (2021)

\bibitem{jia2021variational}
Jia, J., Zhao, Q., Xu, Z., Meng, D., Leung, Y.: Variational bayes' method for
  functions with applications to some inverse problems.
\newblock SIAM Journal on Scientific Computing \textbf{43}(1), A355--A383
  (2021)

\bibitem{kaipio2006statistical}
Kaipio, J., Somersalo, E.: Statistical and computational inverse problems, vol.
  160.
\newblock Springer Science \& Business Media (2006)

\bibitem{khoo2019switchnet}
Khoo, Y., Ying, L.: Switchnet: a neural network model for forward and inverse
  scattering problems.
\newblock SIAM Journal on Scientific Computing \textbf{41}(5), A3182--A3201
  (2019)

\bibitem{kingma2014adam}
Kingma, D.P., Ba, J.: Adam: A method for stochastic optimization.
\newblock arXiv preprint arXiv:1412.6980  (2014)

\bibitem{kingma2013auto}
Kingma, D.P., Welling, M.: Auto-encoding variational bayes.
\newblock arXiv preprint arXiv:1312.6114  (2013)

\bibitem{laloy2018training}
Laloy, E., H{\'e}rault, R., Jacques, D., Linde, N.: Training-image based
  geostatistical inversion using a spatial generative adversarial neural
  network.
\newblock Water Resources Research \textbf{54}(1), 381--406 (2018)

\bibitem{laloy2017inversion}
Laloy, E., H{\'e}rault, R., Lee, J., Jacques, D., Linde, N.: Inversion using a
  new low-dimensional representation of complex binary geological media based
  on a deep neural network.
\newblock Advances in Water Resources \textbf{110}, 387--405 (2017)

\bibitem{li2023deep}
Li, S., Xia, Y., Liu, Y., Liao, Q.: A deep domain decomposition method based on
  fourier features.
\newblock Journal of Computational and Applied Mathematics \textbf{423}, 114963
  (2023)

\bibitem{liao2019adaptive}
Liao, Q., Li, J.: An adaptive reduced basis anova method for high-dimensional
  bayesian inverse problems.
\newblock Journal of Computational Physics \textbf{396}, 364--380 (2019)

\bibitem{lu2021learning}
Lu, L., Jin, P., Pang, G., Zhang, Z., Karniadakis, G.E.: Learning nonlinear
  operators via deeponet based on the universal approximation theorem of
  operators.
\newblock Nature Machine Intelligence \textbf{3}(3), 218--229 (2021)

\bibitem{lye2021iterative}
Lye, K.O., Mishra, S., Ray, D., Chandrashekar, P.: Iterative surrogate model
  optimization (ismo): An active learning algorithm for pde constrained
  optimization with deep neural networks.
\newblock Computer Methods in Applied Mechanics and Engineering \textbf{374},
  113575 (2021)

\bibitem{martin2012stochastic}
Martin, J., Wilcox, L.C., Burstedde, C., Ghattas, O.: A stochastic newton mcmc
  method for large-scale statistical inverse problems with application to
  seismic inversion.
\newblock SIAM Journal on Scientific Computing \textbf{34}(3), A1460--A1487
  (2012)

\bibitem{marzouk2007stochastic}
Marzouk, Y.M., Najm, H.N., Rahn, L.A.: Stochastic spectral methods for
  efficient bayesian solution of inverse problems.
\newblock Journal of Computational Physics \textbf{224}(2), 560--586 (2007)

\bibitem{metropolis1953equation}
Metropolis, N., Rosenbluth, A.W., Rosenbluth, M.N., Teller, A.H., Teller, E.:
  Equation of state calculations by fast computing machines.
\newblock The Journal of Chemical Physics \textbf{21}(6), 1087--1092 (1953)

\bibitem{mo2019deep}
Mo, S., Zabaras, N., Shi, X., Wu, J.: Deep autoregressive neural networks for
  high-dimensional inverse problems in groundwater contaminant source
  identification.
\newblock Water Resources Research \textbf{55}(5), 3856--3881 (2019)

\bibitem{mo2020integration}
Mo, S., Zabaras, N., Shi, X., Wu, J.: Integration of adversarial autoencoders
  with residual dense convolutional networks for estimation of non-gaussian
  hydraulic conductivities.
\newblock Water Resources Research \textbf{56}(2), e2019WR026082 (2020)

\bibitem{padmanabha2021solving}
Padmanabha, G.A., Zabaras, N.: Solving inverse problems using conditional
  invertible neural networks.
\newblock Journal of Computational Physics \textbf{433}, 110194 (2021)

\bibitem{patel2022solution}
Patel, D.V., Ray, D., Oberai, A.A.: Solution of physics-based bayesian inverse
  problems with deep generative priors.
\newblock Computer Methods in Applied Mechanics and Engineering \textbf{400},
  115428 (2022)

\bibitem{povala2022variational}
Povala, J., Kazlauskaite, I., Febrianto, E., Cirak, F., Girolami, M.:
  Variational bayesian approximation of inverse problems using sparse precision
  matrices.
\newblock Computer Methods in Applied Mechanics and Engineering \textbf{393},
  114712 (2022)

\bibitem{raissi2019physics}
Raissi, M., Perdikaris, P., Karniadakis, G.E.: Physics-informed neural
  networks: A deep learning framework for solving forward and inverse problems
  involving nonlinear partial differential equations.
\newblock Journal of Computational Physics \textbf{378}, 686--707 (2019)

\bibitem{ranganath2014black}
Ranganath, R., Gerrish, S., Blei, D.: Black box variational inference.
\newblock In: Artificial Intelligence and Statistics, pp. 814--822. PMLR (2014)

\bibitem{rezende2015variational}
Rezende, D., Mohamed, S.: Variational inference with normalizing flows.
\newblock In: International Conference on Machine Learning, pp. 1530--1538.
  PMLR (2015)

\bibitem{robert1999monte}
Robert, C.P., Casella, G., Casella, G.: Monte Carlo statistical methods,
  vol.~2.
\newblock Springer (1999)

\bibitem{roeder2017sticking}
Roeder, G., Wu, Y., Duvenaud, D.K.: Sticking the landing: Simple,
  lower-variance gradient estimators for variational inference.
\newblock Advances in Neural Information Processing Systems \textbf{30} (2017)

\bibitem{stuart2010inverse}
Stuart, A.M.: Inverse problems: a bayesian perspective.
\newblock Acta numerica \textbf{19}, 451--559 (2010)

\bibitem{sun2020surrogate}
Sun, L., Gao, H., Pan, S., Wang, J.X.: Surrogate modeling for fluid flows based
  on physics-constrained deep learning without simulation data.
\newblock Computer Methods in Applied Mechanics and Engineering \textbf{361},
  112732 (2020)

\bibitem{tarantola2005inverse}
Tarantola, A.: Inverse problem theory and methods for model parameter
  estimation, vol.~89.
\newblock SIAM (2005)

\bibitem{tripathy2018deep}
Tripathy, R.K., Bilionis, I.: Deep uq: Learning deep neural network surrogate
  models for high dimensional uncertainty quantification.
\newblock Journal of Computational Physics \textbf{375}, 565--588 (2018)

\bibitem{tsilifis2016computationally}
Tsilifis, P., Bilionis, I., Katsounaros, I., Zabaras, N.: Computationally
  efficient variational approximations for bayesian inverse problems.
\newblock Journal of Verification, Validation and Uncertainty Quantification
  \textbf{1}(3) (2016)

\bibitem{wan2011bayesian}
Wan, J., Zabaras, N.: A bayesian approach to multiscale inverse problems using
  the sequential monte carlo method.
\newblock Inverse Problems \textbf{27}(10), 105004 (2011)

\bibitem{wang2018randomized}
Wang, K., Bui-Thanh, T., Ghattas, O.: A randomized maximum a posteriori method
  for posterior sampling of high dimensional nonlinear bayesian inverse
  problems.
\newblock SIAM Journal on Scientific Computing \textbf{40}(1), A142--A171
  (2018)

\bibitem{wang2020deep}
Wang, L., Chan, Y.C., Ahmed, F., Liu, Z., Zhu, P., Chen, W.: Deep generative
  modeling for mechanistic-based learning and design of metamaterial systems.
\newblock Computer Methods in Applied Mechanics and Engineering \textbf{372},
  113377 (2020)

\bibitem{wang2021fast}
Wang, S., Bhouri, M.A., Perdikaris, P.: Fast pde-constrained optimization via
  self-supervised operator learning.
\newblock arXiv preprint arXiv:2110.13297  (2021)

\bibitem{warner2015stochastic}
Warner, J.E., Aquino, W., Grigoriu, M.D.: Stochastic reduced order models for
  inverse problems under uncertainty.
\newblock Computer Methods in Applied Mechanics and Engineering \textbf{285},
  488--514 (2015)

\bibitem{xia2022bayesian}
Xia, Y., Zabaras, N.: Bayesian multiscale deep generative model for the
  solution of high-dimensional inverse problems.
\newblock Journal of Computational Physics \textbf{455}, 111008 (2022)

\bibitem{xiu2003modeling}
Xiu, D., Karniadakis, G.E.: Modeling uncertainty in flow simulations via
  generalized polynomial chaos.
\newblock Journal of Computational Physics \textbf{187}(1), 137--167 (2003)

\bibitem{zhihang2023domain}
Xu, Z., Xia, Y., Liao, Q.: A domain-decomposed vae method for bayesian inverse
  problems.
\newblock arXiv preprint arXiv:2301.05708  (2023)

\bibitem{YAN2021114087}
Yan, L., Zhou, T.: Stein variational gradient descent with local
  approximations.
\newblock Computer Methods in Applied Mechanics and Engineering \textbf{386},
  114087 (2021)

\bibitem{yang2017bayesian}
Yang, K., Guha, N., Efendiev, Y., Mallick, B.K.: Bayesian and variational
  bayesian approaches for flows in heterogeneous random media.
\newblock Journal of Computational Physics \textbf{345}, 275--293 (2017)

\bibitem{zhang2018advances}
Zhang, C., B{\"u}tepage, J., Kjellstr{\"o}m, H., Mandt, S.: Advances in
  variational inference.
\newblock IEEE transactions on pattern analysis and machine intelligence
  \textbf{41}(8), 2008--2026 (2018)

\bibitem{zhdanov2002geophysical}
Zhdanov, M.S.: Geophysical inverse theory and regularization problems, vol.~36.
\newblock Elsevier (2002)

\bibitem{zhu2018bayesian}
Zhu, Y., Zabaras, N.: Bayesian deep convolutional encoder--decoder networks for
  surrogate modeling and uncertainty quantification.
\newblock Journal of Computational Physics \textbf{366}, 415--447 (2018)

\bibitem{zhu2019physics}
Zhu, Y., Zabaras, N., Koutsourelakis, P.S., Perdikaris, P.: Physics-constrained
  deep learning for high-dimensional surrogate modeling and uncertainty
  quantification without labeled data.
\newblock Journal of Computational Physics \textbf{394}, 56--81 (2019)

\end{thebibliography}

\end{document}